%% file: main.tex
\documentclass[11pt]{amsart}

\input{preamble}

\begin{document}

\begin{abstract}
To every right-angled Coxeter group one can canonically associate a hyperbolic metric space, called its \emph{maximized hyperbolic space}. This space is Morse recognizing and is obtained from the Davis complex by coning off its standard product regions. We characterize, purely in terms of the defining graph, when the Gromov boundary of a maximized hyperbolic space is nonempty and linearly connected. 

As an application, we prove that among one-ended right-angled Coxeter groups, splitting over a subgroup of an unbounded product subgroup is a quasi-isometry invariant. We also establish several features shared by maximized hyperbolic spaces and curve graphs, including the existence of dead ends and the visibility of boundaries. Moreover, when the boundary is nonempty, it is non-compact if and only if the defining graph contains an induced 4-cycle.
\end{abstract}

\maketitle

\setcounter{tocdepth}{1}
\tableofcontents


\section{Introduction} \label{intro section}

\input{sections/intro}

\section{Background}\label{Background Section}
\input{sections/background}

\section{Disconnectedness of the Gromov boundary} \label{Disconnectivity section}
\input{sections/disconnectedness}

\section{Morse boundaries and Connectedness of Gromov boundaries}
\label{connectivity of the boundary section}
\input{sections/Morse-to-Gromov}

\section[Relation between X and the maximized hyperbolic space]{Relation between \texorpdfstring{$X$}{X} and \texorpdfstring{$\widehat{X}$}{the maximized hyperbolic space}}
\label{Relation with X Section}

\input{sections/relation}

\section[Dead ends in the maximized space]{Dead ends in \texorpdfstring{$\widehat{X}$}{the maximized hyperbolic space}}
\label{Dead Ends section}
\input{sections/dead-ends}

\section{Bounded geodesic image and applications}\label{BGI section}

\input{sections/BGI}

\section{The fan-and-filter construction}\label{Fan and Filter Section}

\input{sections/fan-filter}

\section[Linear Connectedness of the Gromov boundary]{Linear Connectedness of \texorpdfstring{$\partial \widehat{X}$}{the Gromov boundary}}
\label{Main Proof Section}
\input{sections/main-proof}

\section{Application to group splittings}
\label{Application section}

\input{sections/group-splitting}


\bibliographystyle{alpha}
\bibliography{mybib}
\end{document}

%% file: preamble.tex
\usepackage[margin=1.25in]{geometry}
\usepackage{amsmath,amssymb,amsthm,amstext,mathrsfs}
\usepackage[shortlabels]{enumitem}
\usepackage{graphicx}
\usepackage{tikz-cd}
\usepackage{caption}
\usepackage{svg}

\usepackage[colorlinks=true,linkcolor=teal,citecolor=cyan,urlcolor=.,pagebackref=true]{hyperref}
\usepackage[nameinlink,nosort,capitalize]{cleveref} 

\crefname{enumi}{condition}{conditions}
\Crefname{enumi}{Condition}{Conditions}

\usepackage{verbatim}

\usepackage{import}
\usepackage{xifthen}
\usepackage{pdfpages}
\usepackage{transparent}

\usepackage{tabularx}
\usepackage{amscd}
\usepackage{xcolor}
\usepackage{scalerel,stackengine}

\def\diam{{\rm{diam}}}

\newcommand{\nn}{\mathcal N} 

\def\<{\left\langle}
\def\>{\right\rangle}

\def\CC{{\mathcal C}}
\def\CCS{{\mathcal C}S}

\def\FF{{\mathcal{F}}}

\def\lk{\text{lk}}
\def\st{\text{st}}
\def\into{\hookrightarrow}

\def\lab{\text{Label}}

\newtheorem{thm}{Theorem}[section]
\newtheorem{prop}[thm]{Proposition}
\newtheorem{cor}[thm]{Corollary}
\newtheorem{conj}[thm]{Conjecture}

\newtheorem{qn}[thm]{Question}
\newtheorem{lem}[thm]{Lemma}
\newtheorem*{claim*}{Claim}

\theoremstyle{definition} 
\newtheorem{defn}[thm]{Definition}
\newtheorem{exa}[thm]{Example}

\theoremstyle{remark}
\newtheorem{rmk}[thm]{Remark}

\newcommand{\tsh}[1]{\left\{\kern-.9ex\left\{#1\right\}\kern-.9ex\right\}}

\def\R{{\mathbb R}}

\def\Z{{\mathbb Z}}
\def\N{{\mathbb N}}

\def\CC{{\mathcal C}}

\def\A{{\mathcal A}}
\def\W{{\mathcal W}}
\def\B{{\mathcal B}}

\def\fg{{\mathfrak g}}

\def\fS{{\mathfrak S}}

\def\<{\langle}
\def\>{\rangle}

\definecolor{darkgreen}{cmyk}{1,0,1,.2}

\stackMath
\newcommand\reallywidehat[1]{%
\savestack{\tmpbox}{\stretchto{%
  \scaleto{%
    \scalerel*[\widthof{\ensuremath{#1}}]{\kern-.6pt\bigwedge\kern-.6pt}%
    {\rule[-\textheight/2]{1ex}{\textheight}}
  }{\textheight}%
}{0.5ex}}%
\stackon[1pt]{#1}{\tmpbox}%
}

\let\emptyset\varnothing

\title[Linear Connectivity of $\partial \widehat{X}$]{Linear Connectivity of Boundaries\\ 
of Right-Angled Coxeter Groups}
\author{Zhihao Mu}

\address{Department of Mathematics\\
The Graduate Center of City University of New York\\
365 Fifth Ave\\
New York, NY, USA, 10016}
\email{zmu@gc.cuny.edu}

\date{\today}

\makeatletter
\renewcommand{\@settitle}{%
  \begin{center}
    \normalfont\fontsize{16}{20}\selectfont\bfseries
    \@title\par
  \end{center}%
}
\makeatother

%% file: sections/intro.tex
\subsection{The main theorem}
We study the large-scale geometry of right-angled Coxeter groups. A right-angled Coxeter group (RACG) $W_\Gamma$ is determined by a finite simplicial graph $\Gamma=(V,E)$. Its generators are the vertices of $\Gamma$; each is an involution, and two generators commute precisely when the corresponding vertices are adjacent. A central theme is to determine which features of the coarse geometry of $W_\Gamma$ are encoded by the combinatorics of $\Gamma$. For example, $W_\Gamma$ is hyperbolic if and only if $\Gamma$ has no induced 4-cycle \cite{moussong1988hyperbolic}; it is one-ended if and only if $\Gamma$ is connected, not a clique, and has no separating clique \cite{mihalik2009visual,davis2025geometry}. Other graph-theoretic characterizations relate to divergence, thickness, relative hyperbolicity, the Morse boundary, and the structure of quasi-flats \cite{behrstock2012divergence,dani2015divergence,behrstock2017thickness,levcovitz2018divergence,levcovitz2022characterizing,cordes2025connectivity,huang2017top,behrstock2021quasiflats}. For a broader overview, we refer to the survey~\cite{dani2018large}.

We add a new entry to this dictionary. Our main result characterizes, purely in terms of $\Gamma$, when one has connectivity of the Gromov boundary of a particular hyperbolic space associated to $W_\Gamma$. The relevant hyperbolic space is called the \emph{maximized hyperbolic space}, which was introduced by Abbott--Behrstock--Durham in \cite{abbott2021largest}.

The geometry of $W_\Gamma$ is carried by the \emph{Davis complex} $X=X_\Gamma$, the flag cube complex on the Cayley graph of $W_\Gamma$ with the standard generating set $V(\Gamma)$. It is a $\mathrm{CAT}(0)$ cube complex on which $W_\Gamma$ acts properly and cocompactly. When $\Gamma$ contains an induced 4-cycle, $W_\Gamma$ contains a $\mathbb{Z}^2$ subgroup and $X$ contains flats, so $X$ is far from being hyperbolic. These flats lie inside \emph{standard product regions}: convex subcomplexes $gX_\Lambda$, where $g\in W_\Gamma$ and $\Lambda \subset\Gamma$ is a nontrivial join subgraph which has a decomposition $\Lambda = \Lambda_0\ast \Lambda_1$ in which every vertex in $\Lambda_0$ is adjacent to every vertex in $\Lambda_1$ and neither $\Lambda_i$ is a clique. For an RACG, the \emph{maximized hyperbolic space} is obtained by collapsing standard product regions: $\widehat{X} =\widehat{X}_{\Gamma}$ is the graph with vertex set $X^{(0)}$ in which two vertices are connected by an edge when they are adjacent in $X$ or lie in a common standard product region (Definition~\ref{Definition of maximized hyperbolic space}). Coning off these regions removes the non-hyperbolic directions and produces a $\delta$-hyperbolic space $\widehat{X}$~\cite[Corollary~2.16]{behrstock2021quasiflats}.  Moreover, $\widehat{X}$ has infinite diameter precisely when $\Gamma$ is neither a clique nor a nontrivial join; see Lemma~\ref{empty Gromov boundary}. If $\Gamma$ is a clique, $X_\Gamma$ is bounded; if $\Gamma$ is a nontrivial join, then $X_\Gamma$ is a single standard product region, hence $\widehat{X}$ has bounded diameter. Thus $\partial\widehat{X}$ is empty in both cases.

Our main theorem shows that the connectivity of $\partial\widehat{X}$ can be characterized from the defining graph $\Gamma$.

\begin{thm}
\label{main theorem}
Let $W_\Gamma$ be a right-angled Coxeter group with defining graph $\Gamma$, and let $\widehat{X}_\Gamma$ be its maximized hyperbolic space. The Gromov boundary $\partial\widehat{X}_{\Gamma}$ is nonempty and connected if and only if $\Gamma$ satisfies the following:
\begin{itemize}
  \item[\textup{(a)}] $\Gamma$ is connected, is not a clique, and has no separating clique; and
  \item[\textup{(b)}] $\Gamma$ has no separating nontrivial join; and
  \item[\textup{(c)}] for every vertex $v\in \Gamma$, the link of $v$ is not contained in a nontrivial join.
\end{itemize}
Furthermore, when $\partial\widehat{X}_\Gamma$ is connected, it is linearly connected with respect to a visual metric.
\end{thm}

Here, a metric space is \emph{linearly connected} if there is a constant $L\geq 1$ such that for any two points $x,y$, there is a compact connected set containing $x,y$ with diameter at most $Ld(x,y)$.  

The connectedness implication also follows from Cordes--Levcovitz's theorem on the connectedness of Morse boundaries of Coxeter groups \cite{cordes2025connectivity}, together with the density of the Morse-boundary image in $\partial\widehat{X}_\Gamma$. We prove this implication in Corollary~\ref{direct proof of connectedness}, within a general framework that also applies to other groups. Our direct geometric argument gives the quantitative conclusion of linear connectedness.

The three conditions have transparent group-theoretic and geometric meanings: (a) says that $W_\Gamma$ is one-ended; (b) says that $W_\Gamma$ does not split as an amalgam over a special unbounded product subgroup; and (c) says no hyperplane in the Davis complex $X_{\Gamma}$ is contained in a standard product region.

The geometry of the Davis complex $X$ is closely related to its standard product regions and the maximized hyperbolic space $\widehat{X}$. Given a pair of vertices in $X$, their distance in $X$ can be estimated by the sum of the distances of their images in standard product regions under the gate map (a closest point projection) and their images in the maximized hyperbolic space under the identity map. This property is captured by the distance formula, which is a result from \cite{behrstock2017hierarchically}; see Theorem~\ref{Distance Formula}. More generally, right-angled Coxeter groups are \emph{hierarchically hyperbolic groups} (Section~\ref{HHS subsection}), whose geometries are encoded, by  coarse Lipschitz projections, in a collection of hyperbolic spaces with a uniform hyperbolicity constant which are organized by three relations: \emph{nesting, orthogonality and transversality}. Hierarchically hyperbolic groups/spaces are a generalization of hyperbolic groups in the manner of Masur--Minsky's description of a mapping class group through the curve graphs of its essential subsurfaces \cite{behrstock2017hierarchically}. Each HHG/HHS may admit several hierarchically hyperbolic structures. For RACGs, the \emph{maximized structure} in the sense of Abbott--Behrstock--Durham \cite[Theorem~3.7]{abbott2021largest} has $\widehat{X}$ as its nested-maximal (top-level) hyperbolic space, and $W_\Gamma$ acts on $\widehat{X}$ as its largest and universal acylindrical action \cite[Theorem~A]{abbott2021largest}. Moreover, a quasi-isometry on $W_\Gamma$ induces a quasi-isometry on its maximized hyperbolic space $\widehat{X}$ \cite[Theorem 6.10]{goldsborough2023induced}, hence a homeomorphism of the Gromov boundary $\partial\widehat{X}$.
Consequently, the connectivity of $\partial\widehat{X}$ is a quasi-isometry invariant of $W_\Gamma$, and Theorem~\ref{main theorem} characterizes this invariant through the combinatorics of $\Gamma$.

\begin{figure}[h]
\centering

\def\graphheight{1.4}
\pgfmathsetmacro{\pentagonradius}{\graphheight/(1+cos(36))}
\pgfmathsetmacro{\hexagonradius}{\graphheight/sqrt(3)}
\pgfmathsetmacro{\horizontalstretch}{\hexagonradius/(\pentagonradius*cos(18))}
\pgfmathsetmacro{\verticalshift}{\graphheight/2-\pentagonradius}

\begin{tikzpicture}[baseline=0pt, scale=1.5,
    every node/.style={circle, draw, fill=black,
                       inner sep=0pt, minimum size=3.5pt}]

  \begin{scope}[xscale=\horizontalstretch, yshift=\verticalshift cm]
    \foreach \i in {0,...,4}{
      \node (a\i) at ({90-72*\i}:\pentagonradius) {}; 
      \node (b\i) at ({90-72*\i}:0.36)   {}; 
      \node (c\i) at ({54-72*\i}:0.48)   {}; 
    }

    \draw (0,0) circle[radius=0.48];

    \foreach \i/\j in {0/1,1/2,2/3,3/4,4/0}{
      \draw (a\i)--(a\j);             
      \draw (b\i)--(b\j);             
      \draw (a\i)--(b\i);             
      \draw (a\i)--(c\i)--(b\j);      
    }
  \end{scope}

  \node[draw=none, fill=none] at (0,-1.05) {$\Gamma_1$};
\end{tikzpicture}
\hspace{2cm}
\begin{tikzpicture}[baseline=0pt, scale=1.5,
    every node/.style={circle, draw, fill=black,
                       inner sep=0pt, minimum size=3.5pt}]
  \foreach \i in {0,1,2,3,4,5}{
    \node (h\i) at ({60*\i}:\hexagonradius) {};
  }
  \foreach \i/\j in {0/1,1/2,2/3,3/4,4/5,5/0}{
    \draw (h\i)--(h\j);
  }

  \node[label=right:$v$] (v) at (0,0) {};
  \draw (v)--(h3);
  \draw (v)--(h5);

  \node[draw=none, fill=none] at (0,-1.05) {$\Gamma_2$};
    \end{tikzpicture}
\caption{The maximized hyperbolic space $\widehat{X}_{\Gamma_1}$ of $W_{\Gamma_1}$ has connected Gromov boundary, but the Gromov boundary of $\widehat{X}_{\Gamma_2}$ is disconnected since the link of the vertex $v$ is contained in an induced 4-cycle in $\Gamma_2$.}
\label{Examples}

\end{figure}
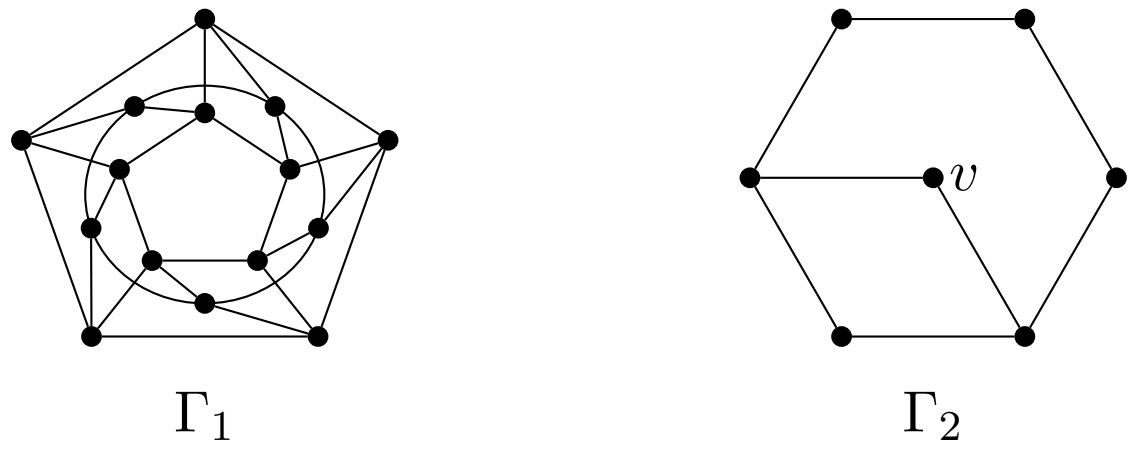

We record some applications of Theorem~\ref{main theorem}. First, we provide new examples of HHGs with connected Gromov boundary on their nested-maximal hyperbolic spaces; see Figure~\ref{Examples}. The classical examples whose top-level hyperbolic space has connected Gromov boundary include one-ended hyperbolic groups and mapping class groups of surfaces of complexity at least two. A recent work of Choi--Talbott also provides infinitely many examples from random quotients of HHGs \cite{choi2026random}.

Second, $\widehat{X}$ distinguishes certain right-angled Coxeter groups up to quasi-isometry. The quasi-isometry classification of RACGs, and the question of which RACGs are quasi-isometric or commensurable to RAAGs, have received considerable attention \cite{behrstock2017counterexample,dani2024right,cashen2025raagedy}. Because connectivity of $\partial\widehat{X}_{\Gamma}$ is a quasi-isometry invariant of $W_\Gamma$, Theorem~\ref{main theorem} immediately gives the following.
\begin{cor}
\label{not qi to racg}
    Let $W_\Gamma$ and $W_\Lambda$ be right-angled Coxeter groups. If $\Gamma$
    satisfies conditions \textup{(a)}--\textup{(c)} of
    Theorem~\ref{main theorem} but $\Lambda$ does not, then $W_\Gamma$ and
    $W_\Lambda$ are not quasi-isometric.
\end{cor}

Every RAAG $A_\Gamma$ is hierarchically hyperbolic. When $\Gamma$ is a connected graph with more than one vertex, $A_\Gamma$ admits universal \cite{kim2014geometry} and largest \cite{abbott2019hyperbolic} acylindrical actions on its extension graph. This graph is quasi-isometric to the contact graph which is a  quasi-tree \cite{hagen2014weak}, so its boundary is totally disconnected; see also Remark~\ref{RAAG boundary}.  We therefore obtain the following corollary.

\begin{cor}
\label{not qi to raag}
    If the defining graph $\Gamma$ satisfies conditions
    \textup{(a)}--\textup{(c)} of Theorem~\ref{main theorem}, then
    $W_{\Gamma}$ is not quasi-isometric to any right-angled Artin group.
\end{cor}
We remark that Corollaries~\ref{not qi to racg} and~\ref{not qi to raag} also follow from Cordes--Levcovitz's theorem on the connectedness of Morse boundaries \cite{cordes2025connectivity}.

Third, in the appendix ~\cite[Theorem A.1]{goldsborough2023induced}, Russell showed that if the nested-maximal hyperbolic space $\CCS$ in an HHG $G$ with a suitable HHS structure is one-ended, then the group of quasi-isometries of $\CCS$ is isomorphic to the group of quasi-isometries of $G$. Combining Russell's result with Corollary~\ref{one-endedness of maximized hyperbolic space}, we obtain the following application to right-angled Coxeter groups.

\begin{cor}
    Let $W_{\Gamma}$ be a right-angled Coxeter group whose defining graph $\Gamma$ satisfies conditions (a)--(c) in Theorem~\ref{main theorem} and $\widehat{X}_{\Gamma}$ be the maximized hyperbolic space of $W_{\Gamma}$. Then $\text{QI}(W_{\Gamma})$ is isomorphic to $\text{QI}(\widehat{X}_{\Gamma})$.
\end{cor}

Fourth, combining the visual decomposition theorem of Mihalik--Tschantz \cite{mihalik2009visual} with Theorem~\ref{splittings are topological} relates algebraic splittings of $W_{\Gamma}$ to the topology of $\partial\widehat{X}$.

\begin{cor}
\label{splitting iff disconnected}
A one-ended right-angled Coxeter group $W_\Gamma$ splits over a subgroup of an unbounded product subgroup if and only if the Gromov boundary $\partial\widehat{X}$ is empty or disconnected.
\end{cor}

\begin{cor}
    In the class of one-ended right-angled Coxeter groups, splitting over a subgroup of an unbounded product subgroup is a quasi-isometric invariant.
\end{cor}

Along the way, we establish several properties that $\widehat{X}$ shares with curve graphs. In Section~\ref{Dead Ends section} we exhibit two RACGs whose maximized hyperbolic spaces contain dead ends. In Section~\ref{BGI section}, using the bounded geodesic image theorem, we show that $\partial\widehat{X}$ is visible, i.e., any two points in $\widehat{X}\cup\partial\widehat{X}$ are joined by a geodesic (Lemma~\ref{Visible} and  Proposition~\ref{Visibility of boundary}),  $\widehat{X}$ is \emph{geodesic almost-extendable} (Proposition~\ref{exist close geodesic ray}), and when $\partial\widehat{X}$ is nonempty, $\partial\widehat{X}$ is non-compact if and only if $\Gamma$ contains an induced 4-cycle (Theorem~\ref{noncompact}).

\subsection{Related work}

\emph{Boundaries of curve graphs and related complexes.}
Theorem~\ref{main theorem} belongs to a broader program relating the topology of boundaries of hyperbolic spaces to the geometry of groups. For mapping class groups, after work of Leininger--Schleimer \cite{leininger2009connectivity} and Leininger--Mj--Schleimer \cite{leininger2011universal}, Gabai \cite{gabai2009almost} proved that the Gromov boundary of the curve graph of a surface of complexity at least two is path-connected and locally path-connected. Wright sharpened this to linear connectedness \cite{wright2024spheres}. In related settings, Long--Tan showed that the Gromov boundary of the fine curve graph is path connected when the genus of the surface is at least 3 and it is linearly connected under certain choices of the visual metric \cite{long2024connectedness}.
Bestvina--Chaika--Hensel  showed that the Gromov boundary of the free factor complex of a free group of rank at least 18 is path connected and locally path connected \cite{bestvina2023connectivity}.

\emph{Boundaries of Coxeter and CAT(0) groups.}
One may instead study the visual ($\mathrm{CAT}(0)$) boundary and the Morse boundary of Coxeter groups. Mihalik--Ruane--Tschantz \cite{mihalik2007local} proved that the visual boundary of a one-ended RACG is locally connected when there are no product separators (Definition~\ref{product separator definition}) and no virtual factor separators (Definition~\ref{virtual factor separator definition}). Cordes--Levcovitz \cite{cordes2025connectivity} proved that a one-ended RACG with \emph{wide-avoidant} defining graph has connected and locally connected Morse boundary. These two graph conditions coincide, and together they are equivalent to conditions (a)--(c) of Theorem~\ref{main theorem} (see Lemma~\ref{conditions are equivalent} for more details). The fan-and-filter construction used in Section~\ref{Fan and Filter Section}, an essential ingredient in our proof, builds on \cite{mihalik2007local,cordes2025connectivity}.

The Morse boundary $\partial_M W_\Gamma$ admits a natural continuous injection into the Gromov boundary $\partial\widehat{X}$ of the maximized hyperbolic space; see \cite[Theorem~6.6]{abbott2021largest} and \cite[Lemma~A.6]{russell2025extensions}. For non-virtually-cyclic HHGs containing a Morse element, Abbott--Behrstock--Russell established density of the Morse-boundary image in the HHS boundary \cite[Corollary~5.10 and Remark~5.11]{abbott2025structure}.

In Section~\ref{connectivity of the boundary section}, we record a general density statement for cobounded Morse-recognizing actions on hyperbolic spaces (Proposition~\ref{dense Morse boundary}). The resulting continuous boundary injection has a dense image, so connectedness of the Morse boundary implies connectedness of the Gromov boundary. We apply the result to the hyperbolic models for finite-rank Coxeter groups constructed in \cite{cashen2026hyperbolic}. Combined with \cite[Theorem~B]{cordes2025connectivity}, this gives an alternative proof of the connectedness implication of Theorem~\ref{main theorem}.

\emph{Splittings and boundary topology.}
The topology of a group boundary is classically tied to its algebraic splittings. By Stallings' theorem on the number of ends \cite{MR415622}, a finitely generated  group splits over a finite subgroup if and only if it has more than one end. Bowditch showed that a one-ended hyperbolic group, which is not virtually cocompact Fuchsian, splits over a two-ended subgroup if and only if its Gromov boundary has a local cut point \cite{bowditch1998cut}. These were extended to $\mathrm{CAT}(0)$ groups by Papasoglu--Swenson \cite{papasoglu2006continua, papasoglu2009boundaries}, to relatively hyperbolic groups by Groff \cite{groff2013quasi}, and for splittings over non-parabolic two-ended subgroups by Haulmark \cite{haulmark2019local} in terms of local cut points of the Bowditch boundary. The relation between splittings and connectivity of the visual boundary of $\mathrm{CAT}(0)$ groups was studied in \cite{mihalik1999cat, mihalik2001cat, hruska2021connectedness, mihalik2007local}. Moreover, Papasoglu showed that a finitely presented one-ended group which is not commensurable to a surface group splits over a two-ended subgroup if and only if its Cayley graph is separated by a quasi-line \cite{papasoglu2005quasi}.

\subsection{Outline of the proof}

For the forward direction,  first consider the bounded cases. If $\Gamma$ is a clique or a nontrivial join, then $\widehat{X}$ has bounded diameter and hence empty boundary. In the remaining cases, failure of one of conditions (a)--(c) gives a splitting over a finite subgroup or over a subgroup of an unbounded product subgroup. Bass--Serre theory then detects separation of $\widehat{X}$ by a bounded set, and the boundary is disconnected. See Section~\ref{Disconnectivity section} for more details.

In the backward direction, we show that under the three conditions (a)--(c) on the defining graph $\Gamma$, the Gromov boundary $\partial \widehat{X}$ is linearly connected (Definition~\ref{linear connected defn}), which implies path connectedness and local path connectedness. Fix a point $o$ in $\widehat{X}$ and take two points $a,b$ in the boundary $\partial\widehat{X}$. In Lemma~\ref{Visible}, we show that there are geodesic rays $\alpha,\beta$ based at $o$ converging to $a,b$, respectively. Then the key idea is to push paths along $\alpha \cup \beta$ away from $o$. To achieve this, we iteratively use Proposition~\ref{spheres are connected}, which states that any two points $x$ and $y$ on a sphere of radius $R$ centered at $o$ in $\widehat{X}$ can be joined by a path lying in the union of the spheres of radii $R$ and $R+1$, with the diameter bounded from above by a uniform constant multiple of $d_{\widehat{X}}(x,y)$.

The proof of Proposition~\ref{spheres are connected} uses the cubical structure of the Davis complex $X$. Identify the vertex sets of $X$ and $\widehat{X}$, and choose geodesic segments $\gamma_x$ and $\gamma_y$ from $o$ to $x$ and $y$. We construct a \emph{truncated filter}: a van Kampen diagram over $W_\Gamma$ whose boundary maps to the union of $\gamma_x$, $\gamma_y$, and a path from $x$ to $y$ in $X^{(1)}$ (Figure~\ref{Filter}). The image of the third path in $\widehat{X}$ lies in the union of the spheres of radii $R$ and $R+1$ centered at $o$.

Mihalik--Ruane--Tschantz introduced filters between geodesic rays to study the visual boundary of RACGs \cite{mihalik2007local}; Cordes--Levcovitz later generalized the construction to study Morse boundaries of Coxeter groups \cite{cordes2025connectivity}. Here we instead construct truncated filters between geodesic segments. A (truncated) filter contains a directed tree whose directed paths map to geodesics in $X^{(1)}$. Furthermore, paths in the directed tree after their last intersection with the side paths have uniformly bounded wide subpaths, with a bound depending only on $\Gamma$. Lemma~\ref{C'} converts this estimate into progress in $\widehat X$. Consequently, the truncation requires only finitely many fan attachments.

\subsection{Open questions}
Mapping class groups are quasi-isometrically rigid \cite{behrstock2012geometry,hamenstaedt2005geometry} (for a direct approach using hierarchical hyperbolicity see \cite{behrstock2021quasiflats, mangioni2025randomquotientsmappingclass}). The proof factors through two steps: a quasi-isometry of the mapping class group induces an automorphism of the curve graph, and, by Ivanov's theorem \cite{ivanov1997automorphisms}, every such automorphism is realized by a mapping class. A motivation for our work is to find hierarchically hyperbolic groups whose nested-maximal hyperbolic space might play the role of the curve graph in such an argument. The contact graph of the standard hierarchical structure of RACGs is unsuitable: as a quasi-tree it can have an enormous automorphism group \cite[Section~7]{behrstock2017hierarchically}. The maximized hyperbolic space $\widehat{X}$, when $\partial\widehat{X}$ is connected, is a more promising candidate, and we expect a quasi-isometric model of it to have a countable automorphism group and to be useful to study the quasi-isometric rigidity problem. This motivates the following question.

\begin{qn}
Is $W_\Gamma$ quasi-isometrically rigid when $\Gamma$ satisfies conditions \textup{(a)}, \textup{(b)}, and \textup{(c)}?
\end{qn}

In view of Corollary~\ref{splitting iff disconnected}, we expect that the connectedness of the boundary can be characterized by the group splittings for general hierarchically hyperbolic groups.

\begin{conj}
Let $G$ be a one-ended hierarchically hyperbolic group with maximized hyperbolic space $\widehat{X}$. Suppose $\widehat{X}$ has infinite diameter. Then $G$ splits over a subgroup of an unbounded product subgroup if and only if $\partial\widehat{X}$ is disconnected.
\end{conj}

Bonk--Kleiner \cite{bonk2005quasi} showed that every one-ended hyperbolic group contains a quasi-isometrically embedded hyperbolic plane, and Mackay--Sisto \cite{Mackay_Sisto_2020} extended this to relatively hyperbolic groups with nilpotent peripheral subgroups. We ask an analogous question for HHGs.

\begin{qn}
Let $G$ be a one-ended hierarchically hyperbolic group with maximized hyperbolic space $\widehat{X}$. Suppose the Gromov boundary of $\widehat{X}$ is  connected. Does $\widehat{X}$ contain a quasi-isometrically embedded hyperbolic plane?
\end{qn}

A quasi-isometrically embedded hyperbolic plane in $\widehat{X}$ induces a quasisymmetrically embedded circle in $\partial\widehat{X}$. In light of the characterization of quasisymmetric images of circles and arcs by Tukia--V\"ais\"al\"a \cite{tukia1980quasisymmetric}, the following question is a natural step.

\begin{qn}
When $\partial \widehat{X}$ is connected, does it contain a quasiarc?
\end{qn}

\subsection{Organization}

Section~\ref{Background Section} reviews RACGs, equivalent graph-theoretic conditions, $\mathrm{CAT}(0)$ cube complexes, Gromov boundaries, and the maximized HHS structure. Section~\ref{Disconnectivity section} proves the forward direction of Theorem~\ref{main theorem}. 
Section~\ref{connectivity of the boundary section} records the density statement for Morse-recognizing actions and its connectedness applications. 
Section~\ref{Relation with X Section} develops the geometry of $\widehat{X}_\Gamma$, beginning with taut geodesics, and Section~\ref{Dead Ends section} gives two examples of dead ends. 
Section~\ref{BGI section} proves bounded geodesic image and establishes visibility and noncompactness of $\partial\widehat{X}$. 
Section~\ref{Fan and Filter Section} develops filters in $X$ and uses them to establish controlled connectivity of spheres in $\widehat{X}$. 
Section~\ref{Main Proof Section} establishes linear connectedness, completing the proof of the main theorem. 
Finally, Section~\ref{Application section} applies the main theorem to group splittings.

\subsection{Acknowledgements}
I am grateful to my advisor, Jason Behrstock, for introducing me to hierarchically hyperbolic spaces and for his continued support and encouragement throughout this project. 
I thank Jacob Russell for encouraging me to think about this problem at the Advances in Hierarchical Hyperbolicity workshop, and for his continued interest in this project. His question about the relationship between the Morse boundary and the Gromov boundary, raised during the Buckeye Metric Geometry Workshop II, led to the material in Section~\ref{connectivity of the boundary section}.
I also thank Michael Mihalik for his interest in this work and, in particular, for discussions of Lemma~\ref{conditions are equivalent} and Remark~\ref{VFS not equivalent to (c)} during the Vanderbilt Colloquium Series on Groups, Geometry, and Dynamics.

I was partially supported by a CUNY Graduate Center Dissertation Fellowship during this project.

\subsection{AI disclosure}
AI tools were used to proofread the manuscript.

%% file: sections/background.tex
We review the background on right-angled Coxeter groups, $\mathrm{CAT}(0)$ cube complexes, Gromov boundaries, and hierarchically hyperbolic structures needed later in the paper.

\subsection{Right-angled Coxeter groups}

Let $\Gamma$ be a simplicial graph with vertex set $V(\Gamma)$ and edge set $E(\Gamma)$. The \emph{right-angled Coxeter group} defined by $\Gamma$ is
\[
W_{\Gamma}=\left\langle V(\Gamma) \mid v^2=1\text{ for }v\in V(\Gamma),(vw)^2=1\text{ for }(v,w)\in E(\Gamma) \right\rangle.
\]

We call $\Gamma$ the \emph{defining graph} of $W_{\Gamma}$ and write $\Sigma_\Gamma$ for the Cayley graph with respect to $V(\Gamma)$. A subgraph $\Delta\subset\Gamma$ is \emph{induced} if it contains every edge of $\Gamma$ whose endpoints lie in $\Delta$. Unless stated otherwise, every subgraph we consider in this article is induced.

The following lemma is usually referred to as the \textit{Deletion Condition}. We refer the reader to \cite[Property 1.4.7]{bjorner2005combinatorics} for a proof.

\begin{lem}[Deletion Condition]
\label{Deletion Condition}
If a word $w=w_1\cdots w_n$ is non-geodesic in $\Sigma_\Gamma$, then there are $0<i<j\leq n$ such that
\[
w = w_1\cdots w_{i-1}w_{i+1} \cdots w_{j-1}w_{j+1}\cdots w_n
\]
in $W_\Gamma$.    
\end{lem}

A subgroup of $W_\Gamma$ is \emph{special} if it is generated by a subset of $V(\Gamma)$. The following standard result justifies the terminology. 

\begin{thm}[\protect{\cite[Proposition~2.4.1]{bjorner2005combinatorics}}]
    Let $W_{\Gamma}$ be a right-angled Coxeter group with defining graph $\Gamma$. Suppose $\Delta \subset \Gamma$ is an induced subgraph. Then the special subgroup generated by $V(\Delta)$ is isomorphic to the right-angled Coxeter group $W_{\Delta}$.
\end{thm}

A subgraph $\Delta$ \emph{separates} $\Gamma$ if $\Gamma\setminus\Delta$ is disconnected. A graph $\Lambda$ is a \emph{join}, written as $\Lambda = \Lambda_0* \Lambda_1$, if $\Lambda_0$ and $\Lambda_1$ are nonempty disjoint subgraphs whose union is $\Lambda$, and every vertex of $\Lambda_0$ is adjacent to every vertex of $\Lambda_1$. A \emph{clique} is a complete graph.

\begin{defn}[Nontrivial Join]
    A join $\Lambda = \Lambda_0* \Lambda_1$ is \emph{nontrivial} if neither $\Lambda_0$ nor $\Lambda_1$ is a clique.
\end{defn}

\begin{rmk}
     The right-angled Coxeter group $W_\Lambda$ whose defining graph is a nontrivial join decomposes as a direct product $W_\Lambda = W_{\Lambda_0} \times W_{\Lambda_1}$, where each factor is infinite. Moreover, the Morse boundary of $W_\Lambda$ is empty and no geodesic in $W_\Lambda$ is contracting.
\end{rmk}

\begin{exa}
    Let $\Gamma$ be a 4-cycle. Its two pairs of opposite vertices give a nontrivial join decomposition. Consequently,
    $W_{\Gamma}\cong D_{\infty}\times D_{\infty}$, where $D_{\infty}$ is the
    infinite dihedral group, and $W_\Gamma$ contains a finite-index copy of
    $\Z^2$.
\end{exa}

We call a subgroup $H$ of $G$ an \textit{unbounded product subgroup} if it decomposes as a direct product $H = H_0\times H_1$ and both $H_0,H_1$ are infinite. If $\Lambda \subset \Gamma$ is a nontrivial join subgraph, the special subgroup $W_\Lambda$ is an unbounded product subgroup of $W_\Gamma$.

Many algebraic and geometric properties of RACGs can be read from their defining graphs. Moussong's theorem classifies hyperbolic right-angled Coxeter groups.

\begin{thm}[\cite{moussong1988hyperbolic}]
    A right-angled Coxeter group $W_{\Gamma}$ is hyperbolic if and only if $\Gamma$ has no induced square.
\end{thm}

The number of ends of $W_{\Gamma}$ is  detectable in $\Gamma$, by Stallings' theorem \cite{MR415622} together with work of Mihalik--Tschantz \cite{mihalik2009visual} or Davis \cite{davis2025geometry}.

\begin{thm}[{\cite[Corollary 16]{mihalik2009visual}}]
\label{one ended}
$W_{\Gamma}$ is one-ended if and only if $\Gamma$ is connected, is not a clique, and has no separating clique.
\end{thm}

We conclude this subsection with a graph-theoretic lemma that will be used repeatedly. For a subgraph $\Lambda\subset\Gamma$, its \emph{link} $\lk(\Lambda)$ is the set of vertices adjacent to every vertex of $\Lambda$, and its \emph{star} is $\st(\Lambda)=V(\Lambda)\cup\lk(\Lambda)$.

\begin{lem}
\label{in the same join}
Let $\Lambda\subset\Gamma$ be a nontrivial join and let $\Sigma\subset\Gamma$. If $\Delta=\lk(\Sigma)\cap\Lambda$ is not a clique, then $\Sigma\cup\Delta$ is contained in a nontrivial join of $\Gamma$.
\end{lem}
\begin{proof}
Write $\Lambda=\Lambda_0*\Lambda_1$, where neither factor is a clique, and put $\Delta_i=\Delta\cap\Lambda_i$. At least one $\Delta_i$ is not a clique; otherwise $\Delta=\Delta_0*\Delta_1$ would be a clique. After interchanging the factors, assume that $\Delta_0$ is not a clique. Because $\Delta\subset\lk(\Sigma)$, every vertex of $\Delta_0$ is adjacent to every vertex of $\Sigma\cup\Delta_1$. Hence
\[
\Delta_0*(\Sigma\cup\Delta_1)
\]
is a join containing $\Sigma\cup\Delta$. If its second factor is not a clique, this join is nontrivial. Otherwise choose a nonclique subgraph
$\Lambda_1'\subset\Lambda_1$ containing $\Delta_1$. Then
\[
\Delta_0*(\Sigma\cup\Lambda_1')
\]
is a nontrivial join containing $\Sigma\cup\Delta$.
\end{proof}

\subsection{Equivalent graph-theoretic conditions}
\label{equivalent graph-theoretic conditions}

We compare the graph-theoretic hypotheses of our main theorem with those used by Mihalik--Ruane--Tschantz and Cordes--Levcovitz. For one-ended right-angled Coxeter groups, we show that conditions (b) and (c) are equivalent both  to the absence of product separators and virtual factor separators and to wide-avoidance. These equivalences will be used in Section~\ref{connectivity of the boundary section} to apply results on Morse boundaries, and in Section~\ref{Fan and Filter Section} to apply the fan-and-filter construction.

\begin{defn}[\protect{\cite[Definition 5.1]{mihalik2007local}}] 
\label{product separator definition}
A \emph{product separator} for $\Gamma$ is an induced subgraph $\Lambda$ such that $\Lambda$ separates $\Gamma$ and the special subgroup $W_\Lambda$ decomposes as $W_\Lambda = W_{\Lambda_0} \times W_{\Lambda_1}$, where both $W_{\Lambda_0}$ and $W_{\Lambda_1}$ are infinite.
\end{defn}

\begin{defn}[\protect{\cite[Definition 3.1]{mihalik2007local}}]
\label{virtual factor separator definition}
A \emph{virtual factor separator} for $\Gamma$ is a triple $(C,C_1,K)$ of induced subgraphs such that
\begin{enumerate}
    \item $C_1\subset C$ and $C$ separates $\Gamma$;
    \item $W_{C_1}$ has finite index in $W_C$;
    \item $K\subset \lk(C_1)$ and $W_K$ is infinite.
\end{enumerate}
\end{defn}

We record a theorem about the visual boundary of $W_\Gamma$ by Mihalik--Ruane--Tschantz.

\begin{thm}[\cite{mihalik2007local}]
    Let $W_\Gamma$ be a one-ended right-angled Coxeter group. If the defining graph $\Gamma$ has no product separator and no virtual factor separator, then the visual boundary of $W_{\Gamma}$ is locally connected.
\end{thm}

We note that, in the right-angled case, the notion of a wide subgraph used in \cite[Definition~3.1 and Remark~3.4.1]{cordes2025connectivity} agrees with that of an induced nontrivial join.

\begin{defn}[\protect{\cite[Definition 3.3]{cordes2025connectivity}}]  We say the graph $\Gamma$ is \emph{wide-avoidant} if given any nontrivial join subgraph $\Delta \subset \Gamma$ and any pair of vertices $s,t\in \Gamma$, there exists a path $\alpha \subset \Gamma$ from $s$ to $t$ such that $\alpha \cap \Delta \subset \{s,t\}$.
\end{defn}

We record a theorem about the Morse boundary of $W_\Gamma$ by Cordes--Levcovitz.

\begin{thm}[{\cite[Theorem B]{cordes2025connectivity}}]
Let $W_{\Gamma}$ be a one-ended right-angled Coxeter group. Then the Morse boundary $\partial_M W_\Gamma$ is connected and locally connected if and only if $\Gamma$ is wide-avoidant.
\end{thm}

The following observation is due to Michael Mihalik.

\begin{lem}
\label{conditions are equivalent}
    Let $\Gamma$ be the defining graph of a one-ended RACG $W_\Gamma$. Then the following are equivalent:
    \begin{itemize}
        \item[(i)] $\Gamma$ has neither product separator nor virtual factor separator;
        \item[(ii)] $\Gamma$ is wide-avoidant;
        \item[(iii)] $\Gamma$ has no separating nontrivial join and for each vertex $v\in \Gamma$, the link of $v$ is not contained in a nontrivial join. 
    \end{itemize}
\end{lem}
\begin{proof}
For \textup{(i)}$\Leftrightarrow$\textup{(iii)}, observe first that a special subgroup is infinite exactly when its defining subgraph is not a clique. Thus product separators are precisely separating nontrivial joins. When no product separator is present,
\cite[Lemma~5.3]{mihalik2007local} identifies virtual factor separators with vertices whose links are contained in nontrivial joins. This proves the equivalence.

We prove \textup{(ii)}$\Leftrightarrow$\textup{(iii)} directly. A separating nontrivial join $\Delta$ prevents vertices in different components of $\Gamma\setminus\Delta$ from being joined while avoiding $\Delta$. Likewise, if $\lk(v)\subset\Delta$ for a nontrivial join $\Delta$, then every path from $v$ to a vertex outside $\Delta\cup\{v\}$ meets $\Delta$ immediately. Hence failure of either condition in \textup{(iii)} implies failure of wide-avoidance.

Conversely, assume \textup{(iii)} and let $\Delta$ be a nontrivial join. Then $\Gamma\setminus\Delta$ is connected. Given $s,t\in V(\Gamma)$, let $s'=s$ if $s\notin\Delta$. If $s\in\Delta$, condition~\textup{(iii)} gives a vertex $s'\in\lk(s)\setminus\Delta$. Define $t'$ similarly. Join $s'$ to $t'$ by a path in $\Gamma\setminus\Delta$ and, when necessary, append the edges to $s$ and $t$, respectively. The resulting path meets $\Delta$ only at those endpoints that already lie in $\Delta$. Thus $\Gamma$ is wide-avoidant.
\end{proof}

\begin{rmk}
\label{VFS not equivalent to (c)}

\begin{figure}[t]
  \centering
  \includegraphics[width=0.7\linewidth]{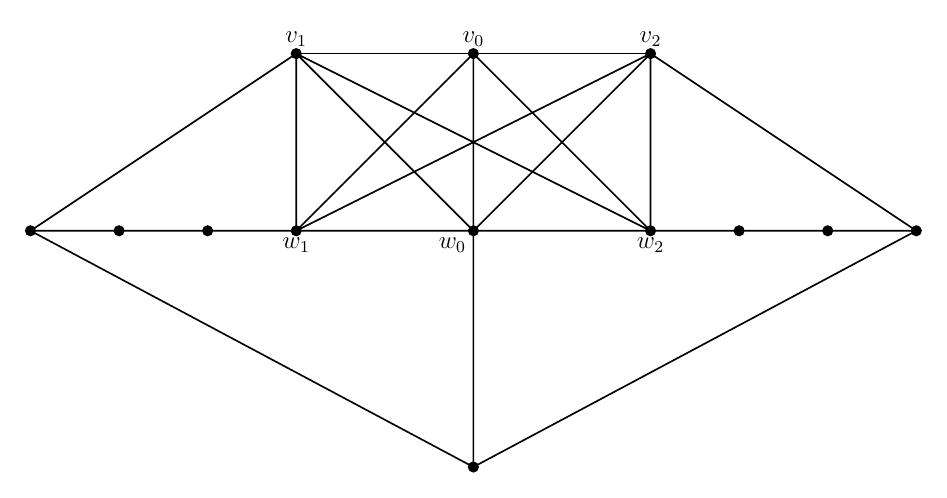}
  \caption{No virtual factor separator but the link of $v_0$ is contained in a nontrivial join.}
  \label{Link but no VFS}
\end{figure}

\begin{figure}[t]
  \centering
  \includegraphics[width=0.5\linewidth]{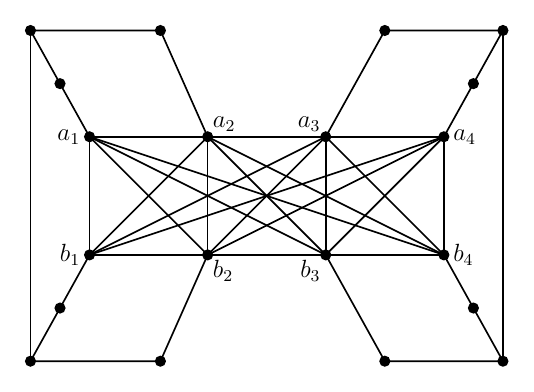}
  \caption{No link of a vertex is contained in a nontrivial join but $\{a_1,a_2,a_3,a_4,b_1,b_2\}$ spans a virtual factor separator.}
  \label{VFS but no Link}
\end{figure}

    The hypothesis excluding product separators is essential in the equivalence above. Without it, the absence of virtual factor separators is not equivalent to \cref{Condition:(c)}. See examples observed by Mihalik in Figures~\ref{Link but no VFS} and~\ref{VFS but no Link}.
\end{rmk}

\subsection[CAT(0) cube complexes]{$\mathrm{CAT}(0)$ cube complexes}
\label{Cube complex subsection}
The geometry of a Coxeter group is encoded by its \emph{Davis complex}; see \cite{davis2025geometry} for the construction and its properties. We write $X=X_\Gamma$ for the Davis complex of $W_\Gamma$. If $W_\Delta$ is a special subgroup of $W_\Gamma$, then $X_\Delta$ embeds isometrically as a convex subcomplex of $X_\Gamma$. Throughout the paper, Davis complexes carry their combinatorial metrics rather than their $\mathrm{CAT}(0)$ metrics.

For right-angled Coxeter groups, the 1-skeleton $X^{(1)}$ is precisely the Cayley graph $\Sigma_\Gamma$. After choosing a basepoint $o\in X$, we label its vertices by elements of $W_\Gamma$ and its edges by generators in $V(\Gamma)$. If $\alpha$ is an edge path in $X$, let $\overline{\alpha}=s_1\cdots s_n$ be the word obtained by reading its edge labels in order, and set
\[
\lab(\alpha)=\lab(\overline{\alpha})
=\{s\in V(\Gamma):s\text{ occurs in }\overline{\alpha}\}.
\]

All edges dual to a given hyperplane of $X$ have the same label $v\in V(\Gamma)$; we say that such a hyperplane is of \emph{type $v$}. Let $e_v$ be the edge labeled $v$ incident to $o$, and let $H_v$ be its dual hyperplane. Then
\[
\operatorname{Stab}_{W_\Gamma}(H_v)=\langle\st(v)\rangle.
\]
More generally, a hyperplane $gH_v$ of type $v$ has stabilizer $g\langle\st(v)\rangle g^{-1}$.

Each hyperplane $H$ separates $X$ into two halfspaces, denoted by $\overleftarrow{H}$ and $\overrightarrow{H}$. Its \emph{carrier} $\nn(H)$ is the union of the closed cubes that meet $H$. There is a cubical isometric embedding $H\times[-\frac12,\frac12]\hookrightarrow X$; the images of $H\times\{\pm\frac12\}$ are the two \emph{combinatorial hyperplanes}, denoted by $H^\pm$.

For 0-cubes $x,y\in X$, the distance $d(x,y)$ equals the number of hyperplanes separating them. The \emph{convex hull} of a subcomplex $Y\subset X$ is the intersection of all convex subcomplexes containing $Y$. By Sageev's construction \cite{sageev1997codimension}, the convex hull of two 0-cubes is cubically isometric to the cube complex dual to the wallspace of hyperplanes separating them.

A graph is \emph{median} if every triple of vertices $a,b,c$ has a unique vertex $m(a,b,c)$ lying on a geodesic between each pair. The 1-skeleton of a $\mathrm{CAT}(0)$ cube complex is a median graph \cite{chepoi2000graphs}. For 0-cubes $x,y\in X$, the median interval between them is the 1-skeleton of their convex hull; we denote it by $[x,y]$.

If $Y\subset X$ is convex, the closest-point projection $\fg_Y:X\to Y$ is well defined and we call it the \emph{gate map}. Thus $\fg_Y(x)$ is the unique 0-cube $y\in Y$ minimizing the distance $d(x,z)$ over all 0-cubes $z\in Y$. Let $\mathcal{W}(Y)$ denote the set of hyperplanes crossing $Y$.

The following is a standard application of the gate map, which is usually referred to as the \emph{bridge lemma} introduced in~\cite{behrstock2012divergence}. See~\cite[Lemma 3.23]{fernos2024homotopy} for a complete proof. 
\begin{lem}[Bridge Lemma]
\label{Bridge Lemma}
    Let $P,Q$ be convex subcomplexes of a CAT(0) cube complex $X$. Let $\fg_P:X\to P$ and $\fg_Q:X\to Q$ be the gate maps. Then the following hold:
    \begin{enumerate}
        \item $\fg_P(Q)$ and $\fg_Q(P)$ are convex subcomplexes.
        \item $\mathcal{W}(\fg_P(Q)) = \mathcal{W}(\fg_Q(P)) = \mathcal{W}(P)\cap \mathcal{W}(Q)$.
        \item The map $\fg_P: \fg_Q(P)\to \fg_P(Q)$ is an isomorphism of CAT(0) cube complexes.
        \item The cubical convex hull $[\fg_P(Q),\fg_Q(P)]$ is a CAT(0) cube complex isomorphic to $\fg_P(Q)\times [x,\fg_P(Q)]$ for a vertex $x\in \fg_Q(P)$.
        \item If $P\cap Q \neq \emptyset$, then $\fg_P(Q) = \fg_Q(P) = P\cap Q$ and $\mathcal{W}(P\cap Q) = \mathcal{W}(P)\cap \mathcal{W}(Q)$.  
    \end{enumerate}
\end{lem}

The next lemma is standard in Coxeter group theory; see
\cite{humphreys1992reflection}. We include a cubical proof. It says that the
intersection of the median interval $[o,x]$ with the 1-neighborhood of $x$ spans a cube.

\begin{lem}
\label{B(g) is a clique lemma}
    Let $\alpha$ be a geodesic from $o$ to $x$ in $X$. Suppose $e$ and $f$ are 1-cubes containing $x$ such that $\alpha e$ and $\alpha f$ are not geodesic. Then $e$ and $f$ span a 2-cube.
\end{lem}
\begin{proof}
     Let $H_e,H_f$ be the hyperplanes dual to $e,f$, respectively. Suppose $H_e$ and $H_f$ do not cross. $H_e$ separates $X$ into two halfspaces $\overleftarrow{H_e}$ and $\overrightarrow{H_e}$. Suppose $x\in \overrightarrow{H_e}$. Since $H_f$ does not cross $H_e$, $H_f\subset \overrightarrow{H_e}$. Since $\alpha e$ is not geodesic, $H_e$ separates $o$ and $x$. In particular, $o \in \overleftarrow{H_e}$. Therefore, $H_f$ does not separate $o$ and $x$, which implies $\alpha f$ is geodesic. Contradiction. Hence $H_e$ must cross $H_f$ and the 1-cubes $e,f$ span a 2-cube.
\end{proof}

\begin{rmk}
\label{B(g) is a clique remark}

Let $\alpha$ be a geodesic from $x$ to $y$ in $X$ and $\overline{\alpha}$ be the word in $W_\Gamma$ corresponding to $\alpha$. Denote by $B(\alpha)  = B(\overline{\alpha})$ the set of vertices $v\in \Gamma$ such that $\overline{\alpha}v$ is non-geodesic. By Lemma~\ref{B(g) is a clique lemma}, $B(\alpha)$ spans a clique in $\Gamma$.
\end{rmk}

Let $J\subset\Gamma$ be a nontrivial join. It is \emph{maximal} if it is not properly contained in another nontrivial join. For a maximal nontrivial join $J$ and $g\in W_\Gamma$, we call $gX_J\subset X_\Gamma$ a \emph{standard product region}.

\subsection{Gromov-hyperbolic spaces and their boundaries}

Let $(X,d)$ be a metric space and fix $o\in X$. The \emph{Gromov product} based at $o$ is
\[
\langle x,y\rangle_o  =\frac{1}{2}\left( d(o,x) + d(o,y) - d(x,y) \right)
\]
for any $x,y\in X$.

A metric space is \emph{Gromov hyperbolic}, or \emph{$\delta$-hyperbolic}, if there is $\delta\geqslant0$ such that
\[
\langle x,y\rangle_o \geq \min\{\langle x,z \rangle_o,\langle y,z\rangle_o\} -\delta
\]
for all $x,y,z,o\in X$.

When $X$ is geodesic, this is equivalent, up to changing the hyperbolicity constant, to requiring every geodesic triangle to be thin. Moreover, $\langle x,y\rangle_o$ coarsely measures the distance from $o$ to a geodesic from $x$ to $y$, with additive error depending only on $\delta$.

A sequence $(x_n)_{n\geq0}$ in $X$ \emph{converges at infinity} if
\[
\liminf_{n,m\to\infty}\langle x_n,x_m\rangle_o=\infty.
\]

This condition is independent of the basepoint. Two sequences $(x_n)$ and $(y_n)$ converging at infinity are \emph{equivalent} if
\[
\liminf_{n,m\to \infty}{\langle x_n,y_m\rangle_o} = \infty.
\]
This is an equivalence relation; write $[(x_n)]$ for the class of $(x_n)$.

\begin{defn}[Gromov boundary]
    Let $(X,d)$ be a Gromov-hyperbolic space. Its \emph{Gromov boundary} is
    \[
    \partial X := \{ [(x_n)]\mid (x_n) \text{ is a sequence converging at infinity in $X$}\}.
    \]
\end{defn}

If $X$ is proper, meaning that every closed ball is compact, then every boundary point is represented by a geodesic ray \cite[III.H, Lemma~3.13]{bridson2013metric}. For a nonproper hyperbolic space, every boundary point is still represented by a quasi-geodesic ray; see, for example, \cite[Proposition~2.3]{long2024connectedness}.

The Gromov product can be extended to the boundary points. For $x\in X$ and $\zeta\in \partial X$,
\[
\langle x,\zeta \rangle_o = \inf \{ \liminf_{n\to \infty } \langle x,x_n\rangle_o : (x_n)\in \zeta \}. 
\]
Also for $\zeta, \eta\in \partial X$, we have
\[
\langle \zeta,\eta \rangle_o = \inf\{ \liminf_{n,m\to \infty} \langle x_n,y_m\rangle_o : (x_n)\in \zeta, (y_m) \in \eta\}.
\]

The sets
\[
U(\zeta,r) = \{\eta \in \partial X \mid \langle \zeta,\eta\rangle_o > r\}
\]
for $r>0$ form a neighborhood basis at $\zeta$ and define the topology on $\partial X$.

The boundary admits a class of mutually quasisymmetric visual metrics; see \cite{vaisala2005gromov,das2017geometry} for definitions in this generality. These metrics relate boundary distance to the Gromov product as follows.

\begin{prop}[{\cite[Proposition~5.16]{vaisala2005gromov},
\cite[Proposition~3.6.8]{das2017geometry}}]
\label{visual metric property}
If $b>1$ is sufficiently close to 1, then for each $o\in X$, a visual metric $d_{\infty}$ on $\partial X$ satisfies the following inequality:
\[
b^{-\langle \zeta,\eta\rangle _o} / 4 \leq d_{\infty}(\zeta,\eta) \leq b^{-\langle \zeta,\eta\rangle_o}
\]
for any $\zeta,\eta\in \partial X$. Furthermore, this metric is compatible with the topology on $\partial X$.
\end{prop}

A $(K,C)$-quasi-isometry of a $\delta$-hyperbolic space extends to a quasisymmetry of its boundary, with distortion controlled by $K$, $C$, and $\delta$ \cite[Theorem 5.35]{vaisala2005gromov}.

\subsection{Linear connectedness}
We recall the definition of linear connectedness.
\begin{defn}[Linear connectedness]
\label{linear connected defn}
    A metric space $(X,d)$ is \emph{$L$-linearly connected}, for $L\geq1$, if
    every $x,y\in X$ lie in a compact connected set $Y\subset X$ with
    $\diam(Y)\leq Ld(x,y)$.  A metric space is \emph{linearly connected} if it is
    $L$-linearly connected for some $L$.
\end{defn}

After increasing $L$ by an arbitrarily small amount, one may require $Y$ to be an arc; see the introduction of \cite{mackay2008existence}. A linearly connected metric space is path-connected and locally path-connected, and linear connectedness is preserved by quasisymmetries.

\subsection{A hierarchically hyperbolic structure on RACGs}
\label{HHS subsection}

Behrstock--Hagen--Sisto introduced hierarchically hyperbolic spaces and groups (HHSs and HHGs) in \cite{behrstock2017hierarchically}; we refer the reader there for the full definition. Special groups in the sense of Haglund--Wise, including right-angled Artin and Coxeter groups, are hierarchically hyperbolic. An HHG may admit several HHS structures. Here we use the maximized structure on RACGs introduced in \cite{abbott2021largest}.

Let $X=X_\Gamma$ be the Davis complex. Recall that an induced nontrivial join $J\subset\Gamma$ determines an unbounded product subgroup $W_J$ and an isometrically embedded subcomplex $X_J\subset X_\Gamma$. When $J$ is maximal, each translate $gX_J$ is a standard product region.

We now define the maximized hyperbolic space $\widehat{X}$, the nested-maximal domain in this HHS structure. We give only the RACG-specific definition; the general construction appears in \cite[Section~3.3]{abbott2021largest}.

\begin{defn}[Maximized hyperbolic space]
\label{Definition of maximized hyperbolic space}
Let $X=X_{\Gamma}$ be the Davis complex of $W_\Gamma$. The \emph{maximized hyperbolic space} $\widehat{X}=\widehat{X}_\Gamma$ is the graph with vertex set $X^{(0)}$ in which distinct vertices $x,y$ are adjacent if they are joined by an edge of $X$ or lie in a common standard product region.
\end{defn}

Give every edge of $\widehat{X}$ length one and denote the resulting path metric by $d_{\widehat{X}}$. We identify the vertex sets of $X$ and $\widehat{X}$ throughout.

\begin{rmk}
    Two vertices $x,y\in\widehat{X}$ are adjacent if and only if either they are joined by an edge in $X$ or the set $\lab([x,y])$ of labels of separating hyperplanes is contained in a nontrivial join. This notation is well defined because $[x,y]$ is determined by the hyperplanes separating $x$ and $y$. Equivalently, any two combinatorial geodesics $\alpha,\beta$ from $x$ to $y$ satisfy $\lab(\alpha)=\lab(\beta)$ and lie in $[x,y]$.
\end{rmk}

If $W_\Gamma$ is hyperbolic, equivalently if $\Gamma$ contains no induced 4-cycle, then $\widehat{X}=X^{(1)}$. In general, $\widehat{X}$ has infinite diameter if and only if $\Gamma$ is neither a clique nor a nontrivial join; see Lemma~\ref{empty Gromov boundary}.

By \cite[Theorem~3.7]{abbott2021largest}, $\widehat{X}$ is $\delta$-hyperbolic. Applying the HHS quasi-flats theorem iteratively shows that every quasi-isometry of $X$ induces a quasi-isometry of $\widehat{X}$; see \cite[Theorem 6.10]{goldsborough2023induced}.

The rest of this subsection describes the HHS structure used only in Section~\ref{BGI section}; it may be skipped on a first reading.

In a $\mathrm{CAT}(0)$ cube complex, two hyperplanes \emph{contact} if no hyperplane separates them, equivalently if their carriers intersect. The \emph{contact graph} $\CC X$ has one vertex for each hyperplane of $X$ and an edge for each contacting pair. Hagen proved that $\CC X$ is a quasi-tree \cite{hagen2014weak}, hence hyperbolic. There is a coarse Lipschitz map
\[
c:X\longrightarrow 2^{\CC X},\qquad
x\longmapsto\{H:x\in\nn(H)\}.
\]
 
 Suppose $Y$ is a convex subcomplex of $X$. Its contact graph $\CC Y$ is an induced subgraph of $\CC X$: every hyperplane of $Y$ is the intersection of $Y$ with a hyperplane of $X$, and two such hyperplanes contact in $Y$ if and only if they contact at a cube of $Y$, by the Helly property.

 Let $\mathfrak{Y}$ be a collection of convex subcomplexes of $X$. The \textit{coned-off contact graph} $\widehat{\CC X}$ of $X$ with respect to $\mathfrak{Y}$ is a graph obtained from the contact graph $\CC X$ by adding an edge between a pair of vertices in an induced subgraph $\CC Y$ for every $Y \in \mathfrak{Y}$.

Two convex subcomplexes $Y,Y'\subset X$ are \textit{parallel} if a hyperplane crosses $Y$ if and only if it crosses $Y'$. 

Now we describe the maximized hierarchically hyperbolic structure on $W_\Gamma$. Let $X= X_\Gamma$ be the Davis complex. Let $\mathfrak{L}$ be the collection of induced subgraphs $\Lambda$ of $\Gamma$ such that neither $\Lambda$ nor $\lk(\Lambda)$ is a clique. Note that if $\Lambda \in \mathfrak{L}$, then $\lk(\Lambda) \in \mathfrak{L}$. Then we collect the parallelism classes of all convex subcomplexes of $X$ of the form $gX_{\Lambda}$ where $g\in W_\Gamma$ and $\Lambda \in \mathfrak{L}$. Note that two subcomplexes $gX_\Lambda$ and $hX_\Delta$ are parallel only if $\Lambda = \Delta$ and $g^{-1}h \in W_{\Lambda \ast \lk(\Lambda)}$.

There is a set of domains, $\fS$, which index a collection of $\delta$-hyperbolic spaces to which $X$ projects. Each domain $U \in \fS$  indexes either $\widehat{X}$ or the parallelism class $[gX_\Lambda]$ or a dummy domain (see \cite[Appendix A]{abbott2021largest}). Every domain is associated with a $\delta$-hyperbolic space. The hyperbolic space associated with $\widehat{X}$ is itself. The hyperbolic space associated with $[gX_\Lambda]$ is the coned-off contact graph $\widehat{\CC gX_\Lambda}$ with respect to contact graphs of $hX_\Delta$ where $\Delta \in \mathfrak{L}, \Delta\subsetneq \Lambda$ and $hX_\Delta$ is parallel to a subcomplex of $gX_\Lambda$. The hyperbolic space associated with a dummy domain is a singleton. We denote by $d_U$ the metric in the hyperbolic space associated with the domain $U \in \fS$. Moreover,  for each domain $U\in \fS$ there is a natural coarsely well-defined Lipschitz projection, denoted by $\pi_U$, from $X$ to the hyperbolic space associated with $U$. When $U = \widehat{X}$, $\pi_U$ is the inclusion map from $X$ to $\widehat{X}$. When $U$ is a dummy-domain, $\pi_U$ is the trivial surjection. When $U = [gX_\Lambda]$, $\pi_U: X \to \widehat{\CC gX_\Lambda}$ is defined as the composition
\[
X \to gX_\Lambda \to \CC gX_\Lambda \to \widehat{\CC gX_\Lambda}
\]
of the gate map $\fg_{gX_\Lambda}:X \to gX_{\Lambda}$, the coarse Lipschitz map $c: gX_\Lambda \to \CC gX_{\Lambda}$ and the inclusion of $\CC gX_{\Lambda}$ into the coned-off space $\widehat{\CC gX_\Lambda}$. We say the projection $\pi_U$ is \emph{coarsely well-defined} because the image $\pi_U(x)$ is actually a subset of $\widehat{\CC gX_\Lambda}$ (coarsely) with uniformly bounded diameter (well-defined).

For two non-dummy domains $U,V\in \fS$, we say $U$ is \emph{nested} in $V$, written $U\sqsubseteq V$, if either $V=\widehat{X}$, or $U=[gX_{\Lambda_0}]$ and $V=[hX_{\Lambda_1}]$ for some $g,h\in W_\Gamma$ with $\Lambda_0\subset\Lambda_1$ and $g^{-1}h\in W_{\Lambda_1}$. The relation $\sqsubseteq$ is a partial order on $\fS$, with maximal element $\widehat{X}$. For each non-dummy domain $U=[gX_\Lambda]\neq\widehat{X}$, fix a vertex $\rho^U_{\widehat{X}}\in gX_\Lambda\subset\widehat{X}$. Two non-dummy domains $U=[gX_{\Lambda_0}]$ and $V=[hX_{\Lambda_1}]$ are \emph{orthogonal}, written $U\perp V$, if $\Lambda_0\subset\lk(\Lambda_1)$ and $g^{-1}h\in W_{\Lambda_1}$. Their union then lies in an unbounded product region. Two domains are \emph{transverse} if they are neither nested nor orthogonal.

 Now we recall the distance formula for $X$ from \cite[Theorem 9.1]{behrstock2017hierarchically}. It will be used in Section~\ref{BGI section}.

 We recall some necessary notation. For any $a,b\in \R$ and $K\geqslant 1,C\geqslant 0$, we say $a \asymp_{K,C} b$ if $\displaystyle \frac{a}{K} - C \leqslant b \leqslant a\cdot K+C$. For any $a,s\in \R$, 
 $\{\!\!\{ a\}\!\!\}_s = \begin{cases}
     a \quad \text{ if $a\geqslant s$}\\
     0 \quad \text{ otherwise}
 \end{cases}$.

 \begin{thm}
 [\protect{\cite[Theorem 9.1]{behrstock2017hierarchically}}](Distance Formula)
 \label{Distance Formula}
     There is a constant $s_0 \geqslant 0$ such that for every $s\geqslant s_0$, there are $K\geqslant 1, C\geqslant 0$ such that for all 0-cubes $x,y \in X^{(0)}$
     \[
     d(x,y) \asymp_{K,C} \sum_{U\in \fS} \{\!\!\{ d_U(\pi_U(x),\pi_U(y))\}\!\!\}_s.
     \]
 \end{thm}

 \begin{rmk}
 \label{finitely many relative domains}
     As a consequence of the distance formula, for each $s\geqslant s_0$, only finitely many domains $U\in\fS$ satisfy $d_U(\pi_U(x),\pi_U(y))\geqslant s$.
 \end{rmk}

%% file: sections/disconnectedness.tex
In this section, we identify exactly when $\widehat{X}_\Gamma$ has an empty or disconnected Gromov boundary. Lemma~\ref{empty Gromov boundary} shows that the boundary is empty precisely when $\Gamma$ is a clique or a nontrivial join. Proposition~\ref{Multi-ended} shows that if the boundary is nonempty but one of conditions \textup{(a)}--\textup{(c)} of Theorem~\ref{main theorem} fails, then $\widehat{X}_\Gamma$ is \emph{multi-ended}: the complement of some closed ball has at least two unbounded components. Finally, in Proposition~\ref{disconnectedness of Gromov boundary} we show that when $\widehat{X}_\Gamma$ is multi-ended, it has a disconnected Gromov boundary.

We first use contracting geodesics in $X_\Gamma$ to characterize when $\widehat{X}_\Gamma$ has infinite diameter and a nonempty Gromov boundary. We recall the definition of a contracting subset from \cite[Definition 4.2]{abbott2021largest}.

\begin{defn}[Contracting]
A subset $\gamma$ in a metric space $X$ is said to be $D-$\emph{contracting} if there exits a map $\pi_\gamma :X \to \gamma \subset X$ and constants $A,D > 0$ satisfying the following conditions:
\begin{enumerate}
    \item for any $x\in \gamma$, we have $\diam_X(x,\pi_{\gamma}(x)) < D$;
    \item if $x,y\in X$ with $d_X(x,y) <1$, then $\diam_X(\pi_\gamma(x), \pi_\gamma(y)) <D$;
    \item for all $x\in X$, if we set $R = A\cdot d(x,\gamma)$, then $\diam_X(\pi_\gamma(B_R(x))) \leq D$.
\end{enumerate}
    
\end{defn}

In this article, we only consider the contracting property of geodesics in a CAT(0) cube complex and the map $\pi_\gamma$ can be considered as the gate map.

In the following lemma, we describe when $W_\Gamma$ has contracting elements.

\begin{lem}[Existence of contracting elements]
\label{existence of contracting element}
    Let $\Lambda \subset \Gamma$ be an induced subgraph of $\Gamma$. Suppose $\Lambda$ is not a clique and is not contained in a nontrivial join. Then the special subgroup $W_\Lambda$ contains a contracting element.
\end{lem}

\begin{proof}
    Let $\mathcal{C}$ be the set of cone vertices of $\Lambda$; that is, a vertex $v\in V(\Gamma)$ is contained in $\mathcal{C}$ if and only if $v\in V(\Lambda)$  and for every $w\in V(\Lambda)\setminus \{v\}$, $w$ is adjacent to $v$. Let $\Lambda' = \Lambda \setminus \mathcal{C}$.

    It is clear that $\mathcal{C}$ spans a clique. Since $\Lambda$ is not a clique, $\Lambda'$ is nonempty. Moreover, $\Lambda'$ has no cone vertex since a vertex adjacent to all vertices of $\Lambda'$ is also adjacent to all vertices of $\mathcal{C}$, and would be a cone vertex of $\Lambda$.

    Now we claim that  $\Lambda'$ is not contained in a nontrivial join of $\Gamma$. Suppose not. That is, there is a nontrivial join subgraph $\Delta = \Delta_0 \ast \Delta_1$ of $\Gamma$ containing $\Lambda'$, where neither $\Delta_0$ nor $\Delta_1$ is a clique. Since $\Lambda'$ is not a clique, at least one of $\Lambda' \cap \Delta_0$ or $\Lambda' \cap \Delta_1$ is not a clique. Without loss of generality, we assume that $\Lambda' \cap \Delta_0$ is not a clique. Every vertex of $\Lambda' \cap \Delta_0$ is adjacent to every vertex of $\Delta_1 \cup \mathcal{C}$. Hence, $\Lambda$ is contained in a nontrivial join $(\Lambda' \cap \Delta_0) \ast (\Delta_1\cup \mathcal{C})$. Contradiction.

    Assume the vertex set $V(\Lambda') = \{v_1,\cdots, v_n\}$. Let 
    \[
    g = v_1\cdots v_n \in W_{\Lambda'}
    \]
    be the Coxeter element of $W_{\Lambda'}$. For every $m\in \N$, $g^n$ is geodesic. Indeed, if $g^m$ is not geodesic for some $n$, then the Deletion Condition (see Lemma~\ref{Deletion Condition}) would provide two occurrences of some letter $v_i$ such that $v_i$ commutes with every letter in between. For an innermost such pair, this would imply that $v_i$ commutes with every vertex in $V(\Lambda')\setminus \{v_i\}$, and hence $v_i$ would be a cone vertex of $\Lambda$. Contradiction.

    The axis $\alpha\subset X_\Gamma$ labeled by
    \[
    \cdots (v_1\cdots v_n)(v_1\cdots v_n)(v_1\cdots v_n)\cdots
    \]
    is an axis for $g$. It is clear that every subpath of $\alpha$ of length at least $n$ contains all labels in $V(\Lambda')$, which is not contained in a nontrivial join. Hence, every subpath of $\alpha$ with length at least $n$ is not contained in a standard product region. By \cite[Theorem 4.2]{charney2015contracting}, $\alpha\subset X_\Gamma$ is contracting, and hence $g\in W_\Gamma$ is a contracting element.
\end{proof}

\begin{rmk}
\label{exist a contracting element if and only if not a join}
    When the defining graph $\Gamma$ is a nontrivial join, every axis in $X_\Gamma$ for any infinite order element of $W_\Gamma$ bounds a half-Euclidean plane. Hence it is not contracting. Therefore, by applying Lemma~\ref{existence of contracting element}, we conclude that $W_\Gamma$ contains a contracting element if and only if $\Gamma$ is neither a nontrivial join nor a clique.
\end{rmk}

Using the preceding remark, we characterize when the Gromov boundary $\partial \widehat{X}$ of the maximized hyperbolic space is nonempty.

\begin{lem}[Nonempty boundary $\partial\widehat{X}_\Gamma$]
\label{empty Gromov boundary}
The maximized hyperbolic space $\widehat{X}_{\Gamma}$ has infinite diameter and a nonempty Gromov boundary if and only if $\Gamma$ is neither a clique nor a nontrivial join. Whenever nonempty, $\partial\widehat{X}_\Gamma$ contains at least two points.
\end{lem}

\begin{proof}
    If $\Gamma$ is a clique, then both $X_\Gamma$ and $\widehat{X}_\Gamma$ are bounded. If $\Gamma$ is a nontrivial join, then $X_\Gamma$ is a product of two unbounded subspaces and is itself a standard product region, so $\diam(\widehat{X}_\Gamma)=1$. Thus the boundary is empty in either case.

    Conversely, suppose $\Gamma$ is neither a clique nor a nontrivial join. By Lemma~\ref{existence of contracting element}, there is a contracting geodesic $\alpha$ in $X_\Gamma$. By~\cite[Theorem D]{abbott2021largest}, the image of $\alpha$ in $\widehat{X}_{\Gamma}$ under the natural identity map is a quasi-geodesic. Thus, $\widehat{X}_{\Gamma}$ has infinite diameter. Moreover, $\partial \widehat{X}_{\Gamma}$ is nonempty and contains at least two points.
\end{proof}

For the rest of the section, assume that $\Gamma$ is neither a clique nor a nontrivial join. A subspace $A$ of a metric space $X$ \emph{coarsely separates} $X$ if, for some $K\in\N$, the complement $X\setminus N_K(A)$ has at least two components, which are  not contained in a finite neighborhood of $A$. The following is a classical result from Bass-Serre theory. We refer the reader to \cite[Lemma 1.8]{scott1977ends} and \cite[Lemma 2.2]{papasoglu2012splittings} for a proof.

\begin{lem}
\label{splitting implies coarse separation}
    If a finitely generated group $G$ splits over a finitely generated subgroup $H$, then $H$ coarsely separates the Cayley graph of $G$.
\end{lem}

If $\Delta$ separates $\Gamma$, then $W_\Gamma$ splits over $W_\Delta$ and $X_\Delta$ coarsely separates $X_\Gamma$. When $\Delta$ is a clique or lies in a nontrivial join, its image in $\widehat{X}_\Gamma$ has bounded diameter and also coarsely separates $\widehat{X}_{\Gamma}$. Adapting the argument of
\cite[Lemma~2.2]{papasoglu2012splittings}, we moreover find quasi-geodesic rays
in two distinct unbounded complementary components.

\begin{prop}
\label{Multi-ended}
    Let $W_\Gamma$ be the right-angled Coxeter group with defining graph $\Gamma$. Suppose $\Gamma$ is not a clique, not a nontrivial join, and one of the conditions (a), (b) and (c) from Theorem~\ref{main theorem} fails. Then there is a constant $K'\geq 0$ and a vertex  $o\in \widehat{X}_{\Gamma}$ such that $\widehat{X}_\Gamma\setminus B_{K'}(o)$  has at least two unbounded components. Moreover, at least two unbounded components contain quasi-geodesic rays.
\end{prop}
\begin{proof}
    If condition~\textup{(a)} fails, then $\Gamma$ is disconnected or has a separating clique (the clique case is excluded by hypothesis). If condition~\textup{(b)} fails, then $\Gamma$ has a separating nontrivial join. Finally, suppose condition~\textup{(c)} fails, so that $\lk(v)\subset J$ for some vertex $v$ and some nontrivial join $J$. Because $\Gamma$ itself is not a nontrivial join, $\Gamma\setminus\st(v)$ is nonempty. Thus $\lk(v)$ separates $v$ from $\Gamma\setminus\st(v)$. In every case there is a separating subgraph $\Delta\subset\Gamma$ that is either a clique or is contained in a nontrivial join. When $\Gamma$ is disconnected, we take $\Delta$ to be the empty clique.
    
    Let $\Gamma_0$ be the union of $\Delta$ and a component of $\Gamma\setminus \Delta$. Let $\Gamma_1 = (\Gamma\setminus \Gamma_0) \cup \Delta$. Hence $\Gamma$ is the union $\Gamma = \Gamma_0 \sqcup_{\Delta} \Gamma_1$, and the right-angled Coxeter group $W_\Gamma$ splits as $W_\Gamma \cong W_{\Gamma_0} \ast_{W_{\Delta}}W_{\Gamma_1}$.

    Let $X$ be the 1-skeleton of the Davis complex $X_{\Gamma}$. We identify $X$ with Cayley graph of $W_{\Gamma}$ with respect to the standard generating set $V(\Gamma)$ and let $X_\Delta \subset X$ be the subgraph spanned by the vertex set $W_{\Delta}$. We fix a vertex $o\in X_{\Delta} \subset X$. Let $T$ be the Bass-Serre tree of $W_\Gamma$ for the splitting $W_\Gamma \cong W_{\Gamma_0} \ast_{W_{\Delta}}W_{\Gamma_1}$. Note that the vertex set of $T$ is $\{gW_{\Gamma_i} \mid g\in W_{\Gamma}, i =0,1\}$, the edge set of $T$ is $\{gW_{\Delta}\mid g\in W_{\Gamma}\}$, and $gW_{\Delta}$ connects $gW_{\Gamma_0}$ and $gW_{\Gamma_1}$. $W_\Gamma$ acts on $T$ by automorphisms. 

    Now we want to show $X_\Delta$ coarsely separates $X$. Let $X'$ be the barycentric subdivision of $X$ and $T'$ be the barycentric subdivision of $T$. Let $e$ be the edge of $T$ stabilized by $W_{\Delta}$. The midpoint $a$ of $e$ is a vertex of $T'$. Since every vertex of $X$ is a group element of $W_\Gamma$, we define a map $p:X^{(0)}\to T'$ by $p(g) = g\cdot a$. For every generator $v\in V(\Gamma)$ , $v\cdot a = a$ if $v\in V(\Delta)$, $v\cdot a$ is the midpoint of the edge connecting $W_{\Gamma_0}$ and $vW_{\Gamma_1}$ if $v\in V(\Gamma_0\setminus \Delta)$, and $v\cdot a$ is the midpoint of the edge connecting $W_{\Gamma_1}$ and $vW_{\Gamma_0}$ if $v\in V(\Gamma_1\setminus \Delta)$. Therefore, $d_{T'}(g\cdot a,gv\cdot a) = 0$ or $2$ for each generator $v\in V(\Gamma)$ and every group element $g\in W_{\Gamma}$. Then every edge $(g,gv)$ in $X$ can be mapped by $p$ to either a point or a path with length $2$ in $T'$. Hence, we can extend $p$ to a simplicial map from $X'$ to $T'$. By construction, $p^{-1}(a) = X_{\Delta}$.  Furthermore, for vertices $g,h\in X$, $d_{T'}(p(g),p(h))\leq d_{X'}(g,h)$. 

    Let $s\in \Gamma_0\setminus \Delta$ and $t\in \Gamma_1\setminus \Delta$. In particular, $s$ and $t$ are non-adjacent. Let $\Pi$ be the maximal subgraph of $\Gamma$ consisting of cone vertices of $\Gamma$. Let $\Lambda = \Gamma\setminus \Pi$. Assume $V(\Lambda) = \{s,t,v_1,\cdots,v_n\}$. Let $g = sv_1\cdots v_n t$ and $h = tv_1\cdots v_n s$. $g^m$ and $h^m$ are geodesic for all $m\in \N$. Let $m\in \N$. Since the first letter of $g$ is $s\in \Gamma_0\setminus \Delta$, the geodesic in $T'$ from $a$ to $p(g^m)$ passes through the vertex $W_{\Gamma_0}$. Similarly, the geodesic from $a$ to $p(h^m)$ passes through $W_{\Gamma_1}$. Hence, $p(g^m)$ and $p(h^m)$ are in different components of $T'\setminus \{a\}$. Let $\beta$ be a path from $g^m$ to $h^m$ in $X$. $\beta$ has a unique image in $X'$. Then $p(\beta)$ passes through $a$. Hence, $\beta$ intersects $X_{\Delta}$. By the construction of $g$ and $h$, $d_{T'}(a,g\cdot a)\geq 2$, $d_{T'}(a,h\cdot a)\geq 2$ and they have cyclically reduced amalgam normal forms of positive syllable length. Hence, $g$ and $h$ act on $T'$ loxodromically with translation length $\geq 2$ and  $d_{T'}(a,g^m\cdot a)\geq m$ and $d_{T'}(a,h^m\cdot a)\geq m$. Since $p$ is distance non-increasing, $d_X(g^m,X_{\Delta})\geq m/2$ and $d_X(h^m,X_{\Delta})\geq m/2$ for all $m\in \N$. Hence, $X_{\Delta}$ coarsely separates $X$.

    Let $\gamma_0$ and $\gamma_1$ be the geodesic rays obtained as limits of $(g^m)$ and $(h^m)$, respectively. The image in $T'$ of every path from $\gamma_0$ to $\gamma_1$ passes through $a$, so every such path meets $X_\Delta$. Thus the rays lie in different deep components of $X\setminus X_\Delta$. Because $\lab(g)=\lab(h)=V(\Lambda)$ is not contained in a nontrivial join, both rays are contracting, as in Lemma~\ref{empty Gromov boundary}.

     Let $\pi:X\to \widehat{X}$ be the inclusion map. Clearly, $\pi(X_{\Delta})$ has uniformly bounded diameter. Indeed, if $\Delta$ is a clique, then $\diam(\pi(X_{\Delta})) \leq M$, where $M$ is the maximum cardinality of a clique in $\Gamma$. And if $\Delta$ is a subgraph of a nontrivial join, then $\diam(\pi(X_{\Delta})) \leq 1$. Then we let $K =  M +2$. We want to show that $\widehat{X}\setminus B_{K}(o)$ has at least two unbounded components.

     Let $x,y$ be vertices in distinct deep components of $X\setminus X_\Delta$ such that the images $\pi(x)$ and $\pi(y)$ in $\widehat{X}$ are not contained in the ball $B_{K}(o)$. We identify the vertex set of $X$ with the vertex set of $\widehat{X}$. Suppose there is a path $x=x_0,\cdots,x_n=y$ from $x$ to $y$ in $\widehat{X}$ avoiding $B_{K}(o)$. For each $i=1,\cdots,n$, the vertices $x_{i-1}$ and $x_i$ in $X$ are either connected by a 1-cube or contained in the same standard product region. Let $\alpha_i$ be a combinatorial geodesic in $X$ connecting $x_{i - 1}$ and $x_i$ for each $i$. Then, the path $\alpha= \alpha_1\cdots \alpha_n$ obtained by concatenation is a path in $X$ from $x$ to $y$. Since $x_i$'s are not in $X_{\Delta}$, there is a pair $x_{i-1}$ and $x_i$ that are contained in different components of $X\setminus X_{\Delta}$ for some $i$ and the intersection $\alpha_i\cap X_{\Delta}$ is nonempty. In particular, $\alpha_i$ is not a 1-cube. Then $x_{i-1}$ and $x_i$ are contained in the same standard product region $P$.

     Let $z\in \alpha_i \cap X_\Delta$. Since $\alpha_i$ is geodesic and $P$ is convex, $\alpha_i\subset P$. Hence, $d_{\widehat{X}}(z,x_i) = 1$. Since $\diam_{\widehat{X}}(\pi(X_{\Delta})) \leq M$, $d_{\widehat{X}}(o,x_i)\leq M+1 \leq K$. Contradiction. Therefore, the images under $\pi$ of distinct deep components of $X \setminus X_{\Delta}$ are separated by $B_K(o)$ in $\widehat{X}$. By~\cite[Theorem D]{abbott2021largest}, $\pi(\gamma_0)$ and $\pi(\gamma_1)$ are quasi-geodesic rays in distinct components of $\widehat{X}\setminus B_K(o)$.
\end{proof}

We make a remark on right-angled Artin groups.

\begin{rmk}
\label{RAAG boundary}
    By an argument similar to the one above, we could show that the maximized hyperbolic space of a right-angled Artin group $A_\Gamma$ is either empty, or has a bounded separating set, hence a disconnected Gromov boundary as follows. If the defining graph $\Gamma$ is disconnected, then $A_{\Gamma}$ has infinitely many ends. In this case, a bounded separating set in $A_{\Gamma}$ has an image in the maximized hyperbolic space which is also separating. If the defining graph $\Gamma$ is a join, then $A_{\Gamma}$ is a direct product of two infinite subgroups. Hence the maximized hyperbolic space has finite diameter and the Gromov boundary is empty. Now suppose $\Gamma$ is connected, not a join and has more than one vertex. Then the star of each vertex generates an unbounded product subgroup, i.e. for a vertex $v\in V(\Gamma)$, $A_{\st(v)} \cong A_{v}\times A_{\lk(v)}$ where both $A_{v}$ and $A_{\lk(v)}$ are infinite. The stabilizer of each hyperplane of type $v$ is a conjugate of $A_{\lk(v)}$. Hence each hyperplane is contained in a standard product region, a left coset of a special subgroup containing $A_{\st(v)}$. Then the image of $A_{\lk(v)}$ has a bounded diameter in the maximized hyperbolic space. On the other hand, let $\Gamma_0 = \st(v)$ and $\Gamma_1 = \Gamma\setminus \{v\}$. This provides a group splitting $A_\Gamma \cong A_{\Gamma_0}\ast_{A_{\lk(v)}} A_{\Gamma_1}$. Hence $A_\Gamma$ is coarsely separated by $A_{\lk(v)}$. Therefore the image of $A_{\lk(v)}$ in the maximized hyperbolic space is a bounded separating subset. 
    
    We remark that the maximized hyperbolic space of an RAAG is quasi-isometric to the contact graph and also to the extension graph (\cite[Section 5]{abbott2021largest}), which is a quasi-tree.
\end{rmk}

We conclude the section with the forward implication of
Theorem~\ref{main theorem}.

\begin{prop}
\label{disconnectedness of Gromov boundary}
    Let $W_\Gamma$ be a right-angled Coxeter group. Suppose that at least one
    of the following holds:
    \begin{itemize}
        \item[\textup{(1)}] $\Gamma$ is a clique, is disconnected, or has a separating clique;
        \item[\textup{(2)}] $\Gamma$ has a separating nontrivial join;
        \item[\textup{(3)}] $\lk(v)$ is contained in a nontrivial join for some $v\in V(\Gamma)$.
    \end{itemize}
    Then the Gromov boundary $\partial\widehat{X}$ of the maximized hyperbolic space is either empty or disconnected. Moreover, $\partial \widehat{X}$ is empty if and only if $\Gamma$ is a clique or a nontrivial join.
\end{prop}
\begin{proof}
    The final assertion follows from Lemma~\ref{empty Gromov boundary}. Assume that $\Gamma$ is neither a clique nor a nontrivial join, and fix $o\in\widehat{X}$. By Proposition~\ref{Multi-ended}, some complement $\widehat{X}\setminus B_K(o)$ has more than one unbounded component. Let $(A_i)_{i\in\N}$ be its unbounded components, and let $\partial A_i$ consist of boundary classes represented by sequences eventually contained in $A_i$. We may assume that $A_1$ and $A_2$ contain quasi-geodesic rays, so $\partial A_1$ and $\partial A_2$ are nonempty. We claim that $\partial\widehat{X}=\bigcup_{i\in\N}\partial A_i$.

    Let $(x_n)$ represent a point of $\partial\widehat{X}$. It suffices to show that the sequence eventually lies in a single component $A_i$. Otherwise, some subsequence visits more than one component. Whenever two terms of this subsequence lie in different components, every geodesic between them intersects $B_K(o)$. Hyperbolicity then gives a constant $K'=K'(K,\delta)$ such that their Gromov product based at $o$ is at most $K'$. This contradicts the fact that $(x_n)$ converges at infinity.

    If $(y_n)$ represents a point of $\partial A_1$ and $(z_n)$ represents a point of $\partial A_i$ with $i\neq1$, then $\langle y_n,z_m\rangle_o\leq K'$ for all sufficiently large $n,m$. Hence any visual metric separates $\partial A_1$ uniformly from the other $\partial A_i$. Thus $\partial A_1$ is a nonempty proper open and closed subset of $\partial\widehat{X}$. Therefore, the boundary is disconnected.
\end{proof}

%% file: sections/Morse-to-Gromov.tex
In this section, we deduce the connectivity of certain Gromov boundaries from the connectivity of Morse boundaries. We record a general formulation for cobounded Morse-recognizing actions on hyperbolic spaces, using known results on boundary maps and standard facts about loxodromic dynamics.

We apply this formulation to maximized hyperbolic spaces of HHGs and to the hyperbolic models for Coxeter groups constructed by Cashen--Chu--Tao \cite{cashen2026hyperbolic}. For right-angled Coxeter groups, which are contained in both categories, this gives an alternative proof of the connectedness conclusion of Theorem~\ref{main theorem}\footnote{In a recent work, Cashen--Chu--Tao classify which Coxeter groups are HHGs \cite{cashen2026hierarchicalhyperbolicitycoxetergroups}.}. The direct geometric argument in Section~\ref{Main Proof Section} establishes the stronger conclusion of linear connectedness.

\subsection{Density for Morse-recognizing actions}

Let $G$ be a finitely generated group and $\operatorname{Cay}(G,A)$ be the Cayley graph of $G$ with respect to a finite generating set $A$. We denote by $\partial_M G$ the Morse boundary of $G$ with the direct limit topology. We refer readers to \cite{cordes2017morse} for details of the Morse boundary.

Suppose that $G$ acts on a geodesic hyperbolic space $Y$ by isometries. Denote by $\partial Y$ the Gromov boundary of $Y$. Fix $y_0\in Y$ and consider the orbit map
\[
    \pi\colon G\longrightarrow Y,
    \quad \text{ defined by } \quad
    \pi(g)=g\cdot y_0.
\]
Note that $\pi$ is Lipschitz since $G$ is finitely generated.

\begin{defn}
\label{Morse recognizing orbit map}
We say that the orbit map $\pi$ is \emph{Morse recognizing} if, for every $(k,c)$-quasi-geodesic $\gamma\colon I\to\operatorname{Cay}(G,A)$, the following hold:
\begin{enumerate}
    \item For each Morse gauge $N$, if $\gamma$ is $N$-Morse, then $\pi\circ\gamma$
    is a quasi-geodesic with constants depending only on
    $N$, $k$, and $c$.

    \item If $\pi\circ\gamma$ is a $(K,C)$-quasi-geodesic,
    then $\gamma$ is Morse with a gauge depending only on
    $k$, $c$, $K$, and $C$.
\end{enumerate}
The interval $I$ may be bounded, a ray, or the whole real line.
\end{defn}

This is the Morse-detectability property from \cite[Definition~A.4]{russell2025extensions}, applied to an orbit map.

\begin{prop}[Density for Morse-recognizing actions]
\label{dense Morse boundary}
Let $G$ be a finitely generated group acting coboundedly on a geodesic hyperbolic space $Y$ by isometries. Suppose that the orbit map $\pi\colon G\to Y$ is Morse recognizing. Then $\pi$ induces a natural $G$-equivariant continuous injection
\[
    \iota\colon\partial_M G\longrightarrow\partial Y
\]
with a dense image.
\end{prop}

\begin{proof}
By \cite[Lemma~A.6]{russell2025extensions}, Morse recognition gives a continuous injection
\[
    \iota\colon\partial_M G\longrightarrow\partial Y.
\]
For a Morse geodesic ray $\gamma$ in $G$, the point $\iota([\gamma])$ is the endpoint of the quasi-geodesic ray $\pi\circ\gamma$. By the definition of the orbit map $\pi$, the induced boundary map is $G$-equivariant.

Now we want to show the image $\iota(\partial_M G)$ is dense in $\partial Y$ by proving  that it contains all fixed points of loxodromic elements and that these points are dense in $\partial Y$.

If $Y$ has finite diameter, then $\partial Y$ is empty. Moreover, $\partial_M G$ is empty, since a Morse geodesic ray in $\operatorname{Cay}(G,A)$ would map to an unbounded quasi-geodesic ray in $Y$. Thus the conclusion holds in this case. Henceforth assume that $Y$ has infinite diameter.

Let $g\in G$ act loxodromically on $Y$ with fixed points $g^+,g^-\in \partial Y$. There exist constants $a> 0$ and $b\geq 0$ such that
\[
     a|m-n|-b \leq d_Y(g^m y_0,g^n y_0)
    \qquad\text{for all }m,n\in\mathbb Z.
\]
Set
\[
    L=\max\bigl(\{1\}\cup
       \{d_Y(y_0,s y_0):s\in A\}\bigr).
\]
Since the orbit map $\pi$ is $L$-Lipschitz, we have
\[
    \frac{a|m-n|-b}{L}
    \leq d_G(g^m,g^n)
    \leq |g|_A\,|m-n|.
\]
Hence, $n\mapsto g^n$ is a quasi-isometric embedding of $\Z$ into $G$. We extend it to a quasi-geodesic $\gamma_g: \R\to \operatorname{Cay}(G,A)$ such that $\gamma_g (n) = g^n$ by concatenating consecutive pairs $g^n, g^{n+1}$ by a fixed geodesic path representing the word $g\in G$. 

Clearly, $\pi\circ \gamma_g$ is a quasi-geodesic in $Y$. By the Morse-recognizing property, $\gamma_g\subset \operatorname{Cay}(G,A)$  is a Morse quasi-geodesic. Its two half-rays determine points $\gamma_g^+,\gamma_g^-\in\partial_M G$, and
\[
    \iota(\gamma_g^+)=g^+,
    \qquad
    \iota(\gamma_g^-)=g^-.
\]
Thus every loxodromic fixed point lies in
$\iota(\partial_M G)$.

Let $\Lambda_Y(G)$ denote the limit set of the action $G \curvearrowright Y$. Coboundedness implies
\[
    \Lambda_Y(G)=\partial Y.
\]
Indeed, any sequence representing a boundary point can be approximated within a uniform distance by points of the orbit $G y_0$.

By \cite[Theorem~4.2]{abbott2019hyperbolic}, an unbounded cobounded action has either a two-point or an infinite limit set. In the two-point case, these points are the fixed points of a loxodromic element, so both belong to $\iota(\partial_M G)$. In the infinite case, the attracting fixed points of loxodromic elements are dense in $\Lambda_Y(G)$ \cite[Proposition~7.4.6]{das2017geometry}. Since these points belong to $\iota(\partial_M G)$ and $\Lambda_Y(G)=\partial Y$, the image of $\iota$ is dense in $\partial Y$.
\end{proof}

\begin{cor}
\label{connected Morse HHS implies connected gromov boundary}
Under the hypotheses of Proposition~\ref{dense Morse boundary}, if $\partial_M G$ is connected, then $\partial Y$ is connected. If $\partial_M G$ is also nonempty, then $\partial Y$ is nonempty.
\end{cor}

\begin{proof}
The continuous image $\iota(\partial_M G)$ is connected. By Proposition~\ref{dense Morse boundary},
\[
    \partial Y=\overline{\iota(\partial_M G)}.
\]
The closure of a connected subset is connected, so $\partial Y$ is connected. Nonemptiness follows from the definition of $\iota$.
\end{proof}

\subsection{Applications to hierarchically hyperbolic groups}

We apply the preceding results to the maximized structure on a hierarchically hyperbolic group. We refer readers to \cite[Section 3]{abbott2025structure} for related definitions.

\begin{cor}
\label{connected Morse implies connected maximized boundary}
Let $G$ be a hierarchically hyperbolic group equipped with its maximized structure, and let $\CCS$ be its maximized hyperbolic space. Then the natural map
\[
    \partial_M G\longrightarrow\partial\CCS
\]
is a continuous injection with a dense image. In particular, if $\partial_M G$ is nonempty and connected, then $\partial\CCS$ is nonempty and connected.
\end{cor}

\begin{proof}
For the maximized structure, the action on $\CCS$ is cobounded, and the top-level projection is Morse-recognizing by \cite[Corollary~6.2]{abbott2021largest}. The projection is coarsely equivariant, so it is at a uniformly bounded distance from an orbit map. This orbit map is therefore also Morse recognizing. The conclusions follow from Proposition~\ref{dense Morse boundary} and Corollary~\ref{connected Morse HHS implies connected gromov boundary}.
\end{proof}

\begin{rmk}
For an HHG containing a Morse element, the density assertion has been made in  \cite{durham2017boundaries} and \cite[Remark 5.11]{abbott2025structure}. To see this, let $\mathfrak T$ be the maximized structure with nested-maximal hyperbolic space $\CCS$ and let $\partial(G,\mathfrak T)$ be its  HHS boundary. Both maps in
\[
    \partial_M G \into \partial\CCS \into\partial(G,\mathfrak T)
\]
are continuous injections with dense images. 
\end{rmk}

\subsection{Hyperbolic models for Coxeter groups}

Proposition~\ref{dense Morse boundary} also applies to hyperbolic models for Coxeter groups. 

A recent result by Cashen--Chu--Tao shows that every Coxeter group admits a hyperbolic space which is Morse recognizing in \cite{cashen2026hyperbolic}. Let $(W, S)$ be a finite-rank Coxeter system with Davis complex $\Sigma$. The hyperbolic model $\widehat{\Sigma}$ is obtained by coning off the wide parabolic subcomplexes of $\Sigma$ \cite[Definition~4.18]{cashen2026hyperbolic}. Note that if the Coxeter group is right-angled, this hyperbolic model $\widehat{\Sigma}$ is quasi-isometric to the maximized hyperbolic space $\widehat{X}$ (Definition \ref{Definition of maximized hyperbolic space}). 

\begin{cor}
\label{Coxeter model connected boundary}
Let $(W,S)$ be a finite-rank Coxeter system, and let $\widehat{\Sigma}$ be the hyperbolic model of \cite{cashen2026hyperbolic}. Then the natural map
\[
    \partial_M W\longrightarrow\partial\widehat{\Sigma}
\]
is a continuous injection with dense image. Consequently, if $\partial_M W$ is nonempty and connected, then $\partial\widehat{\Sigma}$ is nonempty and connected.

\end{cor}

\begin{proof}
Equip $\Sigma^{(1)}$ with its combinatorial metric and identify it with $\operatorname{Cay}(W,S)$. By \cite[Corollary~4.20]{cashen2026hyperbolic}, the space $\widehat{\Sigma}$ is geodesic and hyperbolic, and the inclusion
\[
    \Sigma^{(1)}\longrightarrow\widehat{\Sigma}
\]
is Morse recognizing. Clearly, $W$ acts on $\widehat{\Sigma}$ coboundedly. The conclusions follow from Proposition~\ref{dense Morse boundary} and Corollary~\ref{connected Morse HHS implies connected gromov boundary}.
\end{proof}

Combining this result with \cite[Corollary 4.10]{cordes2026morseaccessibilityapplicationscoxeter} gives the following corollary.

\begin{cor}
\label{Coxeter group connected boundary corollary}
    Let $W$ be a one-ended, wide-spherical-avoidant Coxeter group. Then the Gromov boundary $\partial \widehat{\Sigma}$ of its Morse-recognizing space is connected.
\end{cor}

 For the definition of \emph{wide-spherical-avoidant}, we refer readers to \cite[Definition 3.5]{cordes2025connectivity}.

In the case of right-angled Coxeter groups, combining the preceding corollaries with the characterization of Cordes--Levcovitz recovers the connectedness implication of our main theorem.

\begin{cor}
\label{direct proof of connectedness}
Let $W_\Gamma$ be a right-angled Coxeter group whose defining graph satisfies conditions~\textup{(a)}--\textup{(c)} of Theorem~\ref{main theorem}. Then $\partial\widehat{X}_\Gamma$ is nonempty and connected.
\end{cor}

\begin{proof}
Condition~\textup{(a)} of Theorem~\ref{main theorem} and Theorem~\ref{one ended} imply that $W_\Gamma$ is one-ended. Conditions~\textup{(b)} and~\textup{(c)} then imply that $\Gamma$ is wide-avoidant, by Lemma~\ref{conditions are equivalent}.

By \cite[Theorem~B]{cordes2025connectivity}, the Morse boundary $\partial_M W_\Gamma$ is nonempty and connected. Corollary~\ref{connected Morse implies connected maximized boundary} or Corollary~\ref{Coxeter model connected boundary} shows that $\partial\widehat{X}_\Gamma$ is nonempty and connected.
\end{proof}

The arguments in this section establish the connectedness of the  Gromov boundary. Continuity and density alone do not imply path connectedness, local connectedness (for instance, consider the topologist's sine curve), or linear connectedness.

%% file: sections/relation.tex
The purpose of this section is to establish tools relating the cubical geometry of $X$ to the fine geometry of the maximized hyperbolic space $\widehat{X}$. These tools will be used to establish visibility in Section~\ref{BGI section} and to control the fan-and-filter construction in Section~\ref{Fan and Filter Section}. The main difficulty in studying $\widehat{X}$ is that coning off standard product regions in $X$ simplifies the large-scale geometry but generally destroys properness and can produce infinitely many geodesics between two vertices in $\widehat{X}$; see Example~\ref{exist infinitely many geodesics}. We therefore need controlled representatives of geodesics in $\widehat{X}$, together with estimates that give lower bounds on distance in $\widehat{X}$ in terms of combinatorial information in $X$.

Recall that we identify the vertex set of the Davis complex $X$ with the vertex set of $\widehat{X}$ via the natural inclusion $X \hookrightarrow\widehat{X}$.

This section proceeds in three parts. First, we introduce \emph{taut geodesics} in $\widehat{X}$ (Definition~\ref{taut geodesic definition}) and prove their existence in Proposition~\ref{taut geodesic existence}. The vertices of a taut geodesic lie in the cubical convex hull in $X$ of its endpoints. Since the convex hull of two vertices is finite, only finitely many taut geodesics join any fixed pair of vertices. Thus, taut geodesics allow us to study geodesics in $\widehat{X}$ by using the hyperplanes of a finite cubical convex hull.

Second, we study cubical convex hulls and gate maps with respect to the metric $d_{\widehat{X}}$.  Lemma~\ref{Projections don't increase the distance} shows that gate maps to convex subcomplexes are $1$--Lipschitz with respect to $d_{\widehat{X}}$. Lemma~\ref{distance monotonicity in hull} shows that the distance from a fixed basepoint is nondecreasing along combinatorial geodesics starting at that basepoint. We use these estimates to find a point far from the basepoint in the intersection of two convex hulls whose terminal vertices are close in $\widehat{X}$ (Lemma~\ref{large intersection of convex hull}). We also prove that every point of a cubical convex hull lies uniformly close in $\widehat{X}$ to a taut geodesic joining its endpoints (Lemma~\ref{convex hulls are uniformly thin}).

Finally, we introduce \emph{wide subpaths}, which record the portions of a combinatorial geodesic that lie inside standard product regions and are therefore collapsed by the coning-off procedure. Lemma~\ref{Morse Geodesic C'} shows that a uniform bound on the lengths of the wide subpaths of one combinatorial geodesic gives a uniform bound for every combinatorial geodesic with the same endpoints. Using taut geodesics, Lemma~\ref{C'} then converts this combinatorial information into progress in $\widehat{X}$. More precisely, if the wide subpaths of every combinatorial geodesic from $x$ to $y$ have length at most $C'$, then
\[
d_X(x,y)\geq n(C'+1) \quad\Longrightarrow\quad d_{\widehat{X}}(x,y)\geq n.
\]
These results will control the geometry of the filters that we will construct in Section~\ref{Fan and Filter Section} and ensure that sufficiently deep layers of a filter remain outside of prescribed balls in $\widehat{X}$.

Fix $o\in\widehat{X}$ and write
\[
S_r(o)=\{x \mid d_{\widehat{X}}(o,x)=r\},\qquad
B_r(o)=\{x \mid d_{\widehat{X}}(o,x)\leq r\}.
\]

\subsection{Taut geodesics in $\widehat{X}$}
\begin{defn}[Taut geodesic]
\label{taut geodesic definition}
    Let $x,y\in \widehat{X}$. A geodesic $x=x_0,\ldots,x_n=y$ in $\widehat{X}$ is \emph{taut} if every
    vertex $x_i$ lies in the cubical convex hull $[x,y]\subset X$.
\end{defn}

It is an immediate consequence of the above definition that, since the convex hull of two points in $X$ is finite, for any pair of vertices $x,y \in \widehat{X}$ there are only finitely many taut geodesics from $x$ to $y$ in $\widehat{X}$. 

Now we show the existence of taut geodesics. The proof is inspired by \cite[Lemma 4.11]{hagen2022large}.

\begin{prop}[Existence of taut geodesics]
\label{taut geodesic existence}
Let $x,y\in \widehat{X}$. There is a combinatorial geodesic $\alpha$ in $X$ from $x$ to $y$ and a geodesic $x= x_0,x_1,\cdots,x_n=y$ in $\widehat{X}$ such that each $x_i$ is a vertex of $\alpha$. In particular, this geodesic in $\widehat{X}$ is taut.
\end{prop}
\begin{proof}
    Let $x= x_0,x_1,\cdots, x_n=y$ be a geodesic in $\widehat{X}$ from $x$ to $y$, where each $x_i$ is a vertex and a consecutive pair is joined by an edge. Therefore, for each $1\leq i \leq n$, $x_{i-1}$ and $x_i$ are either adjacent by a 1-cube in $X$ or contained in the same standard product region $P$. Let $\alpha_i$ be a combinatorial geodesic from $x_{i-1}$ to $x_i$ in $X$. Note that $\alpha_i$ is either a 1-cube containing $x_{i-1}$ and $x_i$ or a combinatorial geodesic contained in a product region $P$ from $x_{i-1}$ to $x_i$. Then $\alpha = \alpha_1\cdots \alpha_n$ obtained by concatenation is a combinatorial path from $x$ to $y$ in $X$.

    Define the complexity of the pair $((x_0,\cdots,x_n),(\alpha_1,\cdots,\alpha_n))$  to be the $n$-tuple
    \[
    (|\alpha_1|,|\alpha_2|,\cdots, |\alpha_n|),
    \]
    ordered lexicographically. Suppose that the geodesic $x_0,x_1,\cdots, x_n$ in $\widehat{X}$ and combinatorial geodesics $\alpha_1,\cdots,\alpha_n$ in $X$ have been chosen to minimize the complexity. This means that we prioritize $\alpha_1$ being as short as possible, and among the choices such that $\alpha_1$ is as short as possible, we minimize the length of $\alpha_2$ and continue. We claim that $\alpha$ is a combinatorial geodesic in $X$ from $x$ to $y$.

    It is enough to show that every hyperplane $W$ crosses $\alpha$ at most once. Suppose a hyperplane $W$ crosses $\alpha$ more than once at an edge $e$ and then at an edge $f$. Suppose $e\subset \alpha_i$ and $f\subset \alpha_j$ with $i < j$. Note that $i \neq j$ by construction since $e\neq f$ and each $\alpha_i$ is geodesic.

    \noindent
    \textbf{Case 1: $j=i+1$.} If both $\alpha_i$ and $\alpha_{i+1}$ are
    single edges of $X$, then the two edges incident to $x_i$ have the same
    label and hence coincide. Thus $x_{i-1}=x_{i+1}$, contradicting that
    $x_0,\ldots,x_n$ is geodesic in $\widehat{X}$.

    Suppose, without loss of generality, that $\alpha_i$ is a combinatorial geodesic lying in a standard product region $P$ and $\alpha_{i+1}$ is a 1-cube. By the convexity of the product region, the latter 1-cube has the same label as $e$, hence lies in $P$. Thus $x_{i-1}$ and $x_{i+1}$ lie in the same standard product region and are adjacent in $\widehat{X}$, again contradicting geodesicity.

    Now suppose both $\alpha_i$ and $\alpha_{i+1}$ are not 1-cubes. The hyperplane $W$ separates the CAT(0) cube complex $X$ into two halfspaces $\overleftarrow{W}$ and $\overrightarrow{W}$. Suppose $x_{i-1},x_{i+1} \in \overrightarrow{W}$ and $x_i\in \overleftarrow{W}$. Let $A$ and $C$ be the endpoints of $e$ and $f$, respectively, that lie in $\overrightarrow{W}$, and let $B = x_i$. Let $m = m(A,B,C)$ be the median of $A,B,C$. Note that, by construction, the geodesics $\alpha_i,\alpha_{i+1}$ are contained in standard product regions $P_i,P_{i+1}$, respectively. $P_i$ and $P_{i+1}$ are convex. Since $A,B\in P_i$, $m\in P_i$. Similarly, $m\in P_{i+1}$. Also, since $A,C\in \overrightarrow{W}$ and the halfspace $\overrightarrow{W}$ is convex, $m\in \overrightarrow{W}$. In particular, $m\neq B$ since $B\in \overleftarrow{W}$. By the definition of median, there is a geodesic $\beta_i$ from $A$ to $B$ passing through $m$. Thus, $d(A,m) < d(A,B)$. 
    
    Since $m\in P_i\cap P_{i+1}$, the path $x_0,\cdots,x_{i-1},m,x_{i+1},\cdots,x_n$ is  a geodesic in $\widehat{X}$. Let $\alpha_i'$ be the combinatorial geodesic from $x_{i-1}$ to $m$ inside $P_i$. Clearly, $|\alpha_{i}'| < |\alpha_i|$, hence contradicting the lexicographical minimality of $\alpha$.

    \noindent
    \textbf{Case 2: $j>i+1$.} Choose $e$ and $f$ innermost, so the subpath $\beta$ between them crosses no hyperplane twice and is therefore geodesic. Since each combinatorial hyperplane is convex, $\beta$ is contained in a combinatorial hyperplane of $W$, say $W\times \{\frac{1}{2}\}$.  Sliding $\beta$ across the carrier of $W$ gives a parallel geodesic $\beta'$ in the opposite combinatorial hyperplane $W\times \{-\frac{1}{2}\}$. For $i+1\leq k\leq j-1$, let $x_k'$ be the vertex opposite $x_k$ across $W$. Convexity of the carrier shows that $x_k'$ and $x_{k+1}'$ are adjacent in $\widehat{X}$ exactly when $x_k$ and $x_{k+1}$ are. Replacing the portion of the path between $e$ and $f$ by $\beta'$ preserves the length of the $\widehat{X}$-geodesic but shortens the first affected subpath $\alpha_i$, contradicting lexicographic minimality.

    Therefore, $\alpha$ is a geodesic from $x$ to $y$ in $X$. 
\end{proof}

\begin{rmk}
\label{special taut geodesic}
    By Proposition~\ref{taut geodesic existence}, we may choose a taut geodesic $x = x_0, \cdots, x_n = y$ whose vertices occur, in this order, on a combinatorial geodesic $\gamma$ in $X$. Consequently, each hyperplane of $X$ separates at most one consecutive pair $x_i,x_{i+1}$. 
\end{rmk}

\begin{rmk}
\label{taut geodesic remark}
 Suppose $d_{\widehat{X}} (x,y) = n$. Then by the existence of taut geodesics, there is $x'\in [x,y]$ with $d_{\widehat{X}}(x,x') = n-1$ and $d_{\widehat{X}}(x',y) = 1$. Equivalently, if every $x'\in [x,y] \cap S_1(y)$ is at distance at least $n-1$ from $x$, then $d_{\widehat{X}}(x,y) \geq n$. This will be an essential observation used in Lemma~\ref{C'}.
\end{rmk}

\subsection{Image of cubical convex hulls in $\widehat{X}$}

Now we show that the gate map $\fg$ from $X$ to a convex subcomplex is 1-Lipschitz with respect to the $d_{\widehat{X}}$ metric.

\begin{lem}[Gate maps are $1$-Lipschitz]
\label{Projections don't increase the distance}
    Let $x,y\in X$ with $d_{\widehat{X}}(x,y) = n$. Let $C\subset X$ be a convex subcomplex of $X$. Let $x'=\fg_C(x)$ and $y'=\fg_C(y)$ be the closest point projections. Then $d_{\widehat{X}}(x',y') \leq n$.
\end{lem}

\begin{proof}
    We argue by induction on the distance $n = d_{\widehat{X}} (x,y)$. Suppose $n=1$. By the property of closest point projection, every hyperplane separating $x'$ and $y'$ also separates $x$ and $y$. Thus, $\lab([x',y']) \subset \lab([x,y])$. Then $d_{\widehat{X}}(x',y') \leq 1$.

    Assume $d_{\widehat{X}}(x',y') \leq n-1$ whenever $d_{\widehat{X}}(x,y) \leq n-1$. Now suppose $d_{\widehat{X}}(x,y) = n$. Then there is $z\in [x,y]$ such that $d_{\widehat{X}}(x,z) = n-1$ and $d_{\widehat{X}}(z,y) =1$ by Proposition~\ref{taut geodesic existence}. Let $z' = \fg_C(z)$. By the induction hypothesis, $d_{\widehat{X}}(x',z') \leq n-1$. For the same reason as above, $d_{\widehat{X}}(z',y') \leq 1$. Hence $d_{\widehat{X}}(x',y') \leq d_{\widehat{X}}(x',z') + d_{\widehat{X}}(z',y') \leq n$.  
\end{proof}

The existence of taut geodesics allows us to study $d_{\widehat{X}}(x,y)$ using only the vertices of the cubical convex hull $[x,y] \subset X$.

\begin{lem}[Monotonicity in a convex hull]
\label{distance monotonicity in hull}
    Let $o,x$ be vertices in $\widehat{X}$. For every $y\in [o,x]$, $d_{\widehat{X}}(o,y) \leq d_{\widehat{X}}(o,x)$.
\end{lem}
\begin{proof}
    Since $y \in [o,x]$, we have $[o,y] \subset [o,x]$. Furthermore, $\fg_{[o,y]}(x) = y$. By applying Lemma~\ref{Projections don't increase the distance}, we have
    \[
    d_{\widehat{X}}(o,y) = d_{\widehat{X}}(\fg_{[o,y]}(o) , \fg_{[o,y]}(x))\leq d_{\widehat{X}}(o,x).
    \]

\end{proof}

The next  lemma shows that if two points are close in $\widehat{X}$ then the convex hulls from a fixed point $o$ to each of them have a large intersection. It will be used in the filter construction to control the diameter of the top boundary path in Section~\ref{Fan and Filter Section}.

\begin{lem}[Large intersection of hulls]
\label{large intersection of convex hull}
    Suppose $x,y\in S_R(o)$ and $d_{\widehat{X}}(x,y) = E$ where $E< R$. Then there is $z\in [o,x]\cap [o,y]$ with $d_{\widehat{X}}(o,z) \geq R-E$.

\end{lem}
\begin{proof}
    Let $z:= \fg_{[o,y]}(x)$, the closest point projection. Since $[o,y]\cap [o,x] \neq \emptyset$, by Lemma~\ref{Bridge Lemma}, $z\in [o,x]\cap [o,y]$. Furthermore, by Lemma~\ref{Projections don't increase the distance},
    \[
    d_{\widehat{X}}(z,y) = d_{\widehat{X}}(\fg_{[o,y]}(x),\fg_{[o,y]}(y)) \leq d_{\widehat{X}}(x,y) = E.
    \]
    Hence, $d_{\widehat{X}}(o,z)\geq d_{\widehat{X}}(o,y) - d_{\widehat{X}}(z,y) \geq R-E$.
\end{proof}

The following lemma shows that the image of a convex hull of two vertices in $X$ is uniformly close to a taut geodesic in $\widehat{X}$.

\begin{lem}[Hulls are uniformly close to taut geodesics]
\label{convex hulls are uniformly thin}
    Let $M =\max \{|K| \mid K\subset \Gamma \text{ is a clique}\}$. Let $o,z$ be two 0-cubes in $X$. Then for every $x\in [o,z]$, there is $y$ on a taut geodesic from $o$ to $z$ with $d_{\widehat{X}}(x,y) \leq M+1$.
\end{lem}
\begin{proof}
    Let $M' = \dim([o,z])$, the dimension of the convex hull of $\{o,z\}$. Note that $M'\leq M$.
    Let $o = y_0,\cdots, y_n = z$ be a taut geodesic in $\widehat{X}$ from $o$ to $z$ such that every hyperplane in $X$ separates at most one consecutive pair. For two vertices $a,b$, denote by $\W(a,b)$ the set of hyperplanes separating $a$ and $b$. Then for each $0\leq i \leq n$, define
    \[
    \A_i  = \W(o,x)\setminus \W(o,y_i), \quad
    \B_i = \W(o,y_i) \setminus \W(o,x).
    \]

    \begin{claim*}
        For each $0\leq i \leq n$, $\W(x,y_i) = \A_i \sqcup \B_i$. Moreover, for every hyperplane $V\in \A_i$ and $H\in \B_i$, $V$ crosses $H$.
    \end{claim*}
    \begin{proof}
        Let $W$ be a hyperplane separating $x$ and $y_i$. Since $x,y_i \in [o,z]$, $W$ also separates $o$ and $z$. If $W\notin \A_i$, then it must separate $\{o,x\}$ from $\{z,y_i\}$. Hence, $W\in \B_i$. Also, clearly, $\A_i \cap \B_i = \emptyset$. Thus, $\W(x,y_i) = \A_i \sqcup \B_i$.

        Let $V\in \A_i$ and $H\in \B_i$. Suppose $H$ has two halfspaces $\overleftarrow{H}$ and $\overrightarrow{H}$ such that $\{o,x\} \subset \overleftarrow{H}$ and $\{z,y_i\} \subset \overrightarrow{H}$. Suppose $V$ does not cross $H$. Then $V \subset \overleftarrow{H}$. Hence, $V$ would separate $x$ from $\{o,z\}$. Contradiction.        
    \end{proof}
    Then for each $i$,  $ \lab(\A_i)$ and $\lab(\B_i)$ span a join in the defining graph $\Gamma$. Furthermore, by the choice of taut geodesic, $\A_i \subset \A_{i-1}$ and $\B_{i-1} \subset \B_{i}$ for $i=1,\cdots ,n$.

    Let $k\in \{0,1,\cdots ,n\}$ be the minimal index such that $\lab(\A_k)$ is a clique. Such an index $k$ exists. Indeed, since $x\in[o,z]$, we have $\W(o,x)\subseteq\W(o,z)$, and  $\A_n=\emptyset$. We consider the empty set as a clique.

    Suppose $k = 0$. Then $y_0 = o$ and $x$ are contained in a cube. Then,
    $d_{\widehat{X}}(x,o) \leq d(x,o) \leq \dim([o,z]) = M'$.

    Suppose $k\geq 1$. Then $\lab(\A_k)$ is a clique and $\lab(\A_{k-1})$ is not a clique. 

    \noindent
    \textbf{Case 1:} Suppose $\lab(\B_{k-1})$ is not a clique. Then $\lab(\A_{k-1})\ast \lab(\B_{k - 1})$ is a nontrivial join. Hence, $d_{\widehat{X}}(x,y_{k-1}) \leq 1$. 

    \noindent
    \textbf{Case 2:} Suppose $\lab(\B_{k-1})$ is a clique. Let $x' = \fg_{[y_{k-1},y_k]}(x)$. See Figure~\ref{Hull is uniform close to a taut geodesic picture}.

    \begin{figure}
    \centering
    \includegraphics{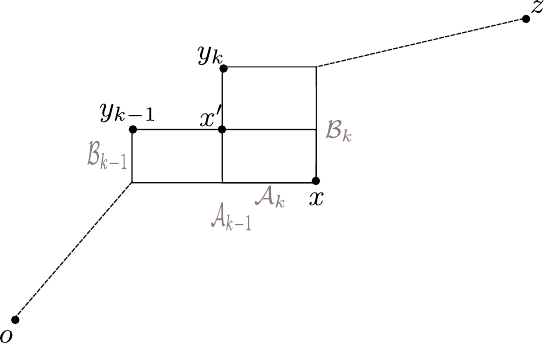}
    \caption{Vertices $x,y_k,y_{k-1}$ in the convex hull of $o$ and $z$.}
    \label{Hull is uniform close to a taut geodesic picture}
    \end{figure}

     \begin{claim*}
        $\W(x,x') = \A_k \sqcup \B_{k-1}$. Furthermore, every hyperplane in $\A_k$ crosses every hyperplane in $\B_{k-1}$.
     \end{claim*}
     
    \begin{proof}
        By the gate-map property, a hyperplane is in $\W(x,x')$ if and only if  it separates $x$ from $\{y_{k-1},y_k\}$.

        \emph{$\A_k \cup \B_{k-1} \subset \W(x,x')$:} Let $V\in\A_k$. Since $\A_k\subseteq\A_{k-1}$, $V$ separates $o$ from $x$, but separates neither $y_{k-1}$ nor $y_k$ from $o$. Thus $V$ separates $x$ from $\{y_{k-1},y_k\}$, and hence $V\in\W(x,x')$.

        \emph{$\W(x,x')\subset \A_k \cup \B_{k-1}$:} Let $W\in \W(x,x')$. Suppose $W$ separates $x$ and $o$. Then it does not separate $o$ and $y_k$; otherwise $W\notin \W(x,x')$. Hence, $W\in \A_k$. Now suppose $W$ does not separate $x$ and $o$. Then it must separate $\{o,x\}$ from $\{y_k,y_{k-1}\}$. Hence, $W\in \B_{k-1}$.

        Therefore, $\W(x,x') = \A_k \cup \B_{k-1}$. Since $\B_{k-1}\subset \B_k$ and $\B_k\cap \A_k = \emptyset$, $\B_{k-1}\cap \A_k = \emptyset$. Hence, $\W(x,x') = \A_k \sqcup \B_{k-1}$. By the previous claim, every hyperplane in $\B_k$ crosses every hyperplane in $\A_k$. Since $\B_{k-1}\subset \B_k$, the furthermore part follows.
    \end{proof}

    Therefore, $\lab([x,x']) = \lab(\A_k) \ast \lab(\B_{k-1})$ is a clique. $d_{\widehat{X}}(x,x')\leq d(x,x') \leq M'$. Since $y_{k-1} $ and $y_k$ are adjacent in $\widehat{X}$, they are either connected by a 1-cube  or contained in the same product region in $X$. Since $x'\in [y_{k-1}, y_k]$, in either case $d_{\widehat{X}} (x',y_k) \leq 1$. Therefore,
    \[
    d_{\widehat{X}}(x,y_k)\leq d_{\widehat{X}}(x,x') + d_{\widehat{X}}(x',y_k)\leq M'+1.
    \]

\end{proof}

\subsection{Wide combinatorial paths in $X$}
\begin{defn}[Wide path]
    We say a combinatorial path $\alpha$ in $X_{\Gamma}$ is \emph{wide} if $\alpha$ is contained in a standard product region. Equivalently, $\alpha$ is wide if and only if $\lab(\alpha)$ is contained in a nontrivial join.
\end{defn}

The next two lemmas relate wide subpaths in $X$ to distance in $\widehat{X}$. Lemma~\ref{Morse Geodesic C'} transfers a bound on wide subpath lengths from one combinatorial geodesic to every combinatorial geodesic with the same endpoints. Lemma~\ref{C'} then converts this bound into a lower bound on distance in $\widehat{X}$. Together, these estimates will control progress in the filter construction in Section~\ref{Fan and Filter Section}.

\begin{lem}[Longest wide subpath]
\label{Morse Geodesic C'}
    Let $\alpha$ be a combinatorial geodesic from $x$ to $y$ such that every wide subpath has length at most $C$ for some constant $C\geq 1$. Then there is a constant $C' = C'(C,\Gamma)$ such that every wide subpath of any combinatorial geodesic from $x$ to $y$ has length at most $C'$.
\end{lem}
\begin{proof}
Let $M$ be the maximal cardinality of a clique in $\Gamma$. We will prove that one may take
\[
    C'=(2|V(\Gamma)|+1)\max\{C,M\}.
\]

Let $\gamma$ be a combinatorial geodesic from $x$ to $y$, and let $\gamma'$ be a wide subpath of $\gamma$. We assume that $\gamma'$ has positive length. There is a nontrivial join $\Lambda \subseteq\Gamma$ containing $\lab(\gamma')$.

Both $\alpha$ and $\gamma$ cross precisely the hyperplanes separating $x$ from $y$, each exactly once. Thus the hyperplanes crossed by $\gamma'$ provide a collection of 1-cubes of $\alpha$. Let $\alpha_1,\ldots,\alpha_n$ be the maximal nonempty subpaths of $\alpha$ formed by concatenating these 1-cubes. Without loss of generality, we assume $\alpha_i$ appears before $\alpha_{i+1}$ when traveling from $x$ to $y$ along $\alpha$. For $1\leq t<n$, let $\beta_t$ be the nonempty subpath
between $\alpha_t$ and $\alpha_{t+1}$. Denote
\[
    A_i=\lab(\alpha_i),
    \qquad
    B_t=\lab(\beta_t).
\]
Clearly, for each $i=1,\cdots,n$, $A_i\subset \Lambda$ and $\alpha_i$ is a wide subpath. Moreover, by the construction, we have
\[
    |\gamma'|=\sum_{i=1}^{n}|\alpha_i|.
\]

\medskip

Let $\mathcal B$ be the set of hyperplanes crossed by $\gamma$ but not crossed by $\gamma'$. Since $\gamma'$ is a subpath of $\gamma$, each hyperplane in $\mathcal B$ separates $\gamma'$ from either $x$ or $y$.
We define
\[
    \mathcal B^x = \{W\in\mathcal B\mid W\text{ separates }x\text{ from }\gamma'\},
\quad
    \mathcal B^y = \{W\in\mathcal B \mid W\text{ separates }y\text{ from }\gamma'\}.
\]
Hence, $\mathcal{B} = \mathcal{B}^x\sqcup \mathcal{B}^y$.
For $1\leq t<n$, let
\[
    B_t^x
    =
    \lab\{W\in\mathcal B^x:W\text{ crosses }\beta_t\},
    \qquad
    B_t^y
    =
    \lab\{W\in\mathcal B^y:W\text{ crosses }\beta_t\}.
\]
Thus $B_t=B_t^x\cup B_t^y$. See Figure~\ref{longest wide subpath picture}.

    \begin{figure}
    \centering
    \includegraphics[scale=0.8]{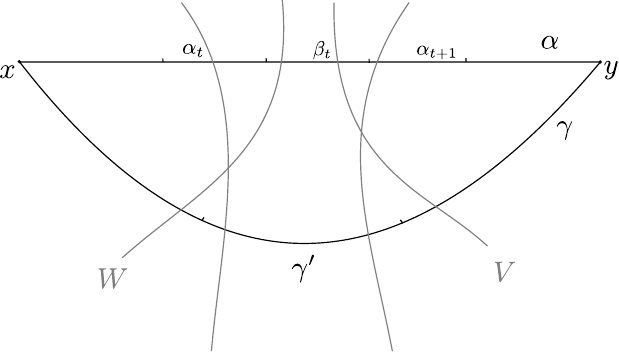}
    \caption{$\lab(W) \in B^x_t$ and $\lab(V) \in B^y_t$.}
    \label{longest wide subpath picture}
    \end{figure}

Suppose that $W\in\mathcal B^x$ crosses $\beta_t$. For every $i\leq t$, the subpath $\alpha_i$ and $\gamma'$ are separated by $W$. Since by construction every hyperplane crossing $\alpha_i$ also crosses $\gamma'$, it must cross $W$ as well.
Hence
\[
    A_i\subseteq\lk(B_t^x) \qquad(\text{for }i\leq t).
\]
Similarly,
\[
    A_i\subseteq\lk(B_t^y) \qquad(\text{for }i>t).
\]

For $1\leq i\leq n$, let 
\[
    X_i=\bigcup_{t=i}^{n-1}B_t^x, \qquad Y_i=\bigcup_{t=1}^{i-1}B_t^y.
\]

The preceding observation gives
\[
    A_i\subseteq\lk(X_i)\cap\lk(Y_i).
\]
Furthermore, by construction we have
\[
    X_1\supseteq X_2\supseteq\cdots\supseteq X_n, \qquad
    Y_1\subseteq Y_2\subseteq\cdots\subseteq Y_n.
\]

\medskip

Partition $\{1,\ldots,n\}$ into maximal consecutive intervals $I_1,\ldots,I_m$ on which the pair $(X_i,Y_i)$ is constant. Each change of this pair removes at least one vertex from $X_i$ or adds at least one vertex to $Y_i$. Since each vertex can be removed once and added once, there are at most $2|V(\Gamma)|$ changes. Thus
\[
    m\leq 2|V(\Gamma)|+1.
\]

Fix an interval $I_\ell=\{p,\cdots,q\}$, and let
\[
    K=\bigcup_{j=p}^{q}A_j, \qquad T=\bigcup_{t=p}^{q-1}B_t.
\]
Let $\tau_\ell$ be the subpath of $\alpha$ beginning with $\alpha_p$ and ending with $\alpha_q$, i.e.
\[
    \tau_\ell = \alpha_p\beta_p\alpha_{p+1} \cdots \beta_{q-1}\alpha_q.
\]
When $p=q$, we just have $\tau_\ell=\alpha_p$. In particular,
\[
    \lab(\tau_\ell)=K\cup T.
\]

We claim that every element of $T$ commutes with every element of $K$. Indeed, for any $b\in T$, there is $p\leq t \leq q-1$ such that either $b\in B^x_t$ or $b\in B^y_t$.  If $b\in B^x_t$, since $X_t = X_q$ by constancy on $I_\ell$, we have $K\subset \lk(b)$. Similarly,  if $b\in B^y_t$, since $Y_{t+1} = Y_p$ by constancy on $I_\ell$, we have $K\subset \lk(b)$. Thus, $K\cup T$ spans a join in $\Gamma$.

Now we claim that
\[
    \sum_{j=p}^{q}|\alpha_j|\leq\max\{C,M\}.
\]

\medskip
\noindent
\textbf{Case 1}: $K$ is not a clique.
Let
\[
    \Delta=\lk(T)\cap \Lambda.
\]
Since $K\subseteq\Delta$, the subgraph $\Delta$ is not a clique. Lemma~\ref{in the same join} implies that $T\cup\Delta$ is contained in a nontrivial join. Hence,
\[
    \lab(\tau_\ell) =K\cup T \subseteq \Delta\cup T
\]
is contained in a nontrivial join. Thus $\tau_\ell$ is a wide subpath of $\alpha$, and
\[
    \sum_{j=p}^{q}|\alpha_j| \leq|\tau_\ell| \leq C.
\]

\medskip
\noindent
\textbf{Case 2}: $K$ is a clique.
Since every element of $T$  commutes with every element of $K$, there is a cube spanned by (possibly translates of) $\alpha_p,\cdots,\alpha_q$ in the convex hull $[x,y]$. Hence,
\[
    \sum_{j=p}^{q}|\alpha_j| \leq |K| \leq M.
\]

\medskip
Therefore,

\[
    |\gamma'|
    = \sum_{\ell=1}^{m} \sum_{j\in I_\ell}|\alpha_j|
    \leq m\max\{C,M\}
    \leq (2|V(\Gamma)|+1)\max\{C,M\}.
\]
\end{proof}

\begin{lem}[Bounded wide subpaths in $X$ force progress in $\widehat{X}$]
\label{C'}
    Let $x,y$ be two 0-cubes in $X_\Gamma$ with the property that every wide subpath of any combinatorial geodesic from $x$ to $y$ has length at most $C'$ for some constant $C' \geq 1$. For each $n\in \N$, if $d(x,y)\geq n(C'+1)$, then $d_{\widehat{X}}(x,y)\geq n$.
\end{lem}
\begin{proof}
    We argue by induction on $n$. Suppose $d(x,y)\geq C'+1$. Then $x\neq y$. Hence, $d_{\widehat{X}}(x,y) \geq 1$.

    Assume the property holds for $n-1$. Now we assume that $d(x,y) \geq n(C'+1)$. Let $z\in [x,y]$ with $d_{\widehat{X}}(y,z)\leq 1$. Then $d(y,z)\leq C'$. Hence $d(x,z) \geq (n-1)(C'+1) +1$. Since $z\in [x,y]$, every combinatorial geodesic from $x$ to $z$ extends to a combinatorial geodesic from $x$ to $y$. Consequently, each of its wide subpaths has length at most $C'$. Then by the induction hypothesis, $d_{\widehat{X}}(x,z)\geq n-1$.  By the arbitrariness of $z$ and Remark~\ref{taut geodesic remark}, $d_{\widehat{X}}(x,y)\geq n$.
\end{proof}

%% file: sections/dead-ends.tex
Birman and Menasco showed that the curve graph of a closed orientable surface of genus at least two contains dead ends \cite{birman2015curve}. In particular, they constructed dead ends at arbitrarily large distance from a fixed vertex. Since the curve graph is a basic model for hyperbolic spaces arising in hierarchical structures, it is natural to ask whether maximized hyperbolic spaces of right-angled Coxeter groups exhibit a similar phenomenon.

A dead end results in the failure of exact geodesic extension. A geodesic segment ending at a dead end cannot be extended further as a geodesic from a basepoint. Thus, exact geodesic extension cannot be assumed in $\widehat X$. The examples below illustrate this obstruction and motivate the coarse extension results proved in the following section.

We first recall the definition of dead ends and then give two examples. The first produces dead ends on the sphere with radius one and pairs of vertices joined by infinitely many geodesics. The second produces dead ends arbitrarily far from a fixed basepoint. Although this resembles the phenomenon in the curve graph, the underlying technique is different: in the curve graph, dead ends arise from separating curves, whereas in our examples they arise from the arrangement of maximal nontrivial joins in the defining graph.

\begin{defn}[Dead ends]
\label{dead ends definition}
Let $\mathcal G$ be a simplicial graph equipped with its path metric $d$, and fix a vertex $o\in\mathcal G$. A vertex $x\in\mathcal G$ is a \emph{dead end} with respect to $o$ if
\[
d(o,z)\leq d(o,x)
\]
for every vertex $z$ adjacent to $x$.
\end{defn}

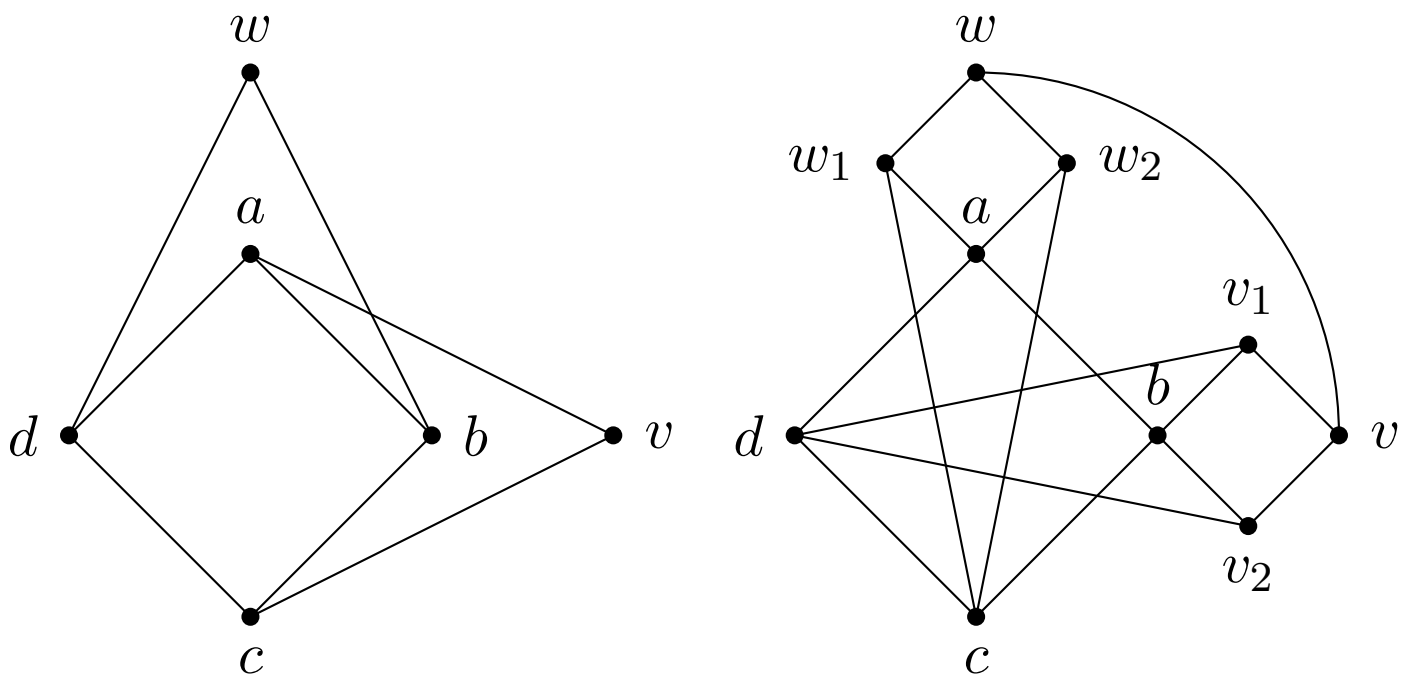
\begin{figure}
    \centering
    \begin{tikzpicture}[scale = 0.6]
    \tikzstyle{vertex} = [circle, fill, inner sep=1.2pt]

    \draw (-5,2) -- (-7,0);
    \draw (-5,2) -- (-3,0);
    \draw (-5,2) -- (-1,0);

    \draw (-5,-2) -- (-7,0);
    \draw (-5,-2) -- (-3,0);
    \draw (-5,-2) -- (-1,0);

    \draw (-5,4) -- (-7,0);
    \draw (-5,4) -- (-3,0);

    \node [vertex] at (-5,2) {};
    \node [above=2pt] at (-5,2) {$a$};

    \node [vertex] at (-3,0) {};
    \node [right=2pt] at (-3,0) {$b$};

    \node [vertex] at (-5,-2) {};
    \node [below=2pt] at (-5,-2) {$c$};

    \node [vertex] at (-7,0) {};
    \node [left=2pt] at (-7,0) {$d$};

    \node [vertex] at (-1,0) {};
    \node [right=2pt] at (-1,0) {$v$};

    \node [vertex] at (-5,4) {};
    \node [above=2pt] at (-5,4) {$w$};


    \draw (2,3) -- (3,4);
    \draw (2,3) -- (3,2);
    \draw (2,3) -- (3,-2);

    \draw (4,3) -- (3,4);
    \draw (4,3) -- (3,2);
    \draw (4,3) -- (3,-2);

    \draw (6,1) -- (1,0);
    \draw (6,1) -- (5,0);
    \draw (6,1) -- (7,0);

    \draw (6,-1) -- (1,0);
    \draw (6,-1) -- (5,0);
    \draw (6,-1) -- (7,0);

    \draw (3,2) -- (1,0);
    \draw (3,2) -- (5,0);

    \draw (3,-2) -- (1,0);
    \draw (3,-2) -- (5,0);

    \draw (3,4) to[out = 0,in = 90] (7,0);

    \node [vertex] at (3,2) {};
    \node [above=2pt] at (3,2) {$a$};

    \node [vertex] at (5,0) {};
    \node [above=2pt] at (5,0) {$b$};

    \node [vertex] at (3,-2) {};
    \node [below=2pt] at (3,-2) {$c$};

    \node [vertex] at (1,0) {};
    \node [left=2pt] at (1,0) {$d$};

    \node [vertex] at (7,0) {};
    \node [right=2pt] at (7,0) {$v$};

    \node [vertex] at (3,4) {};
    \node [above=2pt] at (3,4) {$w$};

    \node [vertex] at (6,1) {};
    \node [above=2pt] at (6,1) {$v_1$};

    \node [vertex] at (6,-1) {};
    \node [below=2pt] at (6,-1) {$v_2$};

    \node [vertex] at (2,3) {};
    \node [left=2pt] at (2,3) {$w_1$};

    \node [vertex] at (4,3) {};
    \node [right=2pt] at (4,3) {$w_2$};
\end{tikzpicture}

   \caption{Defining graphs $\Gamma_1$ (left) and $\Gamma_2$ (right) of RACGs
   whose maximized hyperbolic spaces contain dead ends.}
    \label{Dead Ends Graphs}
\end{figure}

The following example produces dead ends on the sphere with radius one about any chosen basepoint. It also gives pairs of vertices joined by infinitely many geodesics in $\widehat X$.
 
\begin{exa}
\label{exist infinitely many geodesics}
    Let $X=X_{\Gamma_1}$, where $\Gamma_1$ is the left graph in
    Figure~\ref{Dead Ends Graphs}. It has two maximal nontrivial joins:
    \[
    \Lambda_1= \{a,c\}\ast \{b,d,v\}, \quad \Lambda_2 = \{b,d\} \ast \{a,c,w\}.
    \]
    Their intersection $\Delta = \{a,c\}\ast \{b,d\}$ is itself a nontrivial join. Because $v$ and $w$ do not lie in a common nontrivial join, $\widehat{X}$ has infinite diameter.

    \noindent\emph{Dead ends on the sphere with radius 1.} Fix $o\in \widehat{X}$ and take $x\neq o\in \widehat{X}$ such that $\lab([o,x]) \subset \Delta$. Then $d_{\widehat{X}}(o,x) = 1$. Let $z\in \widehat{X}$ be adjacent to $x$. Then $z$ and $x$ are either connected by a 1-cube or contained in the same standard product region in $X$. 

    \noindent
    \textbf{Case 1}: Suppose $z$ and $x$ are adjacent by a 1-cube labeled by $u\in V(\Gamma_1)$. Since for any $u\in V(\Gamma_1)$, $u$ and $\Delta$ are contained in a nontrivial join, $\lab([o,z]) \subset \Delta \cup \{u\}$ is contained in a nontrivial join. Hence $d_{\widehat{X}}(o,z) \leq  1$. 

    \noindent
    \textbf{Case 2}: Otherwise, suppose $d(x,z) \geq 2$ and $\lab([x,z])$ is contained in a nontrivial join. Then $\lab([x,z]) \subseteq \Lambda_i$ for $i=1$ or $2$. Since $\Delta = \lab([o,x]) \subseteq \Lambda_i$ for $i=1,2$ as well, $\lab([o,z])\subseteq \Lambda_i$. Then $d_{\widehat{X}}(o,z) \leq  1$. 

    Therefore, $x$ is a dead end with respect to $o$.
    
    \smallskip
    \noindent\emph{Infinitely many geodesics between a pair of points in $\widehat{X}_{\Gamma_1}$.} We use this example to show the existence of infinitely many geodesics between a pair of points in $\widehat{X}_{\Gamma_1}$ as well. Let $p,q$ be two 0-cubes in $X_{\Gamma_1}$ which are separated exactly by two hyperplanes labeled by $v$ and $w$. Since $v,w$ are not in a nontrivial join in $\Gamma_1$, $d_{\widehat{X}}(p,q) = 2$. Let $z$ be the 0-cube which is adjacent to both $p$ and $q$ in $X_{\Gamma_1}$. Since the convex hull $[p,q]\subset X_{\Gamma_1}$ is the combinatorial geodesic from $p$ to $q$, the path $p,z,q$ in $\widehat{X}_{\Gamma_1}$ is the unique taut geodesic from $p$ to $q$. For every vertex $z' \neq z$ contained in the copy of the flat $X_{\Delta}$ passing through $z$, the path $p,z',q$ in $\widehat{X}_{\Gamma_1}$ is also a geodesic but not taut.
\end{exa}

The next example contains dead ends arbitrarily far from the basepoint.

\begin{exa}
    Let $X_{\Gamma_2}$ be the Davis complex of the RACG $W_{\Gamma_2}$ whose defining graph appears in Figure~\ref{Dead Ends Graphs}. The graph $\Gamma_2$ has four maximal nontrivial joins: 
    \[
    \begin{aligned}
    \Lambda_1 &= \{v_1,v_2\}*\{b,d,v\},&
    \Lambda_2 &= \{w_1,w_2\}*\{a,c,w\},\\
    \Delta_1 &= \{a,c\}*\{b,d,w_1,w_2\},&
    \Delta_2 &= \{b,d\}*\{a,c,v_1,v_2\}.
    \end{aligned}
    \]
    Fix $o\in \widehat{X}_{\Gamma_2}$. Let $\gamma$ be the combinatorial geodesic in $X_{\Gamma_2}$ starting at $o$ labeled by the word $\overline{\gamma} = vwacbd$. Let $x$ be the terminal point of $\gamma$. Applying Proposition~\ref{taut geodesic existence}, we know that $d_{\widehat{X}}(o,x) = 3$ by enumerating paths in $\widehat{X}_{\Gamma_2}$ with vertices in $[o,x]$.  There are two points $y_1,y_2\in [o,x]$ such that $y_1$ and $x$ are separated by two hyperplanes of types $b$ and $d$, and $y_2$ and $x$ are separated by two hyperplanes of types $a$ and $c$. Note that $y_1,y_2\in S_2(o)$.

    We claim that $x$ is a dead end with respect to $o$. Let $z$ be a vertex in $\widehat{X}_{\Gamma_2}$ such that $d_{\widehat{X}}(x,z) = 1$ and $z\notin [o,x]$. Let $\Pi=\lab([x,z])$. The set $\Pi$ is either a singleton or is contained in a nontrivial join; hence it lies in one of $\Lambda_1,\Lambda_2,\Delta_1,\Delta_2$. Consequently, $\Pi$ lies in a common nontrivial join with either $\{b,d\}$ or $\{a,c\}$. If $\Pi$ lies in a common nontrivial join with $\{b,d\}$, then $\lab([z,y_1])\subseteq \Pi\cup\{b,d\}$ is contained in a nontrivial join. Hence $d_{\widehat X}(z,y_1)\leq 1$. Similarly, if  $\Pi$ lies in a common nontrivial join with $\{a,c\}$, then $d_{\widehat{X}}(z,y_2)\leq 1$. Therefore,
    \[
    d_{\widehat{X}}(o,z)\leq3=d_{\widehat{X}}(o,x).
    \]
    Neighbors of $x$ in $[o,x]$ are not farther from $o$, so $x$ is a dead end with respect to $o$.

    Furthermore, for each $n\in \N$, let $x_n$ be the terminal point of the geodesic $\gamma_n$ labeled by the word $(vwacbd)^n$ with initial point $o$. By this construction, $d_{\widehat{X}}(o,x_n) = 2n+1$. (The existence of taut geodesics enables us to focus on combinatorial geodesics in the convex hull when computing the distance in $\widehat{X}_{\Gamma_2}$.) By the same analysis, $x_n$ is a dead end with respect to $o$ for every $n\in \N$.    
\end{exa}

These examples highlight a distinction from the curve graph $\CCS$ of a finite-type surface. Following \cite[Definition~4.2]{wright2024spheres}, a simple closed curve is \emph{essentially non-separating} if it is non-separating or bounds a pair of pants. The bounded geodesic image theorem \cite[Theorem~3.1]{masur1998geometry} implies that such a curve is never a dead end, regardless of the base vertex (see also \cite[Corollary~4.3]{wright2024spheres}). We do not have an analogous characterization for $\widehat{X}$. Proposition~\ref{Bounded Geodesic Image} only forces a geodesic in $\widehat{X}$ to pass uniformly \emph{close} to a fixed point, not \emph{through} that point.

%% file: sections/BGI.tex
The examples in the preceding section show that exact geodesic extension may fail in $\widehat X$. Moreover, since $\widehat X$ need not be proper, the usual Arzelà--Ascoli argument cannot be used to extract limiting geodesic rays. In this section, the bounded geodesic image theorem and the HHS distance formula provide a substitute for the missing compactness.

We first give a proof of the bounded geodesic image statement for $\widehat X$ (Proposition~\ref{Bounded Geodesic Image}). Together with the distance formula and the local finiteness of the Davis complex $X$, it allows us to find several properties of $\partial \widehat{X}$ which are familiar from boundaries of hyperbolic groups. We use them to prove that every boundary point in $\partial \widehat{X}$ can be represented by a geodesic ray in $\widehat{X}$ (Lemma~\ref{Visible}) and that every pair of distinct boundary points can be joined by a bi-infinite geodesic (Proposition~\ref{Visibility of boundary}). These visibility results imply that geodesic ideal triangles in $\widehat{X}$ are uniformly slim and that $\widehat X$ is geodesic almost extendable~(Proposition~\ref{exist close geodesic ray}). Finally, when $\partial\widehat{X}_\Gamma$ is nonempty, we use bounded geodesic image once more to show that $\partial\widehat{X}_{\Gamma}$ is noncompact exactly when $\Gamma$ contains an induced 4-cycle~(Theorem~\ref{noncompact}).

\subsection{Bounded geodesic image}

The bounded geodesic image property for maximized hyperbolic spaces of HHSs with unbounded product regions follows from \cite[Theorem~3.7]{abbott2021largest}. For completeness, we give a direct proof in the case of right-angled Coxeter groups.

Recall that $\fg_K:X\to K$ denotes the gate map to a convex subcomplex $K\subset X$. Let $U\in\fS$ be a domain represented by a convex subcomplex $K\subset X$. If
\[
d_U(x,y)>\dim X,
\]
then $\fg_K(x)$ and $\fg_K(y)$ cannot lie in a common cube (see Lemma~\ref{Bridge Lemma}). Thus, sufficiently large projection distance under the gate map gives the hypothesis in the following proposition.

\begin{prop}[Bounded geodesic image]
\label{Bounded Geodesic Image}
    There is a constant $C = C(\Gamma)$ satisfying the following property. Let $x,y$ be any two 0-cubes in $X$. Suppose there is a domain $U \in \fS$ represented by a convex subcomplex $K = gX_{\Lambda}$ for some $g\in W_\Gamma$ and $\Lambda \subset \Gamma$ such that $\fg_K(x),\fg_K(y)$ are not contained in a cube. Let $\rho^U_{\widehat{X}}\in gX_{\Lambda}\subset\widehat{X}$ be the fixed point in $\widehat{X}$ for the domain $U$. Then for any geodesic $\alpha \subset \widehat{X}$ from $x$ to $y$, we have 
    \[
    \alpha \cap N_{C}(\rho^U_{\widehat{X}}) \neq \emptyset.
    \]
\end{prop}

The proof has two steps. First, we construct a standard product region $P$ containing $[x_K,y_K]$ and show that the gate $\fg_P(x)$ belongs to the cubical interval $[x,y]$. Second, we show that this gate is uniformly close to $\rho^U_{\widehat X}$ and transfer this bound from the cubical interval to an arbitrary geodesic in $\widehat X$.

\begin{proof}
By the construction of the maximized HHS structure, neither $\Lambda$ nor $\lk(\Lambda)$ is a clique. Put $x_K=\fg_K(x)$ and $y_K=\fg_K(y)$. Every hyperplane separating $x_K$ from $y_K$ also separates $x$ from $y$. Moreover,
\[
\Lambda'=\lab([x_K,y_K])\subset\Lambda
\]
is not a clique, because $x_K$ and $y_K$ do not lie in a common cube by assumption. Let $\Delta=\lk(\Lambda')$. The graph $\Delta$ is not a clique, since it contains $\lk(\Lambda)$, and $\Lambda'\subset\lk(\Delta)$. Thus
$\Lambda'\cup\Delta$ lies in a nontrivial join. Let $P$ be a standard product
region containing $[x_K,y_K]$ and associated to a maximal join containing $\Lambda' \ast \Delta$. See Figure~\ref{BGI Projection Figure}.

    \begin{figure}[h]
    \centering
    \includegraphics[width=0.6\textwidth]{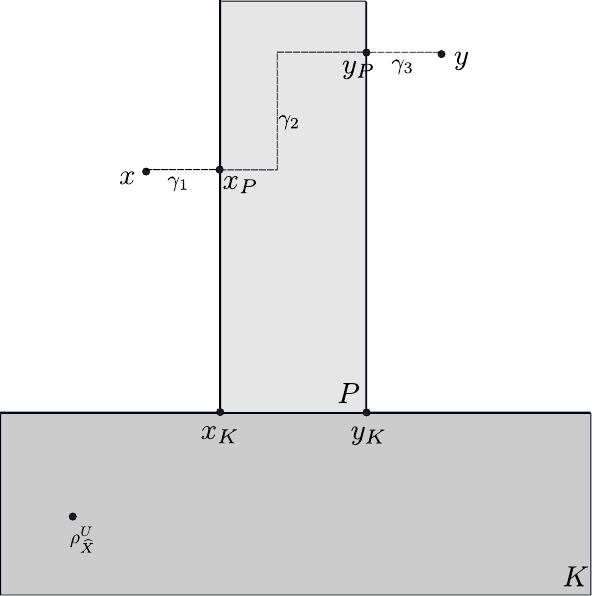}
    \caption{
    The standard product region $P$ and the subcomplex $K$ are in gray. $x_K$ (resp. $y_K$) is the image of $x$ (resp. $y$) under the gate map $\fg_K$. Every geodesic from $x$ to $y$ in $\widehat{X}$ passes within a uniform distance of $\rho^U_{\widehat{X}}$.
    }
    \label{BGI Projection Figure}
    \end{figure}

 Let $x_P=\fg_P(x)$ and $y_P=\fg_P(y)$.  Let $\gamma_1$ be a combinatorial geodesic from $x$ to $x_P$, $\gamma_2$ be a combinatorial geodesic from $x_P$ to $y_P$ and $\gamma_3$ be a combinatorial geodesic from $y_P$ to $y$. Since the product region $P$ is convex, $\gamma_2\subset P$. We claim that the concatenated combinatorial path $\gamma = \gamma_1 \gamma_2 \gamma_3$ is a geodesic from $x$ to $y$. Since each $\gamma_i$ is geodesic for $i=1,2,3$, it is enough to show that no hyperplane crosses both $\gamma_i$ and $\gamma_j$ for $i\neq j$. By the bridge lemma (see for example Lemma~\ref{Bridge Lemma} (4)), any hyperplane  crossing $\gamma_1$ (resp. $\gamma_3$) separates $P$ from $x$ (resp. $y$). Hence it does not cross $\gamma_2$. Now suppose there is a hyperplane crossing $\gamma_1$ and $\gamma_3$. Equivalently, there is a hyperplane separating $x$ and $y$ from $P$. Let $W$ be such a hyperplane which contacts $P$. By the bridge lemma again, 
\[
\lab(W) \in \lk(\lab([x_P,y_P])) \subset \lk(\Lambda') = \Delta.
\]
Hence $W\cap P \neq \emptyset$. Contradiction. Hence $\gamma$ is a combinatorial geodesic from $x$ to $y$ that has nonempty intersection with $P$.

Since $d_{\widehat{X}}(x_K,\rho^U_{\widehat{X}})\leq1$ and $x_K$ and $x_P$ lie in the same standard product region, we have
\[
d_{\widehat{X}}(x_P,\rho^U_{\widehat{X}})\leq2.
\]
Lemma~\ref{convex hulls are uniformly thin} provides a point of a taut geodesic from $x$ to $y$ in $\widehat{X}$ within $M+1$ of $x_P$, where $M$ is the maximal size of a clique in $\Gamma$. Any geodesic in $\widehat{X}$ from $x$ to $y$ fellow travels that taut geodesic within a constant depending only on the hyperbolicity constant $\delta$. Combining these bounds gives the required constant $C=C(\Gamma)$.

\end{proof}

\subsection[Visibility of geodesic rays]{Visibility of geodesic rays in $\widehat{X}$}
\label{Visibility of Geodesic Rays subsection}

In a proper hyperbolic space, any two distinct boundary points are joined by a bi-infinite geodesic. This need not hold in a general non-proper hyperbolic space because the Arzelà--Ascoli theorem is not applicable. Nevertheless, in this subsection, we will show that the points in $\partial\widehat{X}$ are still \emph{visible}, that is, they can be joined by geodesics. Our argument is inspired by Minsky's proof in curve graphs~\cite[Lemma~5.14]{minsky2010classification}.

First, we prove a technical lemma showing that the combinatorial geodesic in $X$ constructed from a taut geodesic in $\widehat{X}$ in Proposition~\ref{taut geodesic existence} is a \emph{hierarchy path}. Let $(X,\fS)$ be a maximized HHS structure on $W_\Gamma$ in the sense of \cite{abbott2021largest}. For each domain $U\in \fS$, the associated hyperbolic space is either a point corresponding to a dummy domain, the nested-maximal space $\widehat{X}$, or a coned-off contact graph, denoted by $\widehat{\CC U}$, with a coarse Lipschitz projection $\pi_U: X \to \widehat{\CC U}$; see Subsection~\ref{HHS subsection}. We call $\pi_U$ a coarse projection because the image $\pi_U(x)$ is a subset of $\widehat{\CC U}$. But $\diam_{\widehat{\CC U}} \pi_U(x) \leq 1$, so $\pi_U$ is coarsely well-defined. We recall the definition of the \emph{hierarchy path} in HHSs.

\begin{defn}[Hierarchy path]
\label{hierarchy path definition}
    For each $D\geq 1$, a path $\gamma$ in $X$ is a \emph{$D$-hierarchy path} if
    \begin{enumerate}
        \item $\gamma$ is a $(D,D)$-quasi-geodesic;
        \item for each domain $U\in \fS$, $\pi_U \circ \gamma$ is an unparametrized $(D,D)-$quasi-geodesic.
    \end{enumerate}
\end{defn}

\begin{lem}[Hierarchy paths for taut geodesics]
\label{hierarch paths for taut geodesics}
    There exists a constant $1\leq D = D(X,\fS)$ with the following property. Let $(x_0,\cdots,x_n)$ be a geodesic in $\widehat{X}$ and suppose that $\gamma$ is a combinatorial geodesic in $X$ from $x_0$ to $x_n$ such that vertices $x_0,\cdots,x_n$ lie on $\gamma$ in the same order. Then $\gamma$ is a $D-$hierarchy path. In particular, combinatorial geodesics constructed in Proposition~\ref{taut geodesic existence} for taut geodesics in $\widehat{X}$ are hierarchy paths.
\end{lem}

\begin{proof}
It is enough to show that for each $U\in \fS$, the image $\pi_U(\gamma)$ is a uniform unparametrized quasi-geodesic.

\medskip
\noindent
\emph{On the $\sqsubseteq$-maximal domain $\widehat{X}$.}
Let $\gamma_i$ be a subpath of $\gamma$ from $x_{i-1}$ to $x_i$. Then $\gamma = \gamma_1\cdots \gamma_n$. Since $\gamma$ is a combinatorial geodesic in $X$, so is $\gamma_i$ for each $i$. Each $\gamma_i$ is either a 1-cube or a geodesic contained in a standard product region (by convexity). Hence,
\[
    \operatorname{diam}_{\widehat X}
        \bigl(\gamma_i \bigr)\leq 1.
\]
Then $\pi_{\widehat{X}}(\gamma)$ is an unparametrized quasi-geodesic.

\medskip
\noindent
\emph{On a proper domain $U\in \fS$.}
Since the hyperbolic space associated with a dummy domain is a point, we may assume $U$ is not a dummy domain. Then the associated hyperbolic space $\widehat{\CC U}$ is the coned-off contact graph. Suppose $U$ is represented by a convex subcomplex $K \subset X$. The coarse projection $\pi_U: X \to \widehat{\CC U}$ is the following composition:
\[
X \xrightarrow{\fg_K} K \xrightarrow{c_U} \CC U \xrightarrow{i} \widehat{\CC U},
\]
where $\fg_K$ is the gate map, $c_U$ sends each vertex to the set of hyperplanes whose carriers contain the vertex, and $i$ is the inclusion map.

We first establish uniform control on the projection of every subpath of $\gamma$. Let $a,z,b$ be vertices on $\gamma$ such that $z\in [a,b]$, where $[a,b]$ is the combinatorial convex hull in $X$ containing $a$ and $b$. Denote
\[
    a_K=\mathfrak g_K(a),\qquad
    z_K=\mathfrak g_K(z),\qquad
    b_K=\mathfrak g_K(b).
\]
Since $z\in[a,b]$, we have
\[
    z_K\in[a_K,b_K].
\]

Choose $H_a\in\pi_U(a)$ and $H_b\in\pi_U(b)$, and let
\[
    \sigma=(H_a = H_0,\ldots,H_m = H_b)
\]
be a geodesic in $\widehat{\CC U}$ from $H_a$ to $H_b$. 
\begin{claim*}
     $\pi_U(z)\subseteq N_2(\sigma)$.
\end{claim*} 
\begin{proof}
If $z_K = a_K$, then it is clear that $\pi_U(z) \subset N_2(H_a)\subset N_2(\sigma)$. Now let $W$ be the hyperplane dual to the last edge of a combinatorial geodesic from $a_K$ to $z_K$. Since $z_K\in[a_K,b_K]$, the hyperplane $W$ separates $a_K$ from $b_K$. Now we construct a path $\alpha$ from $a_K$ to $b_K$ as follows. For each pair of consecutive hyperplanes $H_{i-1}$ and $H_i$ in $\sigma$, we choose a vertex $y_i\in \nn(H_{i-1})\cap \nn(H_i) \cap K$ when they contact, or a vertex $y_i\in \nn(H_i)$ when $H_{i-1}$ and $H_i$ do not contact and they both cross a proper convex subcomplex $F_i\subset K$  whose contact graph is coned off in $\widehat{\CC U}$. Also let $y_0 = a_K$. Then we choose paths $\alpha_i$ from $y_{i-1}$ to $y_i$ contained in $\nn(H_{i-1})$ when $H_{i-1}$ and $H_i$ contact or in $\nn(H_{i-1}) \cup F_i$ for a proper subcomplex $F_i$ when $H_{i-1}$ and $H_i$ do not contact. Such paths exist since $\nn(H_{i-1})$ and $F_i$ are convex. Let $\alpha = \alpha_1\cdots \alpha_m$ be the concatenation. Then $\alpha$ is a combinatorial path from $a_K$ to $b_K$ and contained in $K$ by construction. Since $W$ is a hyperplane separating $a_K$ and $b_K$, $\alpha$ crosses $W$. Therefore, $W$ either contacts $H_i$  or crosses $F_i$ for some $i\in \{0,\cdots, m\}$. In either case, $d_{\widehat{\CC U}}(W, H_i) \leq 1$. Since every hyperplane incident to $z_K$ contacts $W$ at least at $z_K$, this proves the claim.
\end{proof}

Now for an integer $i$, choose
\[
    p_i\in\pi_U(\gamma(i)).
\]
Applying the preceding claim to every subpath gives, whenever $r\leq s\leq t$,
\begin{equation}
\label{hierarchy lift bounded backtracking}
    d_{\widehat{\CC U}}(p_r,p_s) +d_{\widehat{\CC U}}(p_s,p_t)
    \leq d_{\widehat{\CC U}}(p_r,p_t)+4.
\end{equation}
Moreover, consecutive points $p_t,p_{t+1}$ are at distance at most $2$.

Therefore, $(p_t)$ is a uniform unparametrized quasi-geodesic. Indeed, fix $R>4$ and select points successively, starting with $p_0$,  taking the next point at distance at least $R$ from the previously selected point. Include the terminal point if it has not been selected, and denote the resulting sequence by $q_0,\dots,q_{\ell}$. Then consecutive selected points are at distance between $R$ and $R+2$,  except the last pair which are at distance at most $R+2$. Applying \eqref{hierarchy lift bounded backtracking} successively gives, for $0\leq i,j \leq \ell$,
\[
    (R-4)|i-j| - R \leq d_{\widehat{\CC U}}(q_i,q_j) \leq (R+2)|i-j|.
\]
Thus the selected sequence forms a uniform quasi-geodesic. Since each $p_t$ which is not contained in the sequence $q_0,\cdots, q_\ell$ lies in $R-$neighborhood of some $q_{t'}$,  the projection $(p_t)$ is an unparametrized quasi-geodesic with constants independent of $U$ or $\gamma$.

\end{proof}

Now, by using taut geodesics and hierarchy paths, we construct a geodesic ray in $\widehat{X}$ from a basepoint to a boundary point.

\begin{lem}
\label{Visible}
    For each $o\in \widehat{X}$ and $\nu \in \partial \widehat{X}$, there is a geodesic ray $\mu:[0,\infty)\to \widehat{X}$ with $\mu(0) = o$ and $\mu(\infty) = \nu$. 
\end{lem}
\begin{proof}

    Let $(x_0,x_1,\cdots,x_n,\cdots)$ be a sequence of points in $\widehat{X}$ representing the boundary point $\nu \in\partial \widehat{X}$. For each $i$, take a  taut geodesic $c_i$ from $o$ to $x_i$, together with the ordered combinatorial geodesic lift $\alpha_i$ provided in Proposition~\ref{taut geodesic existence}. By Lemma~\ref{hierarch paths for taut geodesics}, the combinatorial geodesics $\alpha_i$ are uniform hierarchy paths.
    
    Our goal is to extract a limiting ray from $c_i$'s. 
    By the $\delta$-hyperbolicity of $\widehat{X}$ and the definition of convergence at infinity, we have the following: 
    \begin{itemize}
        \item[(*)] For each $R > 0$, there exists $n$ so that for $i,j\geq n$ the initial segment of length $R$ in $c_i$ is within $\delta$ of $c_j$, and vice versa.
    \end{itemize}

We plan to improve this coarse stabilization to exact stabilization in the locally infinite graph $\widehat{X}$.

    Let $C$ be the constant from Proposition~\ref{Bounded Geodesic Image} (BGI). Pick $R > 20(\delta+C + 1)$ and let $n$ be as in (*). For $j\geq n$, let $v_j'$ be a vertex of $c_j$ within distance $R/2$ of $o$ in $\widehat{X}$. 

    \begin{claim*}
        The set $\bigcup_{j\geq n}\left( c_j^{(0)} \cap B_{R/2}(o) \right)$ is finite.
    \end{claim*}
    \begin{proof}
    Let $v$ be a vertex of $c_n$ such that $d_{\widehat{X}}(v,v_j') \leq \delta$. Such $v$ exists by (*). Let $W\in \fS$ be a domain (in the HHS structure of $W_\Gamma$) relevant to the pair $(v,v_j')$, that is, $d_{W}(v,v_j') > M$, where $M = \dim(X)$.  

    Recall that there is a fixed point $\rho^W_{\widehat{X}} \in \widehat{X}$. Since $d_W(v,v_j') > M$, every geodesic from $v$ to $v_j'$ in $\widehat{X}$ passes through a $C$-neighborhood of $\rho^W_{\widehat{X}}$ by BGI. In particular, $d_{\widehat{X}}(\rho^W_{\widehat{X}},v_j') \leq \delta + C$ and $d_{\widehat{X}}(o, \rho_{\widehat{X}}^{W}) \leq \frac{R}{2} + (\delta + C)$.  
    
    We claim that for every $j\geq n$, $d_W(x_n,x_j)$ is uniformly bounded from above. Let $w$ and $w_j'$ be points on $c_n$ and $c_j$ that are at least $5C+5\delta + 1$ further from $o$ than $v$ and $v_j'$, respectively, and such that $d_{\widehat{X}}(w,w_j')\leq \delta$. This is possible by (*) since $R$ is large enough. Then we can connect $w$ and $w_j'$ by a geodesic in $\widehat{X}$ all of whose vertices are at distance at least $C$ from $\rho^W_{\widehat{X}}$. Then by applying BGI, $d_{W}(w,w_j') \leq  M$. Since the geodesic segment from $w$ to $x_n$ and the geodesic segment from $w_j'$ to $x_j$ are at least $C$-apart from $\rho_{\widehat{X}}^{W}$, by BGI again, $d_{W}(w,x_n)$ and $d_W(w_j',x_j)$ are uniformly bounded from above. Hence, $d_W(x_n,x_j)\leq 3M$.

    Recall that the image of a hierarchy path in each hyperbolic space is a uniform unparametrized quasi-geodesic. In particular, since $v\in \alpha_n, v_j' \in \alpha_j$, where $\alpha_n,\alpha_j$ are hierarchy paths, there is a constant $E = E(X,\fS)$ such that 
    \[
    d_W(v,x_n) \leq d_W(o,x_n) + E, \qquad d_W(x_j, v_j') \leq d_W(x_j, o) + E.
    \]

    Then for all $j\geq n$, we have
    \begin{align*}
        d_W(v,v_j') 
        &\leq d_W(v,x_n) + d_W(x_n,x_j) + d_W(x_j,v_j') \\
        &\leq d_W(o,x_n) + E + d_W(x_n,x_j) + d_W(o, x_j) +E  \\
        &\leq d_W(o,x_n) + d_W(x_n,x_j) + d_W(o, x_n) + d_W(x_n, x_j) + 2E\\
        &\leq 2d_W(o, x_n) + 6M +2E.
    \end{align*}
    
    For domains not satisfying the relevance inequality, the projection distance is at most $M$. Let $s_0$ be the constant from the Distance Formula (Theorem~\ref{Distance Formula}). Take a threshold $s \geq 2s_0+(6M +2E) + \delta$. By the preceding computation and Remark~\ref{finitely many relative domains}, every domain contributing a nonzero term to the thresholded sum in the distance formula is contained in a fixed finite set of domains
    \[
    \fS_{x_n} = \{W\in \fS\mid   d_W(o, x_n) \geq s_0\}.
    \]
    Thus,
    \[
    \sum_{W\in \fS} \{\!\!\{ d_W(v,v_j')\}\!\!\}_s 
    \leq \sum_{W\in \fS_{x_n}} d_W(v,v_j')
    \leq \sum_{W\in \fS_{x_n}} 2d_W(o,x_n) + (6M +2E).
    \]
    Since the right-hand side does not depend on $j,v_j'$ or $v$, the distance formula provides a uniform upper bound for $d(v,v_j')$. Since $v$ is a vertex of a fixed path $c_n$ and $X$ is proper (locally finite), we have that the set 
    \[
    \bigcup_{j\geq n}\left( c_j^{(0)} \cap B_{R/2}(o) \right)
    \]
    is finite.
    \end{proof}

    Thus, there are only finitely many possibilities for the initial segment of length $R/2$ of $c_j$ for all $j\geq n$. Then we can take a subsequence $c_{j(R)}$ for which these initial segments are all the same. Then increase $R$ and make the usual diagonalization argument. There is a convergent subsequence of $\{c_i\}$; call the limit $\mu$. It accumulates on $\nu$ since its initial segments are equal to the initial segments of $c_j$ for all large $j$ in the subsequence. Thus, $\mu$ is the desired geodesic ray.
\end{proof}

We can now join any two distinct boundary points in $\partial\widehat{X}$ by a bi-infinite geodesic.

\begin{prop}[Visibility of $\partial \widehat{X}$]
\label{Visibility of boundary}
    For every pair of distinct points $a,b\in \partial \widehat{X}$, there is a geodesic $\gamma:(-\infty,\infty) \to \widehat{X}$ such that $\gamma(-\infty) = a$ and $\gamma(\infty) = b$.
\end{prop}

\begin{proof}
Fix $o\in \widehat{X}$.
Let $\nu,\mu$ be geodesic rays based at $o$ converging to $a,b$, respectively. Let $\nu_i=\nu(i),\mu_i=\mu(i)$ be points on $\nu,\mu$, respectively. Let $\gamma_i$ be a taut geodesic from $\mu_i$ to $\nu_i$. Consider the geodesic triangle with sides $(o,\mu_i),(o,\nu_i)$ and $\gamma_i$. There is $T\in \N$ such that for each $i\geq T$, $\mu_i \notin N_{\delta}(\nu)$. Since each geodesic triangle is $\delta$-slim, $\gamma_i \cap B_{\delta}(\mu_T) \neq \emptyset$. Suppose $p_i\in B_{\delta}(\mu_T) \cap \gamma_i$ for each $i\geq T$. In particular, $d_{\widehat{X}}(p_i,p_j) \leq 2\delta$ for any $i,j \geq T$. We claim that there is a constant $\delta'$ satisfying the following.
    \begin{itemize}
        \item[(**)] For any $R > 0$ there exists $N$ such that for all $i\geq N$ the subsegments of $\gamma_i$ of radius $R$ centered on $p_i$ are within $\delta'$ of each other.
    \end{itemize}
    Indeed, we fix $p_i$ with $i\geq T$. By $\delta$-hyperbolicity and the definition of the geodesic ray $\mu$, there is $N_1\in \N$ such that for any $j,k \geq N_1$ the initial segment of length $R$ in $c_j$ is within $\delta$ of $c_k$, and vice versa, where $c_j,c_k$ are taut geodesics from $p_i$ to $\mu_j,\mu_k$, respectively. Similarly, there is $N_2\in \N$ such that for any $j,k\geq N_2$ the initial segment of length $R+2\delta$ in a taut geodesic from $p_i$ to $\nu_j$ is within $\delta$ of a taut geodesic from $p_i$ to $\nu_k$, and vice versa. Take $N = \max\{N_1,N_2,T\}$. 

     Let $\gamma_i^{+}$ be a subpath of $\gamma_i$ from $p_i$ to $\mu_i$ of length $R$. We want to show that $\gamma_i^+$ and $\gamma_j^+$ are within distance $\delta' = 3\delta$ for any $i,j \geq N$. Since $i\geq N \geq N_1$, the initial segment with length $R$ of $c_j$ is within $\delta$ of $\gamma_i^+$, and vice versa. Furthermore, since the distance between $p_i$ and $p_j$ is bounded from above by $2\delta$, the initial segment with length $R$ of $c_j$ is within $2\delta$ of $\gamma_j^+$, and vice versa.
     
     The same argument works for subpaths of $\gamma_i,\gamma_j$ from $p_i,p_j$, respectively, to points on $\nu$. This proves the claim.

    Then an argument similar to the proof of Lemma~\ref{Visible} then produces a subsequence converging to a bi-infinite geodesic with endpoints $a$ and $b$.
\end{proof}

\subsection[Geodesic almost extendability]{Geodesic almost extendability of $\widehat{X}$}
\label{Ideal Triangles subsection}

In a geodesic hyperbolic space, ideal triangles are uniformly slim; see \cite[Lemma 1.84]{kapovich2024trees}. A geodesic triangle is \emph{ideal} if at least one of its vertices lies in the Gromov boundary.

We use the uniform slimness for $\widehat{X}$ and the $W_\Gamma$-action to show that $\widehat{X}$ is \emph{geodesic almost extendable}. Mihalik observed that infinite hyperbolic groups have this property \cite[Lemma~3.1]{bestvina1991boundary}.

\begin{defn}(Geodesic almost extendability)
\label{Almost Extendable Definition}
Let $Y$ be a geodesic metric space. We say that $Y$ is \emph{geodesic almost extendable} if there is a constant $E$ such that for all points $p$ and $q$ in $Y$ there is a geodesic ray emanating from $p$ that passes within a distance $E$ of $q$. The constant $E$ is called a \textit{geodesic almost extendability constant} of $Y$.
\end{defn}

Note that in the literature, geodesic almost extendable spaces are sometimes called \emph{visual} spaces (for example, in \cite{bonk00embeddings}).

\begin{prop}[Geodesic almost extendability]
\label{exist close geodesic ray}
    Let $\widehat{X}$ be the maximized hyperbolic space associated with the right-angled Coxeter group $W_\Gamma$. Suppose that $\widehat{X}$ has nonempty Gromov boundary. Then there is a constant $C= C(\Gamma)$ such that $\widehat{X}$ is geodesic almost extendable with the constant $C$.
\end{prop}

\begin{proof}
Let $p,q\in\widehat{X}$. Lemma~\ref{empty Gromov boundary} and Proposition~\ref{Visibility of boundary} give a bi-infinite geodesic $\alpha$. Choose a vertex $a\in\alpha$. Since $W_\Gamma$ acts transitively on the vertex set, a translate $\beta=g\alpha$ passes through $q$.

Let $\beta^+$ and $\beta^-$ be the endpoints of $\beta$. By Lemma~\ref{Visible}, choose geodesic rays $\gamma^+$ and $\gamma^-$ based at $p$ and ending at $\beta^+$ and $\beta^-$, respectively. These rays together with $\beta$ form an ideal triangle. Since the ideal triangles are uniformly slim for some constant $\delta'$ only depending on $\delta = \delta(\Gamma)$,  $q\in\beta$ lies within $\delta'$ of $\gamma^+\cup\gamma^-$. One of these two rays therefore has nonempty intersection with $B_{\delta'}(q)$, so $\delta'$ is an almost-extendability constant.
\end{proof}

\subsection[Noncompactness of the boundary]{Noncompactness of $\partial \widehat{X}$}
\label{Non-compactness subsection}

In this subsection, we will show that when $\partial\widehat{X}$ is nonempty and $W_\Gamma$ is not hyperbolic, the boundary $\partial\widehat{X}$ is non-compact, as are the boundaries of curve graphs and fine curve graphs. This contrasts with the boundary of a proper hyperbolic space, which is always compact. The proof is inspired by
\cite[Theorem~1.4]{long2024connectedness}.

The strategy is similar to showing that the Gromov boundary of an infinite-valence regular tree is non-compact, i.e. constructing a sequence of boundary points with no convergent subsequence. In the case of $\widehat{X}$, the key is to find an element $g\in W_\Gamma$ whose orbit is bounded in $\widehat X$ but whose action on the hyperbolic space $\widehat{\mathcal{C}U}$ associated with a nested domain  is loxodromic. Bounded geodesic image then shows that a suitable subsequence of the orbit $\{g^np\}$ is uniformly separated in a visual metric on $\partial\widehat X$.

\begin{thm}[Non-compactness of $\partial \widehat{X}$]
\label{noncompact}
    Let $W_{\Gamma}$ be a right-angled Coxeter group and $\widehat{X} = \widehat{X}_{\Gamma}$ be the maximized hyperbolic space. Then the Gromov boundary $\partial \widehat{X}$ is nonempty and non-compact if and only if $\Gamma$ is neither a clique nor a nontrivial join and contains an induced 4-cycle.
\end{thm}

\begin{proof}  

\emph{Necessity.}
By Lemma~\ref{empty Gromov boundary}, the boundary is nonempty exactly when $\Gamma$ is neither a clique nor a nontrivial join. If $\Gamma$ has no induced square, then $W_\Gamma$ is hyperbolic and no standard product regions are coned off. Thus $\widehat{X}=X^{(1)}$ is proper and its Gromov boundary is compact.

\emph{Sufficiency.}
Conversely, suppose that $\Gamma$ is neither a clique nor a nontrivial join and that it contains an induced $4$-cycle labeled by $\{a,b,c,d\}$ in order.  Let $\Lambda=\{a,c\}$. Then $\Lambda$ is not a clique and $\lk(\Lambda)$ contains the non-adjacent vertices $b,d$. Hence $\Lambda$ determines a non-maximal domain $U$ of the maximized HHS structure, and its associated coned-off contact graph $\widehat{\mathcal{C}U}$ has infinite diameter. The element $g=ac\in W_{\Gamma}$ acts loxodromically on $\widehat{\mathcal{C}U}$. On the other hand, $g$ preserves $X_\Lambda$, which lies in the standard product region $X_{\Lambda*\lk(\Lambda)}$ and therefore has diameter at most one in $\widehat{X}$.

Fix a basepoint $o\in\widehat{X}$, a boundary point $p\in\partial\widehat{X}$ and a visual metric $d_{\infty}$ on $\partial\widehat{X}$. Choose a representing sequence $(x_i)$ for $p$. For $n\in\mathbb{Z}$ set $p_n=g^np$. We want to show that the sequence $(p_n)\subset \partial\widehat{X}$ has no convergent subsequence. 

Let $\tau>0$ be the stable translation length of $g$ on $\widehat{\CC U}$. Choose $L$ so that $L\tau>M$, where $M=\dim X$. For any $s,t\in L\mathbb{Z}$ with $s\neq t$ and every $i$,
\[
d_U(g^sx_i,g^tx_i)\geq |s-t|\tau>M.
\]
Proposition~\ref{Bounded Geodesic Image} therefore implies that every geodesic in $\widehat{X}$ from $g^sx_i$ to $g^tx_i$ meets a fixed $C$-neighborhood of $\rho^U_{\widehat{X}}$. In a hyperbolic space, the Gromov product based at $o$ coarsely equals the distance from $o$ to a geodesic between the two points. Consequently there is a constant
\[
D=D\bigl(\delta,C,d_{\widehat{X}}(o,\rho^U_{\widehat{X}})\bigr)
\]
such that
\[
\langle g^sx_i,g^tx_i\rangle_o\leq D
\]
for all $i$ and all distinct $s,t\in L\mathbb{Z}$. Since $(g^sx_i)$ and $(g^tx_i)$ represent $p_s$ and $p_t$, respectively, the definition of the boundary Gromov product gives $\langle p_s,p_t\rangle_o\leq D$. Proposition~\ref{visual metric property} now provides a constant $\epsilon>0$, independent of $s$ and $t$, such that
\[
d_\infty(p_s,p_t)\geq\epsilon
\qquad(s,t\in L\mathbb{Z},\ s\neq t).
\]
Thus $(p_{nL})_{n\geq0}$ has no convergent subsequence. Since the boundary is metrizable, it is non-compact.
\end{proof} 

\begin{rmk}
    By arguing as in the proof of Theorem~\ref{noncompact}, one can show the Gromov boundary $\partial\widehat{X}_{\Gamma}$ is nowhere locally compact when $\Gamma$ is neither a clique nor  a nontrivial join and contains an induced 4-cycle. We leave it as an exercise for readers.
\end{rmk}

%% file: sections/fan-filter.tex
The fan-and-filter construction is the main combinatorial ingredient in the proof of Theorem~\ref{main theorem}. The aim is, by using fans and filters, to show that two points in a sphere with radius $R$ in $\widehat{X}$ can be connected by a path such that (1) vertices of this path are at distance at least $R$ from the center of the sphere; and (2) the diameter of this path has a linear upper bound with respect to the distance between endpoints with the linearity constant only depending on $\Gamma$ (Proposition~\ref{spheres are connected}).

A fan is a van Kampen diagram, which locally pushes a path in $X$ away from one standard product region. Iterating this operation produces a filter, which is also a van Kampen diagram of $W_\Gamma$ whose directed paths are combinatorial geodesics in $X_\Gamma$ and do not have arbitrarily long subpaths in standard product regions. By using filters, one can push a path in $X$ away from certain unions of standard product regions whose image in $\widehat{X}$ forms a closed ball.

We begin with a warm-up to show when a standard product region does not separate $X$.

Throughout this section, we denote by $X=X_\Gamma$ the Davis complex associated with the graph $\Gamma$ and  we assume that the defining graph $\Gamma$ satisfies the following conditions:

\begin{enumerate}[label=(\alph*), ref=(\alph*)]
    \item \label{Condition:(a)}
        $\Gamma$ is connected, is not a clique, and has no separating clique;
    \item \label{Condition:(b)}
        $\Gamma$ has no separating nontrivial join;
    \item \label{Condition:(c)}
        for every vertex $v\in V(\Gamma)$, the link of $v$ is not contained in a nontrivial join.
\end{enumerate}

\subsection{Avoiding one standard product region}
In Proposition~\ref{Single product region does not separate}, we will show that an edge path in $\Gamma$ avoiding a nontrivial join could produce a combinatorial path in $X=X_\Gamma$ avoiding a standard product region. This idea will motivate the construction of fans in Subsection~\ref{Fans}.

Recall that a combinatorial path $\alpha$ in $X$ determines a word $\overline{\alpha}$ by reading edge labels, and $\lab(\alpha) \subset V(\Gamma)$ is the set of letters occurring in that word. If $\alpha$ and $\beta$ are combinatorial geodesics with the same endpoints, they cross the same set of hyperplanes and hence $\lab(\alpha)=\lab(\beta)$. We therefore write $\lab([x,y])$ for this vertex set.

\begin{prop}
\label{Single product region does not separate}
Let $X=X_\Gamma$, where $\Gamma$ satisfies \cref{Condition:(a),Condition:(b),Condition:(c)}. If $P\subset X$ is a standard product region, then $X\setminus P$ is path-connected. 
\end{prop}
\begin{figure}[t]
    \centering
    \includegraphics[scale=0.8]{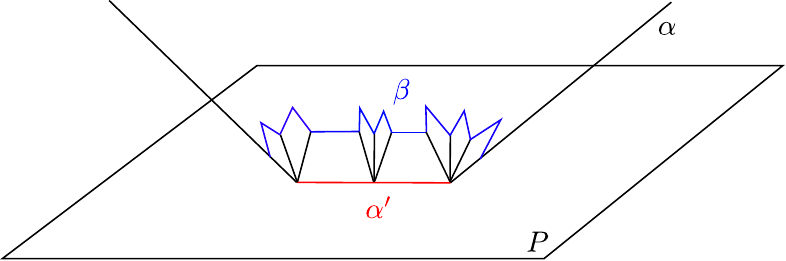}
    \caption{A subpath $\alpha'$ of $\alpha$ contained in a product region $P$
    is pushed to a subpath of $\beta$ across two hyperplane carriers and three
    fans.}
    \label{PushingPathFig}
    \end{figure}

\begin{proof}
Let $x,y\notin P$ and choose a combinatorial geodesic $\alpha$ from $x$ to $y$. There is nothing to prove if $\alpha\cap P=\emptyset$. Otherwise let $P=gX_\Lambda$, where $\Lambda\subset \Gamma$ is a maximal nontrivial join and $g\in W_{\Gamma}$. Convexity of $P$ implies that $\alpha\cap P$ is a (possibly trivial) subpath $\alpha' \subset \alpha$. Write
\[
\overline{\alpha'}=a_1\cdots a_n,
\]
where $a_i\in V(\Gamma)$ for $i=1,\cdots,n$. Denote the successive vertices of $\alpha'$ by $p_0,\cdots,p_n$. Let $e_v$ and $e_w$ be the edges of $\alpha\setminus\alpha'$ incident to $p_0$ and $p_n$, with labels $v$ and $w$ in $V(\Gamma)$, respectively. Note that every 1-cube in $X$ containing a 0-cube in $gX_\Lambda^{(0)}$ lies in $gX_\Lambda$ exactly when its label lies in $\Lambda$. Since $e_v$ and $e_w$ are not in $P$, while $p_0,p_n\in P$, we have $v,w \notin \Lambda$.   Since $\alpha'\subset P$, $a_i\in \Lambda$ for each $i = 1,\cdots,n$. See Figure~\ref{PushingPathFig}.

For each $1\leq i \leq n$, \cref{Condition:(c)} gives a vertex $a_i'\in\lk(a_i)\setminus \Lambda$, because, otherwise, $\lk(a_i)$ would be contained in the nontrivial join $\Lambda$. \cref{Condition:(c)} also guarantees the defining graph $\Gamma$ is not a nontrivial join. In particular, $\Gamma \neq \Lambda$. Then \cref{Condition:(b)} implies that $\Gamma\setminus \Lambda$ is connected. Set $a_0'=v$ and $a_{n+1}'=w$. Hence, for $0\leq i\leq n$, choose an edge path
\[
\tau_i=(a_i'=t_0^i,t_1^i,\ldots,t_{r_i}^i=a_{i+1}')
\]
in $\Gamma\setminus \Lambda$.

At $p_i$, the edge path $\tau_i \subset \Gamma$ determines a chain of squares in $X$. Remove from the boundary of this chain of squares its first and last 1-cubes, both of which are incident to $p_i$ in $X$, and denote the remaining path by $\beta_i$. Since every 0-cube of $\beta_i$ is separated from $p_i$ by a hyperplane with label not in $\Lambda$, $\beta_i \cap P = \emptyset$. 

Since $a_i'\in\lk(a_i)$, the corresponding square labeled by $a_ia_i'a_ia_i'$ translates the $a_i$-edge of $\alpha'$ to an edge $e_i$ connecting the terminal vertex of $\beta_{i-1}$ with the initial vertex of $\beta_i$. Both endpoints of $e_i$ lie outside of $P$ (since $a_i'\notin \Lambda$), and so does $e_i$. The concatenation
\[
\beta = \beta_0e_1\beta_1e_2\cdots e_n\beta_n
\]
has the same endpoints as $e_v\alpha'e_w$ and is disjoint from $P$.  Replacing $e_v\alpha'e_w$ by $\beta$ gives a path from $x$ to $y$ in $X\setminus P$.
\end{proof}

\subsection{Fans}
\label{Fans}
We now define a \emph{fan}, a particular connected van Kampen diagram of $W_\Gamma$. A van Kampen diagram for $W_\Gamma$ is a planar square complex $D$ equipped with a cellular map $D\to X$ whose restriction 
\[
D^{(1)} \to X^{(1)}
\]
preserves edge labels. Hence every square records a commuting relation in $W_\Gamma$. We note that there may be distinct vertices and edges of $D$ with the same image in $X$.

Let $\gamma$ be a path with label $\overline{\gamma}=s_0\cdots s_n$. Its \emph{wide tail} $\gamma'$ is its longest terminal subpath whose labels lie in a nontrivial join. Equivalently, $\gamma'$ has the image contained in a standard product region of $X_\Gamma$. The wide tail may be empty. Let $f_0$ and $f_r$ be edges incident to the terminal vertex of $\gamma$, with labels $v_0,v_r$, such that
$\overline{\gamma f_0}$ and $\overline{\gamma f_r}$ are geodesic.

The following definition is from \cite[Definition 6.1]{cordes2025connectivity} for right-angled Coxeter groups.

 \begin{defn}[Fan]
 \label{Fan Definition}
A \emph{fan} $F$ is a van Kampen diagram consisting of a path $\gamma$, a
sequence of at least three edges $f_0,\ldots,f_r$ incident to the terminal
vertex of $\gamma$, and, for each $0\leq i<r$, a square $c_i$ spanned by
$f_i$ and $f_{i+1}$ with the following properties:
\begin{enumerate}
    \item $\overline{\gamma f_i}$ is geodesic in $W_\Gamma$ for each $0 \leq i \leq r$.
    \item If $\overline{\gamma'}$ is the wide tail of $\overline{\gamma}$, then
    \begin{enumerate}
        \item either $\lab(\gamma')$ is a clique;
        \item or there is a nontrivial join subgraph $\Delta\subset \Gamma$ such that $\lab(\gamma') \subset \Delta$ and $v_i \notin \Delta$, where $v_i = \lab(f_i)$, for each $0 < i < r$.
    \end{enumerate}   
\end{enumerate}
\end{defn}

We call $\gamma$ the \emph{base path}, $f_i$ the \emph{fan edges}. In particular, $f_0,f_r$ are called the \emph{left} and \emph{right} fan edges, and the other $f_i$'s are \emph{interior fan edges}. Orient $\gamma$ toward the common initial vertex of the fan edges and order the fan edges from left to right. A path is \emph{directed} if it follows the resulting orientation. In each square $c_i$, the vertex opposite the terminal vertex of $\gamma$ is called the \emph{top vertex}. Moreover, the edge opposite $f_{i+1}$ is the \emph{top-left} edge and the edge opposite $f_i$ is the \emph{top-right} edge.

The next lemma gives the required fans. For proofs, see \cite[Lemmas~5.5 and~5.6]{mihalik2007local} in the right-angled case and \cite[Lemma~6.4]{cordes2025connectivity} for wide-spherical-avoidant Coxeter groups. 

By Lemma~\ref{conditions are equivalent}, the graph hypotheses (conditions \textup{(a)}--\textup{(c)}) in our main theorem imply the assumptions needed for the fan-and-filter construction in the literature.

\begin{lem}[Existence of fans]
\label{fan}
    Suppose $\Gamma$ satisfies conditions \textup{(a)}--\textup{(c)}.  Let $\gamma$ be an abstract labeled path with $\overline{\gamma}\in W_\Gamma$ being geodesic, and let $e_l,e_r$ be distinct abstract edges incident to the terminal vertex of $\gamma$. Suppose that $\overline{\gamma e_l}$ and $\overline{\gamma e_r}$ are geodesic. Then there is a fan
    \[
    F(\gamma,e_l,e_r)
    \]
    whose base path is $\gamma$ with left and right fan edges $e_l$ and $e_r$, respectively.
\end{lem}

The following lemma is from \cite[Lemma 6.5]{cordes2025connectivity} and also outlined in \cite[Lemma 3.9]{mihalik2007local}. We give a proof in the case of right-angled Coxeter groups for completeness.

\begin{lem}
\label{directed paths are geodesic}
    Let $F$ be a fan with base path $\gamma$. Let $\rho$ be a directed path in $F$ starting at the endpoint of $\gamma$. Then $\overline{\gamma \rho}$ is geodesic.
\end{lem}
\begin{proof}
A directed path beginning at the terminal vertex of $\gamma$ consists of one fan edge, possibly followed by one top fan edge contained in a square. The one-edge case is guaranteed to be geodesic by the definition of the fan. In the two-edge case, write $\overline{\gamma}=s_0\cdots s_n$ and suppose the two edge labels are $s,t$. By the definition of a fan, both $s_0\cdots s_ns$ and $s_0\cdots s_nt$ are geodesic. If $s_0\cdots s_nst$ were not geodesic, the deletion condition would produce an earlier occurrence of $t$ commuting with every letter in between. The same cancellation would shorten $s_0\cdots s_nt$, a contradiction. These two cases exhaust the directed paths in a single fan.
\end{proof}

\subsection{Filters}
\label{Filter Construction Subsection}

Given two geodesic rays based at the same point in $X$, we will build a one-ended van Kampen diagram, called a \emph{filter}, by attaching fans inductively. This filter construction  was used to establish local connectivity of the visual boundary of RACGs in \cite{mihalik2007local} and  connectivity and local connectivity of the Morse boundary of more general Coxeter groups in \cite{cordes2025connectivity}. We study a truncated filter whose image in $\widehat{X}$ produces a path in spheres with the desired radii.

Let $\alpha$ and $\beta$ be combinatorial geodesic rays with a common initial vertex $o\in X$. Consider $\alpha$ and $\beta$ as abstract labeled paths mapping to $X^{(1)}$. Let $\eta$ be their maximal common initial subpath and let $z_0$ be its terminal vertex. Assume that the rays have distinct first edges after $z_0$, and write
\[
\alpha=\eta a_0a_1\cdots,
\qquad
\beta=\eta b_0b_1\cdots.
\]
The path $\eta$ is called the \emph{common stem}.

Begin with the abstract union $\FF_0=\alpha\cup_\eta\beta$ and attach the fan
\[
F(\eta,a_0,b_0).
\]
This is the level-one fan, and the resulting diagram is $\FF_1$. We orient $\alpha$ and $\beta$ from $o$ and $\FF_1$ has the induced orientation.  After deleting the top-left edge of every square in the level-one fan, the remaining edges, together with $\FF_0^{(1)}$, form a spanning tree $T_1$.

Suppose that $\FF_N$ and its spanning tree $T_N$ have been constructed.  If $v$ is the terminal vertex of a fan edge in a level-$N$ fan, there are two outgoing edges immediately above $v$, ordered from left to right based on the orientation of $\FF_{N}$.  Let $\gamma_v$ be the unique path in $T_N$ from $o$ to $v$, and denote the two outgoing edges by $e_v^l,e_v^r$. The simultaneous induction in Lemma~\ref{basic properties of filters} says that both $\gamma_ve_v^l$ and $\gamma_ve_v^r$ are geodesic. By applying Lemma~\ref{fan}, we attach a fan
\[
F(\gamma_v,e_v^l,e_v^r)
\]
at every such vertex. These are the level-$(N+1)$ fans.  The resulting diagram is $\FF_{N+1}$, and $T_{N+1}$ is obtained by deleting the top-left edge of every square in every fan.  Define
\[
\FF=\varinjlim\FF_N, \qquad T=\varinjlim T_N
\]
by taking the direct limit. We call $\FF$ the \emph{full filter} spanning $\alpha$ and $\beta$, and $T$ its \emph{spanning tree}.

We record some properties of full filters. 

\begin{lem}[Basic properties of filters]
\label{basic properties of filters}
For every $N\geq1$,
\begin{enumerate}
    \item $T_N$ is a spanning tree of $\FF_N^{(1)}$, and
    $T_N\subset T_{N+1}$;
    \item every directed path in $\FF_N$ maps to a combinatorial geodesic
    in $X$;
    \item every vertex beyond the common stem has one or two incoming edges in $\FF_N$, depending on whether it is the top vertex of a square in a fan. And exactly one of these incoming edges belongs to $T_N$. Consequently, the unique path in $T_N$ from $o$ to any vertex is directed;
    \item the path $\eta$ is contained in the tree path from $o$ to every vertex in the part of the filter beyond $z_0$;
    \item every terminal vertex of a level-$N$ fan edge has exactly two
    outgoing edges, and these are the left and right fan edges of the unique fan attached there at level $N+1$.
\end{enumerate}
\end{lem}

We refer readers to the discussion before~\cite[Lemma 5.8]{mihalik2007local} or \cite[Lemmas~6.6 and~6.7]{cordes2025connectivity} for a proof.

\begin{figure}[t]
    \centering
    \includegraphics[scale=1]{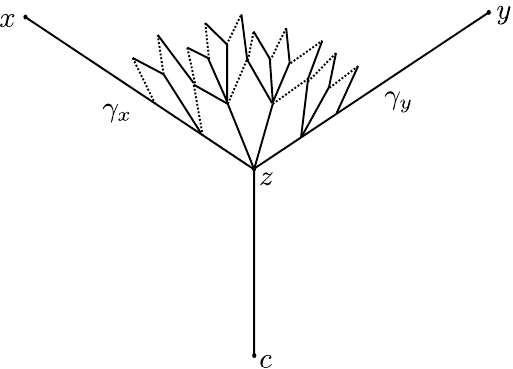}
    \caption{The first two levels of a filter.  The edges of the spanning tree $T_2$ are
    solid, and the deleted top-left edges are dashed.}
    \label{Filter}
\end{figure}

We may also use finite paths instead of rays for a van Kampen diagram. Recall that every finite combinatorial geodesic in $X$ extends to a combinatorial geodesic ray. At each endpoint the set of labels that shorten the current word spans a clique (Remark~\ref{B(g) is a clique remark}), and this clique is a proper subgraph of $\Gamma$ because $\Gamma$ is not a clique. Thus finite side paths of a van Kampen diagram may be extended to rays, and the corresponding full filter may be fixed once and for all.

Recall that we identify the vertex set of the Davis complex $X^{(0)}$ and the vertex set of the maximized hyperbolic space $\widehat{X}^{(0)}$ by the natural inclusion map $X \to \widehat{X}$. We denote the sphere in $(\widehat{X},d_{\widehat{X}})$ with radius $R$ based at a vertex $o$ by 
\[
S_R(o) = \{x\in \widehat{X} \mid d_{\widehat{X}}(o,x)  = R\}.
\]

Fix a full filter $\FF$ spanning $\alpha$ and $\beta$ based at $o$. Let $\gamma_x$ and $\gamma_y$ be finite initial subpaths of $\alpha$ and $\beta$, respectively, ending at $x,y\in S_R(o)$.  We define the \emph{radius-$R$ truncation} $\FF^R$ as a subdiagram of $\FF$ as follows.

Initially, $z_0$ is \emph{eligible for extension} and its prescribed fan from $\FF$ is $F(\eta,a_0,b_0)$.  Thereafter, every endpoint of a fan edge is \emph{eligible for extension}. Lemma~\ref{basic properties of filters} provides its unique next-level fan from the construction of $\FF$. 
We say that an eligible vertex $z$ is \emph{pending} if 
\[
z\notin \{x,y\} \qquad \text{and}\qquad d_{\widehat{X}}(o,z) < R.
\]

We attach the prescribed fan at every pending eligible vertex at every stage, and do not attach a fan at an eligible vertex whose distance from $o$ is at least $R$ with respect to the $d_{\widehat{X}}$-metric.  We also stop at $x$ and $y$.  Denote the finite-stage diagrams and trees by $\FF_N^R$ and $T_N^R$. Proposition~\ref{spheres are connected} will show that the construction terminates after finitely many stages.

At each finite stage $\FF_N^R$, a part of its boundary opposite the common stem is a non-backtracking edge path $\tau_N$ from $x$ to $y$, called the \emph{top boundary}.  The basic filter properties continue to hold for the truncated diagrams, since they are subdiagrams of the fixed full filter.

\begin{lem}[Distance control in a truncated filter]
\label{distance control in a truncated filter}
If $e$ is a directed edge of $T_N^R$ from vertex $u$ to vertex $v$, then
\begin{equation}
\label{distance control in a truncated filter inequality}
d_{\widehat X}(o,u) \leq d_{\widehat X}(o,v) \leq d_{\widehat X}(o,u)+1.
\end{equation}

If the truncation process has no pending eligible boundary vertices at stage
$N$, then
\[
\tau_N^{(0)}\subset S_R(o)\cup S_{R+1}(o).
\]
\end{lem}

\begin{proof}
The tree path from $o$ to $v$ is a combinatorial geodesic containing $u$, so $u\in[o,v]$.  Lemma~\ref{distance monotonicity in hull} gives the first inequality, and the second follows because $u$ and $v$ are adjacent in $X$, hence also adjacent in $\widehat{X}$.

Suppose that no pending eligible boundary vertex remains. Every eligible vertex $p$ of $\tau_N$ is the endpoint of a fan edge based at a pending eligible vertex $p'$. Since $d_{\widehat{X}}(o,p') \leq R-1$, we have $d_{\widehat{X}}(o,p)\leq R$. Furthermore, since $p$ is not pending,  $d_{\widehat{X}}(o,p) = R$, so $p\in S_R(o)$.

Each non-eligible vertex $q\notin \{x,y\}$ of $\tau_N$ is adjacent to an eligible vertex. Hence by the above inequality~\ref{distance control in a truncated filter inequality}, $R\leq d_{\widehat{X}}(o,q) \leq R+1$. Therefore, every vertex of $\tau_N$ lies in $S_R(o)\cup S_{R+1}(o)$.
\end{proof}

Another property needed below is that a directed path away from the side rays $\alpha$ and $\beta$ cannot remain for too long in one standard product region. 
\begin{lem}[Bounded wide paths]
\label{Bounded wide path length}
There is a constant $C_0 = C_0(\Gamma)$ with the following property. Let $\zeta$ be a directed path in the spanning tree $T$.  Suppose that, as a path in the abstract diagram, $\zeta$ meets $\alpha\cup\beta$ in at most its initial vertex.  If $\lab(\zeta)$ is contained in a nontrivial join, then
\[
|\zeta|\leq C_0.
\]
\end{lem}

We refer readers to~\cite[Lemmas~5.8--5.9]{mihalik2007local} for a proof. Cordes and Levcovitz prove the corresponding qualitative boundedness statement for the filters for general Coxeter groups in
\cite[Proposition~6.19]{cordes2025connectivity}.

\subsection{Coarse connectedness of spheres in $\widehat X$}
\label{spheres are connected subsection}

We now apply the filter to show the spheres in $\widehat{X}$ are coarsely connected.  There are two estimates in the proof.  First, the bounded-wide-path lemma forces every branch to reach radius $R$ after finitely many stages in the construction of truncated filters. Second, a long common stem keeps the entire top boundary in a uniformly bounded subset of $\widehat X$.

Recall that the truncated filter $\FF^R_n$ is constructed as a van Kampen diagram whose image is contained in $X$ and we measure the distance under the metric $d_{\widehat{X}}$.

\begin{prop}[Coarse connectedness of spheres]
\label{spheres are connected}
Fix $o\in X^{(0)}=\widehat X^{(0)}$.  For every $E \geq 0$ there is a constant $D=D(\Gamma,E)$ with the following property.  If $R\in\N$ and $x,y\in S_R(o)$ satisfy
\[
d_{\widehat X}(x,y)\leq E,
\]
then there is a combinatorial path
\[
\tau=(x=x_0,x_1,\ldots,x_m=y)
\]
in $X^{(1)}$ such that
\[
\tau^{(0)}\subset S_R(o)\cup S_{R+1}(o)
\]
and
\[
\diam_{\widehat X}\{x_0,\ldots,x_m\}\leq D.
\]
One can take
$
D = 4E + 4M +6$,
where
$
M=\max\{|K|:K\subset\Gamma\text{ is a clique}\}.
$
\end{prop}

\begin{proof}
If $x=y$, take the constant path.  Now we can assume $E\in\N$ and $E\geq1$. If $E<R$, then $d_{\widehat{X}}(x,y)<R$, and Lemma~\ref{large intersection of convex hull} gives a vertex
\[
z\in[o,x]\cap[o,y] \quad\text{with}\quad d_{\widehat X}(o,z)\geq R- d_{\widehat{X}}(x,y)\geq R-E.
\]
If $E\geq R$, take $z=o$.  In either case,
\begin{equation}
\label{common stem lower bound}
d_{\widehat X}(o,z)\geq R-E.
\end{equation}

Choose a combinatorial geodesic from $o$ to $z$ and extend it to combinatorial geodesics $\gamma_x$ and $\gamma_y$ from $o$ to $x$ and $y$, respectively. Let $\eta$ be their maximal common initial subpath, let $z_0$ be its terminal vertex, and let
\[
s=d_{\widehat X}(o,z_0).
\]
The common stem contains the chosen path from $o$ to $z$.  Hence Lemma~\ref{distance monotonicity in hull} and the inequality~\eqref{common stem lower bound} give
\begin{equation}
\label{maximal stem lower bound}
R- E\leq  d_{\widehat X}(o,z)\leq s.
\end{equation}
Since $z_0\in[o,x]\cap[o,y]$, by Lemma~\ref{distance monotonicity in hull} again, we have $s\leq R$.

If $s=R$, no filter is needed.  Every vertex $u$ on either terminal subpath from $z_0$ to $x$ or $y$ satisfies
\[
R=d_{\widehat X}(o,z_0) \leq d_{\widehat X}(o,u) \leq R,
\]
so the concatenation of these two subpaths lies in $S_R(o)$.  Moreover, $z_0\in[o,u]$.  Lemma~\ref{convex hulls are uniformly thin} provides a point $z_u$ on a taut geodesic from $o$ to $u$ with $d_{\widehat X}(z_0,z_u)\leq M + 1$.  Since both $z_0$ and $u$ have distance $R$ from $o$, it follows that $d_{\widehat X}(u,z_0)\leq 2(M+1)$.  The resulting path has diameter at most $4(M+1)$, as required.

Assume from now on that $s<R$.  In particular, both $\gamma_x$ and $\gamma_y$ continue after $z_0$. Extend them to combinatorial geodesic rays $\alpha$ and $\beta$, build the full filter with common stem $\eta$, and construct its radius-$R$ truncation, stopping at $x$ and $y$.

We will show that this truncation terminates after finitely many stages.  Let $C_0 = C_0(\Gamma)$ be the constant from Lemma~\ref{Bounded wide path length} and let $C_1=C_1(C_0,\Gamma)$ be the constant given by Lemma~\ref{Morse Geodesic C'}.  Set
\[
L=(E+2M + 2)(C_1+1).
\]

Let $p$ be a pending eligible vertex in a finite stage of the truncated tree.  If $p$ does not lie on $\gamma_x$ or $\gamma_y$, let $q$ be the last vertex at which the unique tree path from $o$ to $p$ meets the union of paths $\gamma_x\cup\gamma_y$, and let $\zeta$ be the subpath from $q$ to $p$. Every wide subpath of $\zeta$ therefore has length at most $C_0$, by Lemma~\ref{Bounded wide path length}.  The path $\zeta$ maps to a combinatorial geodesic by Lemma~\ref{basic properties of filters}, so $d_X(q,p)=|\zeta|$.   Hence, by Lemma~\ref{Morse Geodesic C'}, every wide subpath of every combinatorial geodesic from $q$ to $p$ has length at most $C_1$.

Suppose that $|\zeta|\geq L$.  Lemma~\ref{C'} shows that
\begin{equation}
\label{filter branch progress}
d_{\widehat X}(q,p)\geq E+2 (M+1).
\end{equation}
Let $t=d_{\widehat X}(o,q)$. By the choice of $q$, $q$ lies in $[z_0,x]$ or $[z_0,y]$. Thus,
\[
R- E \leq s \leq t.
\]
It follows from \eqref{filter branch progress} that
\[
d_{\widehat X}(q,p)\geq R-t+2(M+1).
\]

The directed tree path from $o$ to $p$ is a combinatorial geodesic passing through $q$, so $q\in[o,p]$.  Lemma~\ref{convex hulls are uniformly thin} gives a point $q'$ on a taut geodesic from $o$ to $p$ such that
\[
d_{\widehat X}(q,q')\leq M + 1.
\]
Consequently,
\[
d_{\widehat X}(q',p)\geq R-t+ M +1
\qquad\text{and}\qquad
d_{\widehat X}(o,q')\geq t- (M + 1).
\]
Since $q'$ lies on a geodesic from $o$ to $p$ in $\widehat{X}$, these inequalities imply
\[
d_{\widehat X}(o,p)
=d_{\widehat X}(o,q')+d_{\widehat X}(q',p)
\geq R.
\]
This contradicts the assumption that $p$ is pending.  Thus every directed subpath to a pending eligible vertex has length less than $L$ after its last intersection with $\gamma_x \cup \gamma_y$.

The two paths $\gamma_x$ and $\gamma_y$ are finite, and the abstract spanning tree is locally finite by Lemma~\ref{basic properties of filters}.  The radius-$R$ truncated filter $\FF^R$ can be obtained within an $L$-neighborhood of $\gamma_x \cup \gamma_y$. Hence only finitely many fans are attached in the process, and the truncation terminates after finitely many steps. Let
\[
\tau=(x=x_0,x_1,\ldots,x_m=y)
\]
be the terminal top boundary of $\FF^R$.  Consecutive vertices are adjacent in $X$ (hence adjacent in $\widehat{X}$ as well), and Lemma~\ref{distance control in a truncated filter} gives, for each $0\leq i \leq m$, 
\[
x_i\in S_R(o)\cup S_{R+1}(o).
\]

It remains to control the diameter.  The unique tree path from $o$ to each $x_i$ contains the common stem $\eta$, and hence $z_0\in[o,x_i]$ for each $i = 0,\cdots, m$. Lemma~\ref{convex hulls are uniformly thin} gives a point $z_i$ on a taut geodesic from $o$ to $x_i$ such that
\[
d_{\widehat X}(z_0,z_i)\leq M + 1.
\]
By \eqref{maximal stem lower bound},
\[
d_{\widehat X}(o,z_i)\geq s- (M+1)\geq R-E-(M+1).
\]
Because $d_{\widehat X}(o,x_i)\leq R+1$ and $z_i$ lies on a geodesic in $\widehat{X}$ from $o$ to $x_i$, we obtain
\[
d_{\widehat X}(x_i,z_i)\leq E+M+2.
\]
Therefore
\[
d_{\widehat X}(x_i,z_0)\leq E+ 2M+3
\]
and, for all $i,j$,
\[
d_{\widehat X}(x_i,x_j)\leq 4E+4M+6.
\]
\end{proof}

We conclude this section by a corollary about the number of ends of $\widehat{X}$.  We say a connected, unbounded graph is \emph{one-ended} if the complement of every bounded subset has exactly one unbounded connected component.

\begin{cor}[One-endedness of $\widehat{X}$]
\label{one-endedness of maximized hyperbolic space}
Suppose $\Gamma$ satisfies
\cref{Condition:(a),Condition:(b),Condition:(c)}.  Then the maximized
hyperbolic space $\widehat X_\Gamma$ is one-ended.
\end{cor}

\begin{proof}
By~\cref{Condition:(a),Condition:(c)}, we know that  $\Gamma$ is neither a clique nor a nontrivial join.  Thus $\widehat X_\Gamma$ has infinite diameter by Lemma~\ref{empty Gromov boundary}.

Fix $o\in\widehat{X}_\Gamma$ and $N\in\N$.  Let $u,v\notin B_N(o)$.  Choose geodesics in $\widehat{X}_\Gamma$ from $o$ to $u$ and $v$, and let $u',v'$ be the vertices in the intersection of chosen geodesics with $S_{N+1}(o)$, respectively. The terminal geodesic subpaths from $u'$ to $u$ and from $v'$ to $v$ avoid $B_N(o)$.  Applying Proposition~\ref{spheres are connected} with $R=N+1$ and $E=d_{\widehat{X}}(u',v')$ gives an edge path in
\[
S_{N+1}(o)\cup S_{N+2}(o)
\]
joining $u'$ to $v'$.  Concatenating the above three paths connects $u$ to $v$ in $\widehat X\setminus B_N(o)$. Hence the complement of every bounded ball based at $o$ is connected. Since $W_\Gamma$ acts on $\widehat{X}_\Gamma$ transitively, the complement of every bounded ball in $\widehat{X}_\Gamma$ is connected.

Every bounded subset is contained in such a ball. Therefore its complement has exactly one unbounded component, so $\widehat X_\Gamma$ is one-ended.
\end{proof}

%% file: sections/main-proof.tex
We now combine the results of the preceding sections.  The obstruction established in Proposition~\ref{disconnectedness of Gromov boundary} proves one direction of Theorem~\ref{main theorem}.  For the other direction, we use the controlled connectivity of spheres from Proposition~\ref{spheres are connected} to construct a path with controlled diameter joining a pair of points in the boundary.  The construction proceeds by repeatedly inserting new boundary points between two previously chosen ones.  At each stage the new points are represented by geodesic rays that fellow travel farther from the basepoint, so the diameters of the inserted boundary chains decrease geometrically. 

The proof is inspired by \cite[Proposition 2.2]{wright2024spheres}, \cite[Proposition 3.2]{bestvina1991boundary} and \cite[Proposition 5.2]{Mackay_Sisto_2020}.

\begin{proof}[Proof of Theorem~\ref{main theorem}]
Suppose first that one of conditions \textup{(a)}--\textup{(c)} fails. Proposition~\ref{disconnectedness of Gromov boundary} shows that $\partial\widehat{X}$ is either empty or disconnected. 

Conversely, assume that $\Gamma$ satisfies conditions \textup{(a)}--\textup{(c)}.  We will prove that $\partial\widehat{X}$, equipped with a visual metric $d_{\infty}$, is linearly connected. It is sufficient to show the linear connectivity since linear connectivity implies path connectivity, hence connectivity. The Gromov boundary may admit many visual metrics, but they are quasisymmetric to each other and linear connectivity is preserved under a quasisymmetry. Hence, it is enough to show the linear connectivity for the chosen visual metric $d_{\infty}$.

Condition~\textup{(a)} says that $\Gamma$ is not a clique.  Moreover, $\Gamma$ cannot itself be a nontrivial join, since otherwise $\lk(v)$ would be contained in the nontrivial join $\Gamma$ for every vertex $v$, contrary to condition~\textup{(c)}.  Lemma~\ref{empty Gromov boundary} therefore shows that $\partial\widehat{X}$ is nonempty and contains at least two points.

Let $p,q\in \partial \widehat{X}$. Fix a vertex $o\in \widehat{X}$. We will construct a continuous map $\gamma:[0,1] \to \partial \widehat{X}$ such that $\gamma(0) = p,\gamma(1) = q$ and $\diam(\gamma([0,1]))$ is bounded above by a constant multiple of $d_{\infty}(p,q)$, where the constant is independent of $p$ and $q$. See Figure~\ref{MainThmPic} for an illustration. The key idea is to construct a chain of sequences of points in the boundary.

By Lemma~\ref{Visible}, there are geodesic rays $\alpha$ and $\beta$ based at $o$ and converging to $p$ and $q$, respectively.  Fix $b>1$ and assume the visual metric $d_\infty$ is as in Proposition~\ref{visual metric property}.  The usual fellow-traveling estimate for geodesic rays in a hyperbolic space \cite{vaisala2005gromov,das2017geometry} gives constants $E_0=E_0(\delta)\geq2$ and $A_0\geq1$, depending only on $\delta$ and the chosen visual metric, with the following property.  Let $R$ be the largest integer such that
\[
 d_{\widehat{X}}\bigl(\alpha(R),\beta(R)\bigr)\leq E_0.
\]
Then $R$ differs from the Gromov product $\langle p,q\rangle_o$ by a uniformly bounded constant, and hence
\begin{equation}
\label{visual distance at branching time}
 A_0^{-1}b^{-R}\leq d_\infty(p,q)\leq A_0b^{-R}.
\end{equation}
The integer $R$ exists because the two geodesic rays represent distinct points in the boundary.

We next fix the constants that will be used in the construction.  Let $C$ be the geodesic almost extendability constant from Proposition~\ref{exist close geodesic ray}.  We may assume that $C$ is an integer.  Let
\[
 K=6C+3,
 \qquad
 E_1=K+2C+2.
\]
For $i=0,1$, let $D_i=D(\Gamma,E_i)$ be the constant supplied by Proposition~\ref{spheres are connected}.

We will use the following standard consequence of hyperbolicity and the property of a visual metric.  Suppose that $\rho$ and $\rho'$ are geodesic rays based at $o$, and that they contain points, each at distance at least $T$ from $o$, whose distance is at most $L$.  Then
\begin{equation}
\label{close deep points give close boundary points}
 d_\infty\bigl(\rho(\infty),\rho'(\infty)\bigr)
 \leq A(L)b^{-T},
\end{equation}
where $A(L)$ is a constant depending only on $L$, $\delta$, and the chosen visual metric $d_{\infty}$. Indeed, the Gromov product of the two endpoints is at least $T-O(L+\delta)$, and \eqref{close deep points give close boundary points} follows from Proposition~\ref{visual metric property}.

We now recursively construct increasing finite sets
\[
 I_0\subset I_1\subset\cdots\subset[0,1]
\]
and associated maps
\[
 \gamma_r\colon I_r\longrightarrow\partial\widehat{X}.
\]
For each $a\in I_r$, we choose a particular geodesic ray $\rho_a$ based at $o$ and satisfying $\gamma_r(a)=\rho_a(\infty)$.  Figure~\ref{MainThmPic} illustrates one step of the construction.

\begin{figure}[t]
    \centering
    \includegraphics{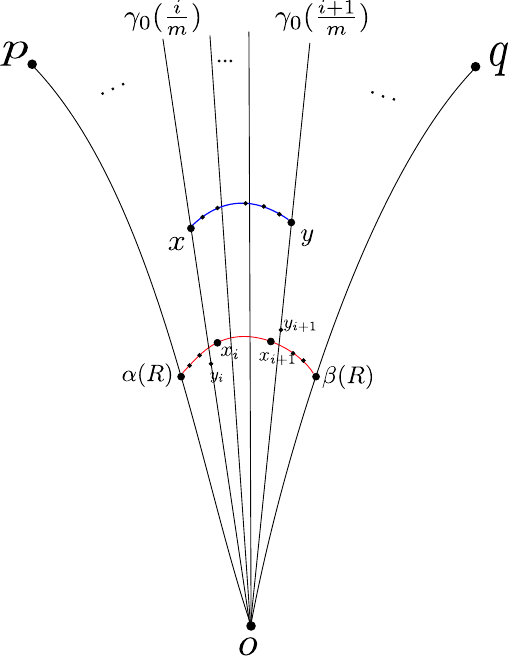}
    \caption{
    At each step, two points on the previously chosen rays are connected by an edge path on the spheres in $\widehat{X}$, provided by Proposition~\ref{spheres are connected}. The path in red is the edge path in the base step and the path in blue is an edge path chosen in the next step.
    }
    \label{MainThmPic}
\end{figure}

\smallskip
\noindent\textbf{The base step.}
Apply Proposition~\ref{spheres are connected} to $\alpha(R)$ and $\beta(R)$.  We obtain an edge path
\[
 \alpha(R)=x_0,x_1,\ldots,x_m=\beta(R)
\]
in $X^{(1)}$ such that
\[
 x_i\in S_R(o)\cup S_{R+1}(o)
 \quad\text{for every }i,
\]
and $d_{\widehat{X}}(x_i,x_j) \leq D_0$ for each $0\leq i,j\leq m$. If $m\geq2$, let $\rho_0=\alpha$ and $\rho_m=\beta$.  For each $0<i<m$, Proposition~\ref{exist close geodesic ray} provides a geodesic ray $\rho_i$ based at $o$ that passes within $C$ of $x_i$.  Define
\[
 I_0=\{i/m:0\leq i\leq m\}, \qquad \gamma_0(i/m)=\rho_i(\infty).
\]
If $m\leq1$, then we just let $I_0=\{0,1/2,1\}$, and 
\[
    \rho_0=\rho_{1/2}=\alpha,
    \qquad
    \rho_1=\beta,
    \qquad
    \gamma_0(a)=\rho_a(\infty)\quad\text{for each }a\in I_0.
\]

Let
\[
 T_0=\max\{0,R-C\}.
\]
We claim that, whenever $a,a'\in I_0$ are consecutive, their assigned rays satisfy
\begin{equation}
\label{ray subdivision invariant}
 d_{\widehat{X}}\bigl(\rho_a(T_0),\rho_{a'}(T_0)\bigr)\leq K.
\end{equation}
To see this when $R\geq C$, choose on each ray a point, say $y_i$, within $C$ of its corresponding $x_i$.  Then 
\[
R - C \leq d_{\widehat{X}}(o,y_i) \leq R + C +1.
\]
Since $x_i$ and $x_{i+1}$ are adjacent,
\[
d_{\widehat{X}}(\rho_a(T_0),\rho_{a'}(T_0)) \leq 2((2C+1)+C) + 1 = K,
\]
for each consecutive pair $a,a'\in I_0$. When $R<C$, both rays begin at $o$ and the estimate at time $T_0=0$ is immediate.

The points $x_i$ have diameter at most $D_0$. Hence, for $0\leq i,j\leq m$, the distance between $y_i$ and $y_j$ is at most $D_0+2C$, and both points lie at distance at least $R-C$ from $o$.  It follows from \eqref{close deep points give close boundary points} that, for some constant $A_1$ (depending only on $D_0$ and $C$),
\begin{equation}
\label{initial boundary chain diameter}
 \diam\bigl(\gamma_0(I_0)\bigr)\leq A_1b^{-R}.
\end{equation}

\smallskip
\noindent\textbf{The inductive step.}
Suppose that $I_r$, $\gamma_r$, and the rays $\rho_a$, for $a\in I_r$, are defined so that \eqref{ray subdivision invariant} holds for every pair of consecutive points of $I_r$ at time
\[
 T_r=T_0+r.
\]

We will construct the finite set $I_{r+1}$ and the corresponding map $\gamma_{r+1}$. Let $(s,t)$ be a component of $[0,1]\setminus I_r$.  Its endpoints belong to $I_r$.  Let
\[
 L_r=T_r+C+1,
 \qquad
 x=\rho_s(L_r),
 \qquad
 y=\rho_t(L_r).
\]
Hence,
\[
 d_{\widehat{X}}(x,y)
 \leq K+2C+2
 =E_1.
\]
We apply Proposition~\ref{spheres are connected} again to obtain an edge path
\[
 x=z_0,z_1,\ldots,z_n=y
\]
in $S_{L_r}(o)\cup S_{L_r+1}(o)$ whose diameter in $\widehat{X}$ is at most $D_1$.

For each $0<i<n$, choose a geodesic ray based at $o$ that passes within $C$ of $z_i$.  If $n \geq2$, subdivide $[s,t]$ into $n$ subintervals equally and assign these rays to the new subdivision points in the order determined by the path.  If $n\leq1$, insert the midpoint of $[s,t]$ and assign to it the ray $\rho_s$.  

We verify that the induction invariant persists.  A point on a newly chosen ray lying within $C$ of some $z_i$ has distance from $o$ between
\[
 T_r+1
 \quad\text{and}\quad
 T_r+2C+2.
\]
For consecutive vertices $z_i,z_{i+1}$, the corresponding nearby ray points are at distance at most $2C+1$. Note that $T_{r+1} = T_r +1$. Then we have
\[
 d_{\widehat{X}}
 \bigl(\rho_a(T_{r+1}),\rho_{a'}(T_{r+1})\bigr)
 \leq 2(2C+1)+(2C+1)
 =K,
\]
where $a,a'$ are consecutive points of a component of $[0,1]\setminus I_r$. When $n \leq 1$, the inserted ray is $\rho_s$, and the estimate follows directly from $d_{\widehat X}(x,y)\leq1$ at time $L_r$.  This proves \eqref{ray subdivision invariant} at step $r+1$.

The points $z_i$ have diameter at most $D_1$, and their nearby ray points lie at distance at least $T_{r+1}$ from $o$.  Consequently, \eqref{close deep points give close boundary points} gives a constant $A_2$, independent of $r,s,t$, such that 
\begin{equation}
\label{boundary subdivision diameter}
 \diam\bigl(\{\gamma_{r+1}(u)\mid u\in I_{r+1}\cap [s,t]\}\bigr)
 \leq A_2b^{-T_{r+1}}.
\end{equation}

\smallskip
\noindent\textbf{Extension of the map.}
Let
\[
 I=\bigcup_{r\geq0}I_r.
\]
Every component of $[0,1]\setminus I_r$ is divided into at least two subintervals at the next step. Thus the largest length of a complementary interval converges to 0; hence $I$ is dense in $[0,1]$.  We define the map
\[
 \gamma\colon I\longrightarrow\partial\widehat{X}
\]
by $\gamma(a) = \gamma_r(a)$ if $a\in I_r$.

For a component $(s,t)$ of $[0,1]\setminus I_r$, by applying \eqref{boundary subdivision diameter} iteratively,  we have
\[
 \diam\bigl(\{\gamma(u) \mid u \in I\cap [s,t]\}\bigr)
 \leq 2\sum_{k\geq r+1}A_2b^{-T_k}.
\]
Denote the series on the right-hand side by $Q_r$.  Then $Q_r\to 0$ as $r\to\infty$ since it is a geometric series with common ratio less than 1.

We claim that $\gamma$ is uniformly continuous.

Given $\varepsilon>0$, first choose $r$ so that $Q_r<\varepsilon/2$, and then choose $\delta>0$ so that every interval of length $\delta$ contains at most one point of $I_r$. This is possible since $I_r$ is a finite subset of $[0,1]$.  If $v,w\in I$ and $|v-w|<\delta$, then they lie either in the closure of one component of $[0,1]\setminus I_r$, or in the closures of two components having a common endpoint in $I_r$.  In the first case their images in $\partial \widehat{X}$ are at distance at most $Q_r$. In the second case, their images are at distance at most $2Q_r$.  Thus $d_\infty(\gamma(v),\gamma(w))<\varepsilon$.
\smallskip

Moreover,
\begin{align*}
 \diam\bigl(\gamma(I)\bigr)
 &\leq A_1b^{-R}
      +2\sum_{r\geq1}A_2b^{-(T_0+r)} \\
 &\leq A_3b^{-R}
\end{align*}
for a uniform constant $A_3$ (recall that $T_0=\max\{0,R-C\}$).

Note that the Gromov boundary $(\partial \widehat{X}, d_{\infty})$ is a complete metric space.  Since $I$ is dense in $[0,1]$, the uniformly continuous map $\gamma$  extends to a unique uniformly continuous map
\[
 \overline\gamma\colon[0,1]\longrightarrow\partial\widehat{X}.
\]
Its image is a compact connected set containing $p$ and $q$.  Finally, \eqref{visual distance at branching time} gives
\[
 \diam\bigl(\overline\gamma([0,1])\bigr)
 \leq A_3b^{-R}
 \leq A_0A_3d_\infty(p,q),
\]
where the constant $A_0A_3$ is independent of $p$ and $q$.  Therefore, $\partial\widehat{X}$ is linearly connected.
\end{proof}

%% file: sections/group-splitting.tex
As an application of Theorem~\ref{main theorem}, in this section, we relate the connectivity of $\partial \widehat{X}_\Gamma$ to a group splitting of $W_\Gamma$. The strategy is the following. The equivalence between the connectivity of $\partial \widehat{X}$ and a graph-theoretic condition on $\Gamma$ is established in Theorem~\ref{main theorem}. The equivalence between the algebraic splitting of $W_\Gamma$ and the same graph-theoretic condition will be obtained by using the visual decomposition theorem from \cite{mihalik2009visual}.

We say a group $G$ \emph{splits} over a subgroup $H\leq G$ if $G$ splits over $H$ as a free product with amalgamation or an HNN extension, that is, there are proper subgroups $A$ and $B$ of $G$ with $A \neq G \neq B$ such that 
\[
G \cong A\ast_H B \qquad \text{or} \qquad G\cong A\ast_H.
\]

We note that every group which splits as an HNN extension has a surjective homomorphism to $\Z$ by sending the stable element to 1. But there is no surjective homomorphism from a Coxeter group to $\Z$, since the abelianization of a Coxeter group is the finite abelian group $(\Z/2\Z)^{|V(\Gamma)|}$. Therefore, no Coxeter group splits as an HNN extension.

The main ingredient of this section is the following proposition where we show that a group splitting over a subgroup of an unbounded product subgroup forces one of the two graph-theoretic obstructions in Theorem~\ref{main theorem}.  The visual decomposition theorem of Mihalik--Tschantz provides the bridge from the arbitrary splitting to a visual splitting on the defining graph.

\begin{prop}
\label{separating subgraph}
If a one-ended right-angled Coxeter group $W_\Gamma$ splits over a subgroup of an unbounded product subgroup, then either $\Gamma$ has a separating nontrivial join or $\lk(v)$ is contained in a nontrivial join for some $v\in V(\Gamma)$.
\end{prop}

\begin{proof}
Suppose $C$ is a subgroup of an unbounded product subgroup $H = H_0\times H_1$ of $W_\Gamma$ such that both factors $H_0$ and $H_1$ are infinite and  $W_\Gamma$ splits over $C$ as an amalgamation $W_\Gamma \cong A\ast_C B$, where $A$ and $B$ are subgroups of $W_\Gamma$. Since $W_\Gamma$ is one-ended, Stallings' theorem implies that $C$ is infinite \cite{MR415622}.

\noindent\emph{Convert to a visual decomposition.}
The visual decomposition theorem \cite[Theorem~1 and Remark~1]{mihalik2009visual} gives a finite reduced visual graph-of-groups decomposition $\Psi$ of $W_\Gamma$ with the following properties. Every vertex group lies in a conjugate of $A$ or $B$, and every edge group lies in a conjugate of $C$. The underlying graph of $\Psi$ is a tree, since otherwise $W_\Gamma$ would be an HNN extension.  Moreover, there is an edge $e$ such that neither of the two components $\Psi_1$ and $\Psi_2$ of $\Psi\setminus e$ contains all the standard generators among its vertex groups.

For $i\in\{1,2\}$, let
\[
    S_i=\{s\in V(\Gamma):s\text{ belongs to a vertex group over }\Psi_i\},
\]
and let $\Gamma_i$ be the induced subgraph of $\Gamma$ spanned by $S_i$. For each standard generator $s\in V(\Gamma)$, the vertices and edges of $\Psi$ whose groups contain $s$ form a nonempty subtree by \cite[Lemma~4]{mihalik2009visual}.  It follows that
\[
    S_1\cup S_2=V(\Gamma).
\]
A generator belongs to both $S_1$ and $S_2$ precisely when it belongs to the special edge group $\Psi(e)$.  Let $\Lambda$ be the induced subgraph spanned by $S_1\cap S_2$.  Then
\[
    \Psi(e)=W_\Lambda,
\]
and neither $S_1\setminus V(\Lambda)$ nor $S_2\setminus V(\Lambda)$ is empty.

Mihalik--Tschantz's separator property \cite[Corollary~5]{mihalik2009visual} implies that $\Lambda$ separates $\Gamma$. In particular, there is no edge between $\Gamma_1\setminus\Lambda$ and $\Gamma_2\setminus\Lambda$, and $\Gamma=\Gamma_1\cup\Gamma_2$.

Collapsing every edge of $\Psi$ other than $e$ gives a nontrivial splitting of $W_\Gamma$ over $W_\Lambda$.  Since $W_\Gamma$ is one-ended, $W_\Lambda$ cannot be finite.  Thus $\Lambda$ is not a clique.

\smallskip
\noindent\emph{The visual edge group is contained in a standard product region.} We claim that $\Lambda$ is contained in a nontrivial join. Suppose not. Lemma~\ref{existence of contracting element} provides a contracting element $g\in W_{\Lambda}$. Then a conjugate 
\[
h = aga^{-1} \in C  \leq H
\]
would be a contracting element in the unbounded product subgroup $H$. Then the centralizer of a finite power of $h$ would contain a $\Z^2$ subgroup.  By~\cite[Theorem 5.4]{bestvina2009characterization}, an element $h\in W_\Gamma$ is contracting if and only if it is rank-one. A finite power of a rank-one element is still rank-one. And by \cite[Corollary 4.6(ii)]{caprace2010rank}, the centralizer of a rank-one element is virtually cyclic. Contradiction.

Hence $\Lambda \subset J$ for some nontrivial join subgraph $J$ of $\Gamma$.

\smallskip
\noindent\emph{Graph-theoretic conditions.} If $J$ separates $\Gamma$, then we are done.  Suppose that $J$ does not separate $\Gamma$.  Then $\Gamma\setminus J$ is empty or connected.  If $\Gamma \setminus J$ is empty, then $\Gamma = J$ is a nontrivial join and for each vertex $v\in V(\Gamma)$, the link $\lk(v)$ is contained in the nontrivial join $\Gamma = J$. 

Now suppose $\Gamma \setminus J$ is nonempty and connected.  The subgraph $\Lambda$ separates $\Gamma$, so there is a component $U$ of $\Gamma\setminus\Lambda$ contained in $J$.  In particular,
\[
    U\subset J\setminus\Lambda.
\]
Choose a vertex $u\in V(U)$.  Every neighbor of $u$ belongs either to $U$ or to $\Lambda$, because $U$ is a component of $\Gamma\setminus\Lambda$. Consequently,
\[
    \lk(u)\subset U\cup\Lambda\subset J.
\]
\end{proof}

We can now characterize these splittings through the connectivity of $\partial\widehat{X}$.

\begin{thm}
\label{splittings are topological}
    Let $W_{\Gamma}$ be a one-ended right-angled Coxeter group. The following are equivalent:
    \begin{enumerate}
        \item $W_\Gamma$ splits over a subgroup of an unbounded product subgroup;
        \item $\Gamma$ has a separating nontrivial join, or $\lk(v)$ is
        contained in a nontrivial join for some $v\in V(\Gamma)$;
        \item $\partial\widehat{X}_\Gamma$ is empty or disconnected.
    \end{enumerate}
\end{thm}

\begin{proof}
    For one-ended RACGs, condition~\textup{(a)} of Theorem~\ref{main theorem} holds. Hence, the equivalence of \textup{(2)} and \textup{(3)} is given  by Theorem~\ref{main theorem}. \textup{(1)} implies \textup{(2)} by Proposition~\ref{separating subgraph}. It remains to prove that \textup{(2)} implies \textup{(1)}. A separating nontrivial join gives the required visual splitting immediately.

    Now suppose $\lk(v)$ is contained in a nontrivial join. If $\Gamma\setminus\st(v)$ is nonempty, then $\lk(v)$ separates $v$ from that set, and $W_\Gamma$ splits over $W_{\lk(v)}$. Otherwise, $V(\Gamma)=\{v\}\cup\lk(v)$ and $\Gamma$ is a nontrivial join. Because $W_\Gamma$ is one-ended, $\Gamma$ is not a clique, so $\lk(v)$ contains non-adjacent vertices $x,y$. Let $\Delta=\Gamma\setminus\{x,y\}$. Then $\Delta$ lies in a nontrivial join and separates $\Gamma$, hence
    \[
    W_\Gamma\cong W_{\Delta\cup\{x\}}*_{W_\Delta}W_{\Delta\cup\{y\}}.
    \]
\end{proof}